\documentclass[draft,12pt,leqno]{amsart}

 \usepackage{amsmath,amssymb,amscd,amsthm}
 \usepackage{latexsym}

 \newtheorem{theorem}{Theorem}
 \newtheorem{lemma}{Lemma}
 \newtheorem{proposition}{Proposition}

 \newtheorem{corollary}{Corollary}
 \newtheorem{claim}{Claim}
 
\newcommand{\q}{\quad}
\newcommand{\qq}{\quad\quad}
\newcommand{\qqq}{\quad\quad\quad}

\newcommand{\norm}[2]{{\left\| #1 \right\|}_{#2}}
\newcommand{\f}[2]{\frac{#1}{#2}}
\newcommand{\dpr}[2]{\langle #1,#2 \rangle}

\newcommand{\pp}[2]{\f{\p #1}{\p #2}}

\newcommand{\al}{\alpha}

\newcommand{\ga}{\gamma}
\newcommand{\Ga}{\Gamma}
\newcommand{\de}{\delta}

\newcommand{\De}{\Delta}
\newcommand{\ve}{\varepsilon}
\newcommand{\Ve}{\Xi}

\newcommand{\ka}{\kappa}
\newcommand{\la}{\lambda}

\newcommand{\La}{\Lambda}
\newcommand{\si}{\sigma}

\newcommand{\vp}{\varphi}

\newcommand{\Om}{\Omega}

\newcommand{\rn}{{\mathbf R}^n}

\newcommand{\rone}{\mathbf R^1}

\newcommand{\rthree}{\mathbf R^3}

\newcommand{\sn}{\mathbf S^{n-1}}

\newcommand{\cs}{\mathcal S}

\newcommand{\cz}{\mathcal Z}

\newcommand{\cc}{\mathcal C}

\newcommand{\cl}{\mathcal L}

\newcommand{\intl}{\int\limits}

\newcommand{\suml}{\sum\limits}
\newcommand{\supl}{\sup\limits}

\newcommand{\p}{\partial}

\newcommand{\beq}{\begin{equation}}
\newcommand{\eeq}{\end{equation}}
\newcommand{\beqna}{\begin{eqnarray*}}
\newcommand{\eeqna}{\end{eqnarray*}}
\newcommand{\beqn}{\begin{equation*}}
\newcommand{\eeqn}{\end{equation*}}
\newcommand{\bp}{\begin{proof}}
\newcommand{\ep}{\end{proof}}
\newcommand{\bprop}{\begin{proposition}}
\newcommand{\eprop}{\end{proposition}}
\newcommand{\bt}{\begin{theorem}}
\newcommand{\et}{\end{theorem}}
\newcommand{\bex}{\begin{Example}}
\newcommand{\eex}{\end{Example}}
\newcommand{\bc}{\begin{corollary}}
\newcommand{\ec}{\end{corollary}}
\newcommand{\bcl}{\begin{claim}}
\newcommand{\ecl}{\end{claim}}
\newcommand{\bl}{\begin{lemma}}
\newcommand{\el}{\end{lemma}}

\begin{document}

\title
[Magnetic Schr\"odinger equation]
{Strichartz estimates for the magnetic Schr\"odinger equation}

\author{Atanas Stefanov}

\address{Atanas Stefanov\\
Department of Mathematics \\
University of Kansas\\
1460 Jayhawk Blvd\\ Lawrence, KS 66045--7523}

\thanks{Supported in part by  NSF-DMS 0300511}
\date{\today}


\begin{abstract}
We prove global,  scale invariant 
Strichartz estimates for the  linear 
magnetic Schr\"odinger equation with small time  dependent 
magnetic field. This is done  by  constructing an 
 appropriate parametrix. 
As an  application, we show a global regularity type result for 
 Schr\"odinger maps in  dimensions $n\geq 6$.  
\end{abstract}

\maketitle

\section{Introduction}
\label{sec:1}
In this paper, we investigate the global behavior of  
certain quantum dynamical systems in the presence of magnetic field. 
To describe the relevant equations, 
introduce the magnetic Laplace  operator 
$$
\De_{\vec{A}}= \suml_{j=1}^n (i \p_j+A_j)^2
$$
The magnetic Schr\"odinger  equation is
\begin{equation}
\label{eq:nac}
\left|
\begin{array}{l}
u_t - i (\De_{\vec{A}}+V) u=F,\\
u(0,x)=f(x)\in L^2(\rn)
\end{array}
\right.
\end{equation}
In the physically 
important case of a real valued $\vec{A}$ and $V$, one has 
conservation of charge, $\norm{u(t,\cdot)}{L^2}=\norm{u(0,\cdot)}{L^2}$. 
More generally,  by  a result of Leinfelder and Simader, \cite{Leinfelder} if
$A\in L^4_{loc}(\rn)$ and $div(A)\in L^2_{loc}(\rn)$,   $V$ is 
relatively bounded with bound less than one
 with respect to $\De$ , one has 
that the operator $\De_{\vec{A}}+V$ is essentially 
self-adjoint on $C_0^\infty(\rn)$. 
In particular, the  spectrum is real and 
one can define functional calculus. 

In this paper, we shall be concerned mainly with the case of 
{\it time dependent vector 
potentials } $\vec{A}$, which are small and real-valued. 
This  is 
dictated by certain  partial differential equations,   appearing  naturally in 
geometry and physics. More specifically, we have 
in mind
 the Schr\"odinger map equation in Hodge 
gauge \cite{Shatah}, \cite{Stefanov}, the Ishimori system, 
\cite{KPV00}, \cite{nahmod}, the Maxwell-Schr\"odinger system, 
\cite{Ginibre2},\cite{nakamura},  \cite{tsutsumi}  and  
several  other models, related to 
the Landau-Lifshitz theory of electromagnetism.

In the case of magnetic-free field ($\vec{A}=0$), 
great progress has been made to address the question for global/local 
existence and uniqueness 
for solutions of 
\eqref{eq:nac}, \cite{Cazenave}, \cite{Sulem}. In particular, 
when $V$ is small and $n\geq 3$, 
one can use the standard Strichartz estimates to 
show by a perturbation 
argument that the corresponding equation has an unique global 
solution 
under reasonable assumptions on 
the right hand-side and the data $f$. In the same spirit, 
one can obtain local well-posedness results for large $V$. 

In the magnetics-free case, the Strichartz estimates 
are well-known  and play a fundamental role in proving 
the existence and uniqueness results alluded to  above. \\
Introduce
$$
\norm{u}{L^q_t L^r}:= \left(\intl_0^\infty \left(\intl_{\rn} |u(t,x)|^r dx\right)^{q/r} dt\right)^{1/q}.
$$
for every pair $q,r\geq 1$ and similarly  the mixed Lebesgue spaces 
$L^q_t L^{r_1}_{x_1}\ldots L^{r_n}_{x_n}$. 

We say that  a pair of indices $(q,r)$  is {\it Strichartz admissible} if 
$2\leq q,r \leq \infty$, $2/q+n/r=n/2$ and 
$(q,r,n)\neq (2,\infty, 2)$.  Then, by a classical result of 
Strichartz, \cite{Strichartz}, later improved by 
Ginibre-Velo, \cite{Ginibre} and finally Keel-Tao, \cite{KeelTao}, we have 
\begin{eqnarray}
\label{eq:200}
& & \norm{e^{ i t \De} f}{L^qL^r}\leq C\norm{f}{L^2}\\
& &
\label{eq:202} \|\int_0^t e^{ i (t-s) \De} 
F(s,\cdot) ds\|_{L^qL^r}\leq C\norm{F}{L^{\tilde{q}'}L^{\tilde{r}'}},
\end{eqnarray}
where $(\tilde{q}, \tilde{r})$ is another Strichartz admissible pair and 
$q'=q/(q-1)$. 

Clearly \eqref{eq:200} and  \eqref{eq:202} are equivalent
to $
\norm{u}{L^qL^r}\leq C \norm{f}{L^2}+ \norm{F}{L^{\tilde{q}'}L^{\tilde{r}'}},
$
whenever $u$ is a solution to the free Schr\"odinger equation with initial data $f$ and forcing term $F$. 
 
Another equivalent formulation is that there exists a constant $C$, so that 
for all test functions 
$\psi:$ 
$$
 \norm{\psi}{L^qL^r}\leq C(\norm{\psi(0,\cdot)}{L^2}+
\norm{(\p_t\pm i\De)\psi}{L^{\tilde{q}'}L^{\tilde{r}'}}).
$$
 In the sequel, we will  make 
extensive use of all these 
points of view.

For  the case of small (but non-zero) potential $\vec{A}$, we have  
$\De_{\vec{A}} =-\De+ 2 i \vec{A}\cdot \nabla +
(i div(\vec{A})+\sum_j A_j^2)\cdot=-\De+2 i \vec{A}\cdot\nabla  + 
\textup{small potential}$, we can effectively treat the 
 the magnetic Schr\"odinger equation in the form
\begin{equation}
\label{eq:1}
\left|
\begin{array}{l}
\p_t u - i \De u+\vec{A}(t,x) \cdot \nabla u = F \q (t,x)\in 
{\mathbf R}^+\times \rn\\
u(0,x)=f(x)
\end{array}
\right.
\end{equation}
where the terms in the form $(i div(\vec{A})+\suml_j A_j^2)u$ 
are subsumed in the right hand side.

Next, we explain the relevance of the magnetic Strichartz estimates to the 
Cauchy problem for  
\begin{equation}
\label{eq:477}
\p_t u- i \De u+\vec{A}(u)\cdot  \nabla u= F(u).
\end{equation}
If $\vec{A}=0$, we can clearly use  \eqref{eq:200}, \eqref{eq:202} to 
set up an iteration scheme for the {\it semilinear problem }\eqref{eq:477} 
in a ball $B=B(0,R)$ in the  ``Strichartz space''
$\cap_{(q,r)}  L^q(0,T) L^r$ to solve for arbitrary $L^2$ 
data, provided one can show $\norm{F(u)}{L^1(0,T) L^2}\leq 
T^\de M(\norm{u}{L^q(0,T)L^r})$ for some bounded 
function $M$.  Choosing $R\sim \norm{f}{L^2}$ and 
$T: T^\de M(R)<< R$ makes such scheme successful to show that the solution 
exists for some time $T=T(\norm{f}{L^2})$. There are of course, 
remaining unresolved by this approach, including 
globality of such solutions\footnote{Note that such approach usually 
guarantees the global existence of small solutions.}, 
smoothness of the solution  etc. 

Clearly, to study \eqref{eq:477} with $\vec{A}\neq 0$, one cannot use 
the standard Strichartz estimates, by the obvious derivative loss. 
One of the goals of this paper is to derive global scale invariant Strichartz 
estimates under appropriate smallness assumptions on the vector 
potential $\vec{A}$. The pioneering work of Barcelo-Ruiz-Vega, 
\cite{barcelo} has addressed some of these issues\footnote{Strictly speaking, 
the 
results in \cite{barcelo} yield scale invariant 
 smoothing estimates, but standard methods allows one to derive 
Strichartz estimates from the results there.}, but was   
restricted to (essentially) radial vector potentials $\vec{A}$. 
To the best of our knowledge, the results in Theorem \ref{theo:pisa} below 
 are the first 
global estimates of such type for Schr\"odinger equations, 
that work for general non-radial potentials $\vec{A}$. 

Let us explain the general scheme for applying such Strichartz 
estimates to concrete quasilinear PDE's. 
Suppose, we have such estimates for the {\it linear gradient 
Schr\"odinger equation} \eqref{eq:1}, provided $\norm{\vec{A}}{Y_T}\leq \ve$, 
 for some concrete 
Banach space $Y_T$ appearing in Section \ref{sec:09}. 
 We apply the magnetic 
Strichartz estimates  to the nonlinear 
equations of the type \eqref{eq:477} as follows. For small 
initial data $f$, run an iteration scheme in the ball 
$B_X(0,R)$ in an  appropriate 
Strichartz space\footnote{Usually one solves the equations 
\eqref{eq:477} for data $f$ in some smooth Sobolev space $H^s$  and 
very often 
in Besov variants of the Strichartz space.}  $X$, 
see Section \ref{sec:11} for precise definitions. 
This is 
possible if 
\begin{itemize}
\item one can ensure {\it a priori} the smallness condition 
$\norm{\vec{A}(u)}{Y_T}\leq \ve$ for {\it all functions} $u$, that are solutions to \eqref{eq:1} 
satisfying $\norm{u}{X}<R$ and for all times $T\leq T_0=T_0(R)\leq \infty$. 
\item $\norm{F(u)}{L^1(0,T) L^2}\leq R/2$, whenever $u$ is a function with 
$\norm{u}{X}<R$ and for all times $T\leq T_0=T_0(R)\leq \infty$. 
\end{itemize}

\subsection{Strichartz estimates for the magnetic Schr\"odinger operator}
\label{sec:09}
The existence and uniqueness problem for \eqref{eq:1} has been studied 
extensively by 
many authors in the mathematics literature. 
We should first point to the pioneering work of Doi, \cite{Doi1}, \cite{Doi2}, 
who has devised a 
method to obtain solutions via energy estimates. The approach then relies 
on cleverly exploiting 
the properties of pseudodiferential operators of order zero  
to obtain {\it a priori} control of 
$\norm{u(t,\cdot)}{L^2}$ in terms of 
$\norm{f}{L^2}$ and $\norm{F}{L^1_T L^2}$. \\
 We also mention the far 
reaching generalization of 
Doi's results, due to Kenig-Ponce-Vega, \cite{KPV98}\footnote{In fact, in the 
proof of the Strichartz 
estimates for \eqref{eq:1}, we shall need a particular 
local existence result from \cite{KPV98}, see Section \ref{sec:4}.}. 
The authors have considered more general equations and 
were  able to
derive {\it a priori} estimates for the $L^2$ norms of the solution as well 
as the 
validity of 
a local smoothing effect, phenomenon well-known for the potential free case. 

Note that \eqref{eq:1} has the important 
scaling invariance 
$u\to u^\la(t,x)=u(\la^2 t, \la x)$, 
$A\to A^\la(t,x)= \la A(\la^2 t, \la x)$, 
$F\to F^\la(t,x)=\la^2 F(\la^2 t, \la x)$.  That is, 
whenever $(u,A, F)$ satisfy \eqref{eq:1}, so does 
$(u^\la, A^\la, F^\la)$ with initial data 
$f^\la(x)=f(\la x)$. \\
We describe the space $Y$ of  vector potentials $\vec{A}$, so 
that the corresponding magnetic \\ Schr\"odinger operator satisfies the 
Strichartz estimates. Let $SU(\rn)$ be the special unitary group acting on $\rn$ and $x(t):\rone_+\to \rn$ be arbitrary measurable function. Define 
\begin{eqnarray*}
& &\norm{\vec{A}}{Y_0} := \norm{\nabla \vec{A}}{L^1L^\infty}+
\norm{\vec{A}}{L^2L^\infty} +
(\suml_l 
2^{2l(1+h)}\norm{A_l}{L^1 L^{n/h}}^2)^{1/2}\\
& &\norm{\vec{A}}{Y_1}:=  \suml_k 2^{k(n-1)} \norm{A_k}{L^\infty_t L^1_x}\\
& &\norm{\vec{A}}{Y_2}:= \suml_k 2^{k(n-1)/2}\supl_{U\in SU(n)} \supl_x 
\| A_{k}(t,x+ Uz)\|_{L^\infty_t L^2_{z_2, \ldots, z_n}L^1_{z_1}}+ \\
& &
+ \suml_k 2^{k(n-5)/2}\supl_{U\in SU(n), x(t)}
\| (|\p^2 A_{k}|+|\p_t A_k|)(t,x(t)+Uz)\|_{L^\infty_t L^2_{z_2, 
\ldots, z_n}L^1_{z_1}})\\
& &
\norm{\vec{A}}{Y_3} := \suml_k 2^{k(n+3)/2}\supl_{U\in SU(n)} \supl_x 
\| A_{k}(t,x+ Uz)\|_{L^1_t L^2_{z_2, \ldots, z_n}L^1_{z_1}}+\\
& &+
\suml_k 2^{k(n-1)/2}\supl_{U\in SU(n), x(t)}
\| (|\p^2 A_{k}|+|\p_t A_k|)(t,x(t)+Uz)\|_{L^1_t L^2_{z_2, 
\ldots, z_n}L^1_{z_1}}).
\end{eqnarray*}
In the case $n\geq 4$, we can replace $Y_1$ by a 
bigger space $\tilde{Y_1}$ 
(with the smaller norm) 
$$
\norm{A}{\tilde{Y}_1}=\suml_k 2^{k(n-1)/p_0} 
\supl_{U\in SU(n), x}\|A_k(t,x+Uz)\|_{L^\infty_t L^{p_0}_{z_2, 
\ldots, z_n}L^1_{z_1}},
$$
for some $p_0: p_0<(n-1)/2$. 
\begin{theorem}
\label{theo:pisa}
Let $n\geq 2$. Then, there exists an $\ve=\ve(n)>0$, so that whenever 
$\vec{A}:\rone\times\rn\to \rn$ is a  
real-valued 
vector potential with $\|\vec{A}\|_{Y_0\cap Y_1\cap Y_2\cap Y_3}
\leq \ve$, (which can be relaxed to 
$\|\vec{A}\|_{Y_0\cap \tilde{Y}_1\cap Y_2\cap Y_3}
\leq \ve$, if $n\geq 4$)  \eqref{eq:1}
 has an unique global solution, whenever the initial data 
$f\in L^2(\rn)$ and the forcing term $F\in L^1L^2$. 
In addition, there exists a constant $C=C(n)$, 
so that the {\it a priori} estimate 
\begin{equation}
\label{eqn:478}
\supl_{q,r - Str.}\norm{u}{L^qL^r}\leq 
C(\norm{f}{L^2}+
\norm{F}{L^1L^2}).
\end{equation}
holds true\footnote{For the two dimensional case, the constant $C$ does 
also  depend 
on $(q,r)$. More specifically the constant blows up as $(q,r)\to 
(2,\infty)$.}. Moreover, for every $\psi\in\cs$ 
\begin{equation}
\label{eq:1009}
\supl_{q,r - Strichartz admissible} \norm{\psi}{L^qL^r}\leq 
C(\norm{\psi(0,\cdot)}{L^2}+\norm{(\p_t - i 
\De+\vec{A}\cdot \nabla) \psi}{L^1L^2}).
\end{equation}
One also 
 has the $l^2$ Besov space version 
\begin{equation}
\label{eq:1011}
(\suml_k \norm{\psi_k}{L^qL^r}^2)^{1/2}
\leq C(\norm{\psi(0,\cdot)}{L^2}+
\norm{(\p_t - i \De+\vec{A}\cdot \nabla) \psi}{L^1L^2}).
\end{equation}
\end{theorem}
{\bf Remark}
\begin{itemize}
\item Note  that for all $Y_j$, $j=0,1,2,3 $ we have that 
$\|\vec{A}^\la\|_{Y_j}=\|\vec{A}\|_{Y_j}$, that is the spaces are scale 
invariant with respect to the natural scaling $A\to A^\la$. 
\item In the case $n=1$,  the 
theorem holds as well. 
Our results are 
however far from optimal, as shown recently by Burq-Planchon, \cite{Burq}. 
It seems that in the one dimensional case one only 
needs to require 
$\sup_t \norm{A(t, \cdot)}{L^1(\rone)}<\infty$, if $A$ 
is a real-valued potential. 
\item If $A:\rone\times \rone\to \cc$ is complex valued, 
and satisfies  \\ $\sup_{t,x}|\int_{-\infty}^x A(t,y)dy|<\ve$ and 
$\|\int_{-\infty}^x (\p_t - i \p_y^2) A(t,y)dy\|_{L^1_t L^\infty_x}<\ve$, 
one has the results of Theorem \ref{theo:pisa}. This is 
shown in \cite{Stefanov2}, together with some applications and uniqueness 
issues,  consult \cite{Stefanov2} for more details. 
Note that by recent examples on  ill-posedness for 
derivative Schr\"odinger equations 
in $\rone$ (due to M. Christ, \cite{Crist20}), some 
smallness assumptions are necessary even for a local well-posedness. 
\end{itemize}
\subsection{Some corollaries}
We present some corollaries of Theorem \ref{theo:pisa}. 
Observe that by Bernstein inequality (Lemma \ref{le:yut}), one can bound 
\begin{eqnarray*}
& & \norm{A}{Y_0\cap Y_1\cap Y_2\cap Y_3}\leq C_n 
(\suml_k 2^{k(n-1)} \norm{A_k}{L^\infty L^1}+
2^{k(n-3)} \norm{\p_t A_k}{L^\infty L^1})+ \\
& &+C_n (\suml_k 2^{k(n+1)} \norm{A_k}{L^1 L^1}+
2^{k(n-1)} \norm{\p_t A_k}{L^1 L^1}).
\end{eqnarray*}
We thus have 
\begin{corollary}
\label{cor:90}
There exists a small positive $\ve>0$, so that whenever 
a real valued vector potential  $\vec{A}$ satisfies 
$$
\suml_k 2^{k(n-1)} \norm{A_k}{L^\infty L^1}+
2^{k(n-3)} \norm{\p_t A_k}{L^\infty L^1}+  2^{k(n+1)} \norm{A_k}{L^1 L^1}+
2^{k(n-1)} \norm{\p_t A_k}{L^1 L^1}\leq \ve,
$$
the conclusions \eqref{eqn:478}, \eqref{eq:1009} and \eqref{eq:1011} hold true. 
\end{corollary}
\noindent For the case of {\it time independent magnetic  potential} $\vec{A}$, we can formulate the following immediate corollary of Theorem \ref{theo:pisa}.
\begin{corollary}
\label{cor:1}
Let $n\geq 2$. Then there exists an 
$\ve=\ve(n)$, so that 
the magnetic Schr\"odinger equation \eqref{eq:1} has an 
unique global solution, provided $\vec{A}=\vec{A}(x)$ is real-valued 
vector function  and 
\begin{eqnarray*}
& &\suml_k 2^{k(n+1)} \norm{A_k}{L^1}
<\infty \\
& & \suml_k 2^{k(n-1)} \norm{A_k}{L^1_x}\leq \ve.
\end{eqnarray*}
Moreover, the solution satisfies 
$$
\supl_{q,r - Str.}\norm{u}{L^q(0,T) L^r}\leq 
C(T,n)(\norm{f}{L^2}+\norm{F}{L^1 L^2}).
$$
\end{corollary}
\noindent Corollary \ref{cor:1} follows easily by just applying Corollary  
\ref{cor:90} with a vector potential of the form $\tilde{A}(t,x):=
\vec{A}(x) \chi(t/\de)$ 
for some appropriate smooth cutoff function $\chi$ and a small $\de$. 
This will produce a solution in
 a small time interval, say $(0,\de/2)$, which is iterated and so on. 
 The smallness of the potential $\tilde{A}$ is achieved by the smallness of 
$\de$ (used to satisfy the requirements $A\in Y_0\cap Y_2\cap Y_3$)  
and by Sobolev embedding and the condition 
$\suml_k 2^{k(n-1)} \norm{A_k}{L^1_x}\leq \ve$ in the case $A\in Y_1$. 
Such a  proof provides an upper bound $C(t,n)\sim C^{T/c\ve}$, 
which is not  optimal in general. 
\subsection{Strichartz estimates with derivatives} 
The Strichartz estimates described in Theorem \ref{theo:pisa} 
can be of course  extended to control the norms of the 
solution $u$ in Besov type norms involving derivatives. One way to do that 
is  considering the Littlewood-Paley reduction of the equation to 
a fixed frequency $k$, applying the regular Strichartz estimates 
(either \eqref{eq:1009} or \eqref{eq:1011} with 
appropriate $p_1, p_2, q_1, q_2$), then multiplying by the 
corresponding power of $2^{ks}$ and square summing in $k$. The result is 
\begin{theorem}
\label{theo:3}
Let $n\geq 2$ and $\vec{A}$ satisfies 
the assumptions in Theorem \ref{theo:pisa}. Then there exists a 
constant $C=C(n)$ ($C=C(n,q,r)$, if $n=2$), so that for 
every $s>0$, initial data $f\in \dot{H}^s$ and 
forcing term $F\in L^1_t\dot{H}^s$, the global solution $u$ 
of \eqref{eq:1} satisfies 
\begin{equation}
\label{i:1}
(\suml_k 2^{2ks} \norm{u_k}{L^q_t L^r_x}^2)^{1/2}\leq C (\norm{f}{\dot{H}^s}+ 
\norm{F}{L^1_t\dot{H}^s})+C \norm{\nabla u}{L^2_t L^{n}_x}
(\suml_k 2^{2ks} \norm{A_k}{L^2_t L^{2n/(n-2)}_x}^2)^{1/2}. 
\end{equation}
for every Strichartz admissible pair $(q,r)$. 
\end{theorem}
\subsection{Strichartz estimates in 
 $L^2_t L^{2(n-1)/(n-3)}_{x_2, \ldots, x_n} 
L^2_{x_1}$}

We  present an extension of Theorem \ref{theo:pisa}, 
which allows us  to control a larger set of norms. 
\begin{proposition}
\label{prop:3}
Let $n\geq 4$ and $\vec{A}$ satisfies 
the assumptions in Theorem \ref{theo:pisa}. 
Then there exists a constant $C_n$, so that 
the solutions of \eqref{eq:1}, satisfy 
\begin{equation}
\label{e:7}
(\suml_k \supl_{U\in SU(\rn), x(t)} 
\norm{u_k(t,x(t)+U z)}{L^2_tL^{2(n-1)/(n-3)}_{z_2, \ldots, z_n} 
L^2_{z_1}}^2)^{1/2}
\leq C(\norm{f}{L^2}+
\norm{F}{L^1L^2}).
\end{equation}
For $n=3$, take any $(q,r):1/q+1/r=1/2$ and $q>2$. Then  there exists 
a constant $C_q$ (which may blow up as $q\to 2$), so that 
$$
(\suml_k \supl_{U, x(t)} 
\norm{u_k(t,x(t)+U z)}{L^q_t L^{r}_{z_2,z_3} 
L^2_{z_1}}^2)^{1/2}
\leq C_q(\norm{f}{L^2}+
\norm{F}{L^1L^2}).
$$
For $n=2$, 
$$
(\suml_k \supl_{U, x(t)} 
\norm{u_k(t,x(t)+Uz)}{L^4_t L^{\infty}_{z_2} 
L^2_{z_1}}^2)^{1/2}
\leq C(\norm{f}{L^2}+
\norm{F}{L^1L^2}).
$$
\end{proposition}
We also have a generalization of Proposition \ref{prop:3} to 
the setting of Theorem \ref{theo:3}, that is involving derivatives. 
Namely, say when $n\geq 4$, one has the {\it a priori} estimate 
\begin{eqnarray*}
& & 
(\suml_k \supl_{U\in SU(\rn), x(t)} 
2^{2ks}\norm{u_k(t,x(t)+U z)}{L^2_tL^{2(n-1)/(n-3)}_{z_2, \ldots, z_n} 
L^2_{z_1}}^2)^{1/2}
\leq C(\norm{f}{\dot{H}^s}+
\norm{F}{L^1\dot{H}^s})+ \\
& &+C \norm{\nabla u}{L^2_t L^{n}}
(\suml_k 2^{2ks} \norm{A_k}{L^2_t L^{2n/(n-2)}_x}^2)^{1/2}.
\end{eqnarray*}

Such results are needed to connect the ``solutions space'' 
with the space of admissible  vector potentials $Y$. More precisely, 
in the applications, we have $A=A(u)$, where 
the relation is usually in the form 
$A=Q(u, \bar{u})$, where $Q$ is a bilinear (or multilinear) 
form acting on the solution and its conjugate. For example, in the 
Schr\"odinger map case (see Section \ref{sec:smaps} below), 
we have schematically $A=|\nabla|^{-1} (u \bar{u})$ 
and for the Maxwell-Schr\"odinger system (see Section \ref{sec:MS} below)
we have $A=\square^{-1}(u \nabla u)$. 

Ignoring the derivatives (and the invertibility of $\square$ in the MS case)
for a second, 
we see that by H\"older's inequality, to have estimates of the form 
$$
\norm{A(u)}{Y_2\cap Y_3}\leq C \norm{u}{X}^2
$$
where $X$ is the solutions space, we must rely on mixed Lebesgue 
estimates like the one in Theorem \ref{theo:3}. Moreover, $X$ must be intersected with a space given by one of the norms involved  in the left hand side of \eqref{e:7}. 

\subsection{Applications to Schr\"odinger maps}
\label{sec:smaps}
In this section, we present a global regularity type result for the 
so-called Modified Schr\"odinger map system (MSM), which was derived in 
\cite{Stefanov}.  According to Theorem 2.2, \cite{Stefanov}, 
the Schr\"odinger map problem, with with target $\sn$, 
was shown to be equivalent (modulo a Lorentz type gauge transformation), 
to a {\it overdetermined} system of Schr\"odinger equations 
with attached consistency conditions.  

We will not discuss here, whether the MSM and 
 the Schr\"odinger map problem are equivalent, 
and how the (properties of the) solutions to one 
relate to the solutions of the other,  
with the acknowledgement that these are by no means unimportant 
or irrelevant issues. We will concentrate instead on the 
question of existence of solutions for MSM, which is 
mathematically more tractable. \\
Consider the MSM, which  takes the form 
\begin{eqnarray*}
\pp{u_j}{t} &=& i\De u_j-2\suml_{k=1}^n  a_k \pp{u_j}{x_k}-
\left(\suml_{k=1}^n  a_k^2\right)u_j+ \\
&+& 2(\suml_{k=1}^n \ Im (\overline{u_j}u_k)u_j) -i a_0u_j \qqq j=1,\ldots,n;
\end{eqnarray*}
where 
\begin{eqnarray*}
a_k &=& \suml_{l=1}^n \pp{\ka_{lk}}{x_l};\\
d\ \ka &=& 0;\\
\De \ka_{kj} &=& -4\ Im(u_k\overline{u_j}) \qq\begin{array}{l}
j=0,1,\ldots,n\\
k=1,\ldots,n,
\end{array}; \\
\De a_0 &=& -4 \suml_{j=1}^n\suml_{k=1}^n \left[\pp{}{x_k}\pp{}{x_j}
Re(u_k\overline{u_j})-\f{1}{2}\left(\pp{}{x_k}\right)^2u_j\overline{u_j}
\right].
\end{eqnarray*}
In short, we will consider 
the following system of Schr\"odinger equations
\begin{equation}
\label{eqnm:1}
\left|
\begin{array}{l}
\p_t u - i \De u+ A(u)\cdot \nabla u =N(u) \\
u(0,x)=f,
\end{array}\right. 
\end{equation}
where $A$ is a real-valued vector potential 
($A=(a_1, \ldots, a_n)$ in MSM), 
\begin{eqnarray*}
& & A(u)=\p^{-1} Q_1(u,\bar{u}) \\
& & a_0(u)= \p^2\De^{-1} Q_2(u, \bar{u})\\
& & N(u)= Q_3(u,u, \bar{u})+Q_4(A, A, u)+Q_5(a_0, u).
\end{eqnarray*}
Here, $Q_1, \ldots, Q_5$ are multilinear forms of their arguments, i.e. \\
$Q_j(u^1, \ldots, u^r)=\sum_j  c_{k_1, \ldots, k_r}^j u^{k_1}...
u^{k_r}$ for some constants $c$. 
We have also  used the notation $\p^s$ to denote a multiplier type operator, 
whose smooth symbol satisfies  $|s(\xi)|\sim |\xi|^s$. \\
All the results that we obtain for \eqref{eqnm:1} cover the MSM system, 
which is our main motivation. 
\begin{theorem}(Global regularity of MSM in high dimensions)
\label{theo:MSM} Let $n\geq 6$, $s_0=n/2-1$ and $s\geq (n+1)/2$. 
Then, there exists 
$\ve>0$, so that whenever $g\in H^s$, with 
$\norm{g}{\dot{H}^{s}}\leq 1$, 
the solution to \eqref{eqnm:1} with initial data $f=\ve g$ 
exists globally and satisfies 
\begin{eqnarray*}
& &\supl_{t} \norm{u(t, \cdot)}{H^s}\leq C \ve 
\end{eqnarray*}
for some constant $C$ depending only on the dimension and $s$. 

For $n=5$, there is an appropriate Besov spaces analogue. 
\end{theorem}

\subsection{Maxwell-Schr\"odinger system}
\label{sec:MS}
Another system of nonlinear PDE's, for which the Strichartz estimates of 
 Theorem \ref{theo:pisa} are applicable is the Maxwell-Schr\"odinger 
system. That is 
\begin{equation}
\label{eqnm:2}
\left|
\begin{array}{l}
i \p_t u + \De_A u= g(|u|^2) u \\
\square A ={\mathbf P} Im(\bar{u} \nabla_A u) ,
\end{array}\right. 
\end{equation}
where ${\mathbf P}$ is the 
Leray projection onto the divergence free vector fields and 
 $g(|u|^2)$ is either the Hartree 
interaction
$g(|u|^2):= \int G(x-y) |u|^2 (y)dy$ or\footnote{Here $G$ is the 
Green function in $n$ dimension.}  simply $g(|u|^2)=|u|^2$.

 This system was studied  by Tsutsumi,  \cite{tsutsumi}, where 
he  constructs the wave operator on 
a class of small scattered states. In particular, 
Tsutsumi showed global existence for a particular class of small data. 

Recently, Nakamura and Wada,  \cite{nakamura} have considered the 
MS system \eqref{eqnm:2}as well. They have obtained local well-posedness 
with data $u_0 \in H^{5/3}(\rthree)$  by using energy estimates approach. 
For related results and recent developments for \eqref{eqnm:2}, one might consult the recent work of Ginibre and Velo, \cite{Ginibre2}. 

We do not state any precise results here for the 
Maxwell-Schr\"odinger system, but it is clear that a variant of Proposition \ref{prop:3} for solutions of the linear wave equation (among other things) 
will be required. 
These issues will be addressed in a forthcomming paper, \cite{Stefanov3}.

A short outline of the paper is as follows. In Section \ref{sec:prelim}, 
we give some definitions from harmonic analysis as well as some
facts from the abstract Strichartz estimates theory 
due to Keel and Tao, \cite{KeelTao}. 
In Section \ref{sec:4}, 
we give some classical energy estimates and Littlewood-Paley reductions, 
which reduce the problem to the existence of parametrix 
construction.  
In Section \ref{sec:par}, 
we motivate and construct the parametrix and then we prove some of 
its main properties. 
In Section \ref{sec:98}, we describe an important 
angular decomposition for the phase of the parametrix and as a corollary 
we show the  crucial pointwise estimates, 
which are used throughout in the sequel. 
In Section \ref{sec:6}, we show that the parametrix satisfies 
 $L^q_t L^r_x$ estimates.
In Section \ref{sec:7}, we show that the 
parametrix almost satisfy the magnetic Schr\"odinger equation. 
In Section \ref{sec:8}, we prove the 
Strichartz estimates stated in Proposition \ref{prop:3}. In Section \ref{sec:11}, we show the global regularity for the modified Schr\"odinger maps. 
Some of the 
technical lemmas used in Sections \ref{sec:par} and \ref{sec:7} 
are formulated  and proved in  the Appendix. 
\section{Preliminaries}
\label{sec:prelim}
\subsection{Fourier transform and Littlewood-Paley projections}
\noindent 
Define the Fourier transform  by 
$$
\hat{f}(\xi)=\intl_{\rn} f(x)e^{-2\pi i x\cdot  \xi} dx
$$ 
and hence 
$$
 f(x)=\intl_{\rn} \hat{f}(\xi)e^{2\pi i x\cdot  \xi} d\xi
$$ 
Introduce a 
positive, smooth and  even function $\chi:\rone\to \rone$, supported in 
 $\{\xi:|\xi|\leq 2\}$ and $\chi(\xi)=1$ for all
$|\xi|\leq 1$.  
 Define $\vp(\xi)=\chi(\xi)-\chi(2\xi)$, 
which is supported in the annulus $1/2\leq |\xi|\leq 2$. Clearly 
$\sum_{k\in \cz} \vp(2^{-k} \xi)=1$ for all $\xi\neq 0$. \\
The $k^{th}$ Littlewood-Paley projection is defined as a multiplier 
type operator by \\
$\widehat{P_k f}(\xi)=\vp(2^{-k}\xi) \hat{f}(\xi)$. 
Note that the kernel of $P_k$ is integrable, smooth and real 
valued for every $k$. 
In particular, it commutes with differential operators. \\
Also of interest will be the properties 
of products under the action of $P_k$. We have that 
for any two (Schwartz ) 
functions $f, g$
\begin{eqnarray*}
P_k( f g) &=& \suml_{l\geq k-2} P_k ( f_l g _{l-2\leq \cdot \leq l+2}) +
\textup{symmetric term}+ \\
&+& P_k( f_{\leq k-4} g_{k-1\leq \cdot\leq  k+1}) +\textup{symmetric 
term}=\\
&=& f_{\leq k-4} g_k + [P_k, f _{\leq k-4}] g_{k-1\leq \cdot k+1} + 
\textup{symmetric terms} \\
&+& \suml_{l\geq k-2} P_k ( f_l g _{l-2\leq \cdot \leq l+2}) +
\textup{symmetric term}.
\end{eqnarray*}
\begin{lemma}(Bernstein inequality)
\label{le:yut}
Let $f$ be Fourier supported in a rectangle  $Q\subset \rn$. 
Then for every $1\leq p\leq q\leq \infty$, one has 
$$
\norm{f}{L^q}\leq C_n|Q|^{1/p-1/q}\norm{f}{L^p}
$$
If $Q={\xi: |\xi|\sim 2^k}$, one can extend to 
 mixed $L^pL^q$ spaces. Suppose $p_1>p_2\geq r$. 
Then 
$$
\norm{f}{L^{p_1}_{x_2, \ldots, x_n} L^{r}_{x_1}}\leq 
C_n 2^{k(n-1)(1/p_2-1/p_1)} 
\norm{f}{L^{p_2}_{x_2, \ldots, x_n} L^{r}_{x_1}}. 
$$
\end{lemma}
\begin{proof}
The first statement is standard. \\
For the second statement, it is equivalent to the boundedness of \\
$P_k:L^{p_2}_{x_2, \ldots, x_n} L^{q_2}_{x_1}\to L^{p_1}_{x_2, 
\ldots, x_n} L^{q_1}_{x_1}$ with bound $ C_n 2^{k(n-1)(1/p_2-1/p_1)}$. 
One can rescale to the case $k=0$, 
since these estimates are scale invariant. Since $P_0$ has integrable 
kernel, we have $P_0:L^q_{x_2, \ldots, x_n} L^r_{x_1}\to L^q_{x_2, 
\ldots, x_n} L^r_{x_1}$ for every $1\leq q,r\leq \infty$. 

On the other hand, an application of the Bernstein inequality 
in the $(n-1)$ variables $x_2, \ldots, x_n$
\begin{eqnarray*}
& & \norm{P_0 u}{L^\infty_{x_2, \ldots, x_n} L^1_{x_1}}\leq 
\norm{P_0 u}{ L^1_{x_1} L^\infty_{x_2, \ldots, x_n}}\leq 
\norm{P_0 u}{ L^1_{x_1} L^1_{x_2, \ldots, x_n}}\leq 
\norm{u}{  L^1_{x_2, \ldots, x_n}L^1_{x_1}}.
\end{eqnarray*} 
A complex interpolation between the last estimate and 
$P_0:L^qL^1\to L^qL^1$ (for every $1\leq q\leq \infty$) 
yields $P_0:L^pL^1\to L^{\tilde{p}}L^1$, whenever 
$1\leq \tilde{p}<p\leq\infty$. Interpolation between the last estimate and 
$p_0:L^\infty L^\infty\to L^\infty L^\infty$ yields 
$P_0: L^{rp}_{x_2, \ldots, x_n} L^r_{x_1}\to
 L^{r\tilde{p}}_{x_2, \ldots, x_n} L^r_{x_1} $ for every 
$r\leq \tilde{p}<p\leq\infty 1$. This is the second 
statement  of Lemma \ref{le:yut} for an appropriate choice of $p, \tilde{p}$. 
\end{proof}
We  also need the following technical lemma in the sequel
\begin{lemma}
\label{le:sequences}
Let $\{a_l\}$, $\{b_l\}$ are two sequences and $h>0$. Then 
$$
\left(\suml_k 2^{2 h  k} \left(\suml_{l\geq k-2} 2^{-h l} a_l b_l \right)^2\right)^{1/2}\leq 
C_h \norm{a}{l^\infty} \norm{b}{l^2}.
$$
\end{lemma}
\begin{proof}
Fix the sequence $\{a_l\}$ and consider the linear operator 
(mapping sequence into a sequence)
$$
(Tb)_k := 2^{\ve(k-l)} \suml_{l\geq k-2} a_l b_l
$$
We will show that $T:l^1\to l^1$ and $T:l^\infty\to l^\infty$. Indeed, 
\begin{eqnarray*}
& &\norm{Tb}{l^1}\leq \suml_l |a_l|| b_l| \suml_{k\leq l+2} 2^{\ve(k-l)}\lesssim \ve^{-1} 
\norm{a}{l^\infty}\norm{b}{l^1},\\
& & \norm{Tb}{l^\infty}\leq \supl_k \supl_l|a_l|  \supl_l|b_l|
\suml_{l\geq k-2} 2^{\ve(k-l)}\lesssim \ve^{-1} 
\norm{a}{l^\infty}\norm{b}{l^\infty}.
\end{eqnarray*}
It follows that for $1\leq p\leq \infty: T:l^p\to l^p$ with norm no bigger
 than $C_\ve\norm{a}{l^\infty}$, 
hence the statement of the lemma.
\end{proof}
\subsection{Keel-Tao theory}
It is well-known that decay and energy estimates 
imply Strichartz estimates in the context of various  dispersive 
equations. We would like to state an abstract result due to M. Keel 
and T. Tao, \cite{KeelTao}, which proved out to be 
very usefull in this context. Let us recall, that this method 
in conjunction with the Hausdorf-Young inequality was 
used by Ginibre and Velo in  their proof of the Strichartz 
estimates for the linear 
Schr\"odinger equation away from the endpoint. 

The abstract version
of Keel and Tao  has the (somewhat) 
unexpected consequence that the endpoint Strichartz estimate follows only 
from decay and energy estimates. 

Let $H$ be a Hilbert space, and $(X, dx)$ be a measure space and 
$U(t):H\to L^2(X)$ be a bounded operator. Suppose that $U$ satisfies 
\begin{eqnarray}
\label{e:100}
\norm{U(t) f}{L^2_x}\leq C\norm{f}{H} \\
\label{e:101} \norm{U(t) U(s)^*f}{L^\infty_x}\leq  C 
|t-s|^{-\si}\norm{f}{L^1}.
\end{eqnarray}
Suppose also that $(q,r)$ are $\si$ admissible, that is 
 $(q,r): q,r\geq 2; 1/q+\si/r=\si/2$ and 
$(q,r,\si)\neq (2,\infty,1)$. 
\begin{proposition}(Keel-Tao, \cite{KeelTao})
\label{prop:12} Let $(q,r)$ and $(\tilde{q}, \tilde{r})$ 
be both $\si$ admissible
 and $U(t)$ obeys \eqref{e:100} and \eqref{e:101}. Then 
\begin{eqnarray}
\label{e:102}
\norm{U(t)f}{L^q_t L^r_x}\leq C\norm{f}{L^2} \\
\label{e:103} \norm{\intl_0^t U(t)(U(s)^*) 
F(s,\cdot) ds}{L^qL^r}\leq C \norm{F}{L^{\tilde{q}'}_t L^{\tilde{r}'}_x}.
\end{eqnarray}
\end{proposition}
{\bf Remark} Note that by the Ginibre-Velo original argument 
(see also \cite{KeelTao}), the Strichartz estimates \eqref{e:102}, 
\eqref{e:103} follows only assuming the energy bound \eqref{e:100} 
and the ``modified decay bound''
\begin{equation}
\label{e:115}
\norm{U(t) U(s)^*f}{L^{p'}_x}\leq C_p |t-s|^{-(2/p-1)\si}\norm{f}{L^p},
\end{equation}
for any $\si: (1-2/p)\si >1$. Note that \eqref{e:115} follows by interpolation between 
\eqref{e:100} and \eqref{e:101}, and so it is in general easier 
 to establish. In fact, we use the modified decay bound, 
instead of the $L^1\to L^\infty$ decay bound 
\eqref{e:101} in order to reduce the smoothness 
assumptions\footnote{We show \eqref{e:101} under the condition $A\in Y_1$ and then interpolate with \eqref{e:100} to obtain \eqref{e:115} under the less restrictive condition 
$A\in \tilde{Y}_1$.} on our vector potentials, 
see Section \ref{sec:6}.\\
\noindent We  need an extension of Proposition \ref{prop:12}, 
which follows from the same proof as in \cite{KeelTao}. 
\begin{proposition}
\label{prop:14}
Let $W(t)$ is an operator defined on all Schwartz 
functions on $\rn$ and it  satisfies 
\begin{eqnarray*}
& & \norm{W(t) f}{L^2_{x_2, 
\ldots, x_n} L^2_{x_1}}=\norm{W(t) f}{L^2_x}\leq C\norm{f}{L^2_x} \\
& & \norm{W(t)W(s)^* f}{L^\infty_{x_2, 
\ldots, x_n} L^2_{x_1}}\leq C |t-s|^{-\si} 
\norm{f}{L^1_{x_2, 
\ldots, x_n} L^2_{x_1}}.
\end{eqnarray*} 
Then 
\begin{equation}
\label{eq:op}
\norm{W(t) f}{L^q_t  L^r_{x_2, \ldots, x_n}L^2_{x_1}}\leq C 
\norm{f}{L^2},
\end{equation}
for all  $\si$ admissible pairs $(q,r)$. 
The usual averaging argument then implies the ``retarded estimate''
\begin{eqnarray*}
& & \norm{\intl_0^t W(t)W(s)^* 
F(s, \cdot) ds}{L^q_t  L^r_{x_2, \ldots, x_n}L^2_{x_1}}\leq 
C\norm{F}{L^1 L^2}.
\end{eqnarray*}
\end{proposition}
\noindent Indeed,  for the ``retarded estimate'', assume \eqref{eq:op} to get 
\begin{eqnarray*}
& & \norm{\intl_0^t W(t)W(s)^* 
F(s, \cdot) ds}{L^q_t  L^r_{x_2, \ldots, x_n}L^2_{x_1}} \leq 
\norm{\intl_0^t \norm{W(t)W(s)^* 
F(s, \cdot)}{ L^r_{x_2, \ldots, x_n}L^2_{x_1}} ds}{L^q_t}\lesssim \\
& & \lesssim  \intl_0^\infty \norm{W(t)[W(s)^* 
F(s, \cdot)]}{ L^q_t L^r_{x_2, \ldots, x_n}L^2_{x_1}} ds \lesssim 
\intl_0^\infty  \norm{W(s)^* F(s, \cdot)}{L^2} ds\leq C \norm{F}{L^1_t L^2_x},
\end{eqnarray*}
where  we have used \eqref{eq:op} in the form 
$\norm{W(t)G(s,\cdot) }{L^q_t  L^r_{x_2, \ldots, x_n}L^2_{x_1}}\lesssim 
\norm{G(s, \cdot)}{L^2}$ as well as  the energy estimate 
$\norm{W(s)^* F}{L^2}\leq C \norm{F}{L^2_x}$.

\section{Proof of Theorem \ref{theo:pisa}}
\label{sec:4}
\subsection{Energy estimates and Littlewood-Paley reductions}
\label{sec:3}
To start our argument, we shall need the following $L^2$ 
existence result, due to 
Kenig-Ponce-Vega, \cite{KPV98}, which generalizes 
an earlier work of Doi, \cite{Doi1}, \cite{Doi2}. We state it here 
only in the particular case of interest to us, namely  Schr\"odinger equation with 
first order perturbations.
\begin{proposition}(Kenig-Ponce-Vega, \cite{KPV98})
\label{prop:KPV}
For the equation
$$
\left|
\begin{array}{l}
\p_t u-i \De u+\vec{b}(t,x)\cdot \nabla u= F\\
u(x,0)=f
\end{array}
\right.
$$
there is an unique global solution, provided $\vec{b}\in C^N$, 
$|Im\ \vec{b}|\leq <x>^{-m}$ for some large integers $N,m$. 
Moreover the solutions 
are smooth, provided $f, F$ are smooth  and for 
every $T>0$
$$
\norm{u}{L^\infty(0,T)L^2}\leq C_T (\norm{f}{L^2}+\norm{F}{L^1(0,T)L^2}).
$$
\end{proposition}
\noindent In the case, when $\vec{b}$ is real valued and 
$\|\vec{b}\|_{L^1L^\infty}<1/2$, 
one can derive {\it a priori} estimates for 
$\norm{u}{L^\infty(0,T)L^2}$ that are $T$  independent. \\
This is very standard energy estimate. 
Indeed, multiply both sides by $\bar{u}$, integrate 
in the spatial variable and take real 
part\footnote{By the smoothness of  the solutions all the operations are justified.}. 
We obtain 
\begin{eqnarray*}
\p_t \int |u|^2 dx+ \int \vec{b}\cdot\nabla |u|^2dx=\f{1}{2}\int (F \bar{u}+\bar{F} u)dx
\end{eqnarray*}
Integrate by parts and then integrate in $(0,T)$ to get
\begin{eqnarray*}
\norm{u(T,\cdot)}{L^2}^2-\norm{f}{L^2}^2&=&\intl_0^T \int div(\vec{b}) 
|u|^2dxdt+ \intl_0^T \int (F \bar{u}+\bar{F} u)dx dt
\leq\\
&\leq& \|\nabla \vec{b}\|_{L^1L^\infty}\supl_{0\leq t\leq T}\norm{u(t,\cdot)}{L^2}^2 
+\norm{F}{L^1L^2}\supl_{0\leq t\leq T}\norm{u(t,\cdot)}{L^2}.
\end{eqnarray*}
whence since $\|\nabla\vec{b}\|_{L^1L^\infty}<1/2$, 
$$
\norm{u(T,\cdot)}{L^\infty L^2}\leq C( \norm{f}{L^2}+\norm{F}{L^1L^2}).
$$
Thus, we have shown the following 
\begin{proposition}
\label{prop:5}
Let $f, F$ be smooth functions. Let also 
$\vec{A}$ be a smooth, real-valued potential with $\|\nabla\vec{A}\|_{L^1L^\infty}<1/2$. 
Then the Schr\"odinger equation 
$$
\left|
\begin{array}{l}
\p_t u-i \De u+\vec{A}(t,x)\cdot \nabla u= F\\
u(x,0)=f
\end{array}
\right.
$$
has an unique global solution and moreover there exists an absolute constant $C$, 
so that for every $T>0$ 
\begin{equation}
\label{eq:2300}
\norm{u(T,\cdot)}{L^\infty L^2}\leq C( \norm{f}{L^2}+\norm{F}{L^1L^2}).
\end{equation}
\end{proposition}
We may restate  Proposition \ref{prop:5} in a 
slightly different manner. Namely, 
the linear operators $U_A(t,s):L^2(\rn)\to L^2(\rn)$, 
where  $U_{A}(t,s)f$ is the unqiue solution $u$ of 
$$
\left|
\begin{array}{l}
\p_t u-i \De u+\vec{A}(t,x)\cdot \nabla u= 0\\
u(s,x)=f
\end{array}
\right.
$$
are well-defined. Moreover, by uniqueness and 
since the equation is time reversible, we can define $U_A(s,t):=U_A(t,s)^{-1}$, which is the solution 
operator to the same equation with data at time $t$ backwards in time to $s$.\\
The Duhamel's formula may be used to write the unique solution to \eqref{eq:1} as 
\begin{equation}
\label{eq:duh}
u(t,x)=U_A(t,0)f+\intl_0^t U_A(t,s) F(s,\cdot) ds.
\end{equation}
\noindent 
Next, we take a Littlewood-Paley projections  of \eqref{eq:1}. 
We get the equations 
\begin{equation}
\label{eq:260}
\left|
\begin{array}{l}
\p_t u_k - i \De u_k + \vec{A}_{\leq k-4} \cdot \nabla u_k= F_k + E^k\\
u_k(0)=f_k
\end{array}
\right.
\end{equation}
where $E^k$ is the error term
\begin{eqnarray}
\label{eq:1201}
E^k &=&   [P_k, \vec{A}_{\leq k-4}] \nabla u_{k-1\leq \cdot k+1} +
\suml_{l\geq k-2} P_k(\vec{A}_l \cdot\nabla u_{l-2\leq \cdot\leq l+2})+\\
&+& \nonumber \suml_{l\geq k-2} P_k(\vec{A}_{l-2\leq\cdot\leq l+2}  \cdot\nabla u_{l})+
P_k ( \vec{A}_{k-1\leq\cdot\leq k+1}\cdot \nabla u_{\leq k-4}).
\end{eqnarray}
Denote 
\begin{eqnarray*}
& & \cl \psi:= \psi_t - i \De \psi+A\cdot \nabla \psi, \\
& & \cl^k\psi:= \psi_t - i \De \psi+A_{\leq k-4} \cdot \nabla \psi.
\end{eqnarray*}
Our next observation is a naive Strichartz estimate for $\cl$, 
which will be the  starting point in a continuity argument later on. 
\begin{proposition}
\label{prop:7}
For a fixed integer  $k_0$, there exists a time 
$T_0=T_0(k_0)\leq \infty$, so that 
whenever $0<T< T_0$, and for every $\psi\in \cs$  
\begin{equation}
\label{eq:1800}
\norm{P_{<k_0}\psi}{L^q(0,T)L^r}\leq C(T,k_0)(\norm{\psi(0,\cdot)}{L^2}+
\norm{\cl\psi}{L^1L^2}).
\end{equation}
Moreover, $C(T,k_0)$ depends on  $T$ in a continuous way.
\end{proposition}
\begin{proof}
The proof is based on the standard Strichartz estimates for the 
linear \\ Schr\"odinger equation, as we treat the term 
$\vec{A}\nabla\psi$ as a perturbation.  
By H\"older's inequality and since 
$\norm{\nabla P_{<k_0}\psi}{L^p}\lesssim 2^{k_0}\norm{\psi}{L^p}$   
\begin{eqnarray*}
& &\norm{P_{<k_0}\psi}{L^q(0,T)L^r}\leq C(\norm{\psi(0)}{L^2}+ \|\vec{A}
\nabla P_{<k_0}\psi\|_{L^1(0,T)L^2}+
\norm{\cl P_{<k_0}\psi}{L^1L^2}) \\
&\leq&  C_1(\norm{\psi(0)}{L^2}+\sqrt{T}2^{k_0}\|\vec{A}\|_{L^2_t L_x^\infty}
\norm{\psi}{L^\infty(0,T)L^2})+\\
&+& C_1(
\norm{P_{<k_0}\cl\psi}{L^1L^2} 
+\norm{[P_{<k_0},\vec{A}]\nabla \psi}{L^1L^2})\leq \\
& &\leq C_1(\norm{\psi(0)}{L^2}+\norm{\cl\psi}{L^1L^2})+ \\
& &+
C_1(\sqrt{T}2^{k_0}\|\vec{A}\|_{L^2_t L_x^\infty}+
\norm{\nabla \vec{A}}{L^1L^\infty})\norm{\psi}{L^\infty(0,T)L^2} \\
& & \leq C_1(\norm{\psi(0)}{L^2}+\norm{\cl\psi}{L^1L^2})+ 
C_1\ve(\sqrt{T}2^{k_0}+1)\norm{\psi}{L^\infty(0,T)L^2}) .
\end{eqnarray*}
According to \eqref{eq:2300}, the last expression is controlled by 
$$
C_1(1+C \ve (\sqrt{T}2^{k_0}+1)) (\norm{\psi(0)}{L^2}+\norm{\cl\psi}{L^1L^2}).
$$
Hence, \eqref{eq:1800} holds with $C(T,k_0)= C_1(1+C \ve (\sqrt{T}2^{k_0}+1))$.
\end{proof}
\noindent Fix $k_0$ and a small $\ve=\ve(n)$ (to be chosen later). 
Set  $0<T^*\leq \infty$ to be 
the maximum time,  so that for all 
 $0<T<T^*$ and for all $\psi\in \cs$: 
\begin{equation}
\label{eq:iteration}
\norm{P_{<k_0}\psi}{L^q(0,T)L^r}\leq \ve^{-1}
(\norm{\psi(0,\cdot)}{L^2}+\norm{\cl \psi}{L^1L^2}).
\end{equation}
If $T^*=\infty$, \eqref{eq:iteration} holds 
for the fixed $k_0$. \\
We will show that  $T^*<\infty$ leads to a 
contradiction, provided $\ve=\ve_n$ was 
chosen suitably small.  We need to consider  separately the homogeneous 
and inhomogeneous problems.
\subsection{Homogeneous problem for $\cl$}
Consider
$$
\left|
\begin{array}{l}
\cl g=0\\
g(0,x)=f(x),
\end{array}
\right.
$$
where $g$ is constructed as $g(t,x)=U_{\vec{A}}(t,0)f$. 
We will show that there exists a 
constant $C$, 
depending on the dimension\footnote{The constant $C$ 
may also depend on the Strichartz pair $(q,r)$ in dimension two, 
although we suppress that dependence.} $n$, but not on $T$, $k_0$, or $f$, 
 so that for any $T<T^*$ and all $k\leq k_0$
\begin{equation}
\label{eq:3099}
(\suml_k \norm{P_{<k_0}P_k g}{L^q(0,T)L^r}^2)^{1/2}\leq C \norm{f}{L^2}.
\end{equation}
An elementary computation shows that $g_k$ solves 
\begin{equation}
\label{eq:1020}
\left|
\begin{array}{l}
\cl g_k=
-[P_k , \vec{A}] \nabla g, \\
g_k (0,x)=f_k.
\end{array}
\right.
\end{equation}
\subsection{A priori $L^2$ estimates for \eqref{eq:1020}}
First, we show an {\it a priori} estimate 
for \\ $(\suml_k \norm{g_k}{L^\infty L^2}^2)^{1/2}$ which generalizes
 \eqref{eq:2300}. \\
Apply \eqref{eq:2300} for the solutions of \eqref{eq:1020}. We get 
\begin{eqnarray*}
& &\norm{g_k}{L^\infty L^2}\leq C(\norm{f_k}{L^2}+
\norm{[P_k , \vec{A}] \nabla g}{L^1L^2}).
\end{eqnarray*}
To tackle the term $\norm{[ P_k , \vec{A}] \nabla g}{L^1L^2}$, 
we split $\vec{A}$ in low and high frequencies portions  
respectively. We have 
$$
\suml_k \norm{[ P_k , \vec{A}] \nabla g}{L^1L^2}^2\leq
\suml_k \norm{[ P_k , \vec{A}_{\leq k-4}] \nabla g}{L^1L^2}^2+
\suml_k \norm{[ P_k , \vec{A}_{>k-4}] \nabla g}{L^1L^2}^2.
$$
Since $[ P_k , \vec{A}_{\leq k-4}] \nabla g= 
[ P_k , \vec{A}_{\leq k-4}]P_{k-2<\cdot<k+2} \nabla g$, 
we estimate this  portion by
\begin{eqnarray*}
\suml_k \norm{[P_k,\vec{A}_{\leq k-4}] \nabla g}{L^1L^2}^2\lesssim 
\norm{\nabla \vec{A}}{L^1L^\infty}^2 \suml_k
\norm{g_{k-2<\cdot<k+2}}{L^\infty L^2}^2\lesssim  \ve^2
\suml_k
\norm{g_k}{L^\infty L^2}^2.
\end{eqnarray*}
For the high-frequency portion, we have 
$$
[ P_k , \vec{A}_{>k-4}] \nabla g=-\vec{A}_{>k-4}\nabla g_k+ 
P_k(\vec{A}_{>k-4} \nabla g).
$$
 The first term is trivially estimated by 
$$
\suml_k \norm{\vec{A}_{>k-4}\nabla g_k}{L^1L^2}^2\lesssim 
\norm{\nabla \vec{A}}{L^1L^\infty}^2\suml_k\norm{g_k}{L^\infty L^2}^2.
$$ 
The rest is estimated 
 by
\begin{eqnarray*}
 & & \suml_k \norm{P_k(\vec{A}_{>k-4} \nabla g)}{L^1L^2}^2\leq 
\suml_k \norm{P_k(\vec{A}_{k-4<\cdot<k+4} \nabla g)}{L^1L^2}^2+\\
& &+
\suml_k \norm{P_k(\vec{A}_{>k+4} \nabla g)}{L^1L^2}^2
\lesssim \suml_k \norm{\vec{\nabla A}_{k-4<\cdot<k+4}}{L^1L^\infty}^2
\norm{g}{L^\infty L^2}^2+\\
& &+
 \suml_k\left(\suml_{l\geq k+4} \norm{P_k( \vec{A}_l \nabla 
g_{l-2\leq \cdot \leq l+2})}{L^1L^2}\right)^2\lesssim\\
& &\lesssim \ve^2\norm{g}{L^\infty L^2}^2+ 
\suml_k\left(\suml_{l\geq k+4} 2^{h k} \norm{P_k( \vec{A}_l \nabla 
g_{l-2\leq \cdot \leq l+2})}{L^1L^{2n/(n+2h)}}\right)^2\lesssim \\
& &\lesssim \ve^2\norm{g}{L^\infty L^2}^2+\suml_k 2^{2 h k} 
\left(\suml_{l\geq k+4} 2^l \|\vec{A}_l\|_{L^1 L^{n/h}}
\norm{g_{l-2\leq \cdot \leq l+2}}{L^\infty L^2}\right)^2\lesssim \\
& &\lesssim \ve^2\norm{g}{L^\infty L^2}^2+ \suml_k 2^{2 h k} 
\left(\suml_{l\geq k+4} 2^{-h l} \|\vec{A}_l\|_{L^1 W^{n/h, 1+h}}
\norm{g_{l-2\leq \cdot \leq l+2}}{L^\infty L^2}\right)^2
\end{eqnarray*}
for some $h: 0<h<1/4$. 
By Lemma \ref{le:sequences}, we get that the second
 expression above is bounded by 
$
C\sup_l \norm{g_{l-2\leq \cdot \leq l+2}}{L^\infty L^2}^2 
\suml_l \|\vec{A}_l\|_{L^1 W^{n/h, 1+h}}^2.$\\
Altogether, we obtain the estimate 
\begin{eqnarray*}
& & \suml_k \norm{g_k}{L^\infty L^2}^2\leq C \norm{f}{L^2}^2+ C\ve^2 
\suml_k \norm{g_k}{L^\infty L^2}^2,
\end{eqnarray*}
whence for $\ve$ small enough
\begin{equation}
\label{eq:1021}
(\suml_k \norm{g_k}{L^\infty L^2}^2)^{1/2}\leq 2 C \norm{f}{L^2},
\end{equation}
which generalizes \eqref{eq:2300}, since $\norm{g}{L^\infty L^2}\leq 
(\suml_k \norm{g_k}{L^\infty L^2}^2)^{1/2}$.

We  also need the following crucial lemma, 
whose proof we postpone for  Section \ref{sec:par}. 
\begin{lemma}(existence of parametrix)
\label{le:parametrix}
Given a  potential 
$\norm{\vec{A}}{Y_1\cap Y_2\cap Y_3}\leq \ve$ and integer $k$ and $T>0$, 
and for every function $f_k\in L^2(\rn)$ with $supp\ \widehat{f_k}
\subset \{|\xi|\sim 2^k\}$, 
one can find a 
function 
$v_k:[0,T)\times \rn\to \cc$, so that $supp\ 
\widehat{v_k}\subset \{|\xi|\sim 2^k\}$ and 
\begin{eqnarray*}
& &\norm{v_k(0,x)-f_k}{L^2}\leq C\ve \norm{f_k}{L^2} \\
& &\norm{v_k}{L^q(0,T)L^r}\leq C\norm{f_k}{L^2}\\
& &\norm{\cl v_k}{L^1L^2}\leq C\ve \norm{f_k}{L^2},
\end{eqnarray*}
for some $C$ independent of $f, k, T$.
\end{lemma}

\noindent Assuming for the moment the validity
 of Lemma \ref{le:parametrix}, we show \eqref{eq:3099}.
Take $v_k$ as in  Lemma \ref{le:parametrix}.  
Applying the {\it a priori} estimate \eqref{eq:iteration} yields 
\begin{eqnarray*}
& &\norm{P_{<k_0} P_k g}{L^q(0,T)L^r}\leq 
\norm{P_{<k_0}(P_k g-v_k)}{L^q(0,T)L^r}+\norm{P_{<k_0} v_k}{L^q(0,T)L^r} 
\leq \\
& & \leq \ve^{-1}(\norm{f_k-v_k(0,x)}{L^2}+\norm{\cl P_k g}{L^1L^2}+
\norm{\cl v_k}{L^1L^2})+C\norm{v_k}{L^q(0,T)L^r}\leq \\
& &\leq \ve^{-1}(C\ve\norm{f_k}{L^2}+
\norm{[ P_k , \vec{A}] 
\nabla g}{L^1L^2})+C\norm{f_k}{L^2}.
\end{eqnarray*}
For the proof of \eqref{eq:1021}, we have already estimated 
$$
\suml_k \norm{[ P_k , \vec{A}] 
\nabla g}{L^1L^2}^2\leq C\ve^2 \suml_k \norm{g_k}{L^\infty L^2}^2,
$$
which by \eqref{eq:1021} is bounded by $C \norm{f}{L^2}^2$. 
Squaring and summing the estimates 
for \\ $\norm{P_{<k_0} P_k g}{L^q(0,T)L^r}$ yields 
$$
\norm{P_{<k_0} g}{L^q(0,T)L^r}\leq 
(\suml_k \norm{P_{<k_0} P_k g}{L^q(0,T)L^r}^2)^{1/2}\leq C\norm{f}{L^2}.
$$
This proves 
\eqref{eq:3099}. 
 Rewriting \eqref{eq:3099} in terms of the operators $U$ yields 
$$
(\suml_k \norm{P_{<k_0} P_k U(t,0)f}{L^q(0,T)L^r}^2)^{1/2}\leq C\norm{f}{L^2},
$$ 
or more generally since 
$U(t,s)=U(t,0)U(0,s)$
\begin{equation}
\label{eq:3100}
(\suml_k \norm{P_{<k_0} P_k U(t,s)f_s}{L^q(0,T)L^r}^2)^{1/2}\leq 
C\norm{U(0,s)f_s}{L^2}\leq C \norm{f_s}{L^2}.
\end{equation}
which is sometimes more convenient to use. In particular, we obtain 
\begin{equation}
\label{eq:3105}
\norm{P_{<k_0} U(t,s)f_s}{L^q(0,T)L^r}\leq 
C \norm{f_s}{L^2}.
\end{equation}
\subsection{The inhomogeneous problem for $\cl$}
We derive estimates similar to \eqref{eq:3099} for the 
inhomogeneous problem associated with  $\cl$. Namely, we consider 
$$
\left|
\begin{array}{l}
\cl w=G\\
w(0,x)=0,
\end{array}
\right.
$$
We show that the solution $w$ 
(constructed by the Duhamel's formula \eqref{eq:duh}) satisfy 
\begin{equation}
\label{eq:3101}
(\suml_k \norm{P_{<k_0}P_k w}{L^q(0,T)L^r}^2)^{1/2}
\leq C\norm{G}{L^1(0,T)L^2},
\end{equation}
whenever $T<T^*$ and the constant $C$ is dependent only on  $n$.\\
Our first step as in the homogeneous case is to project the inhomogeneous 
equation by the Littlewood-Paley operator $ P_k$. We have that 
\begin{equation}
\label{eq:3102}
\cl  P_k w= P_k G - [ P_k, A] \nabla w,
\end{equation}
whence exactly as in the homogeneous case 
(and by using the smallness assumptions on $\vec{A}$), 
we conclude 
\begin{eqnarray*}
& &
(\suml_k \norm{P_k w}{L^\infty L^2}^2)^{1/2}
\leq C (\suml_k \norm{G_k}{L^1L^2}^2)^{1/2}\leq C \norm{G}{L^1L^2},\\
& & (\suml_k \norm{[ P_k, A] \nabla w}{L^1L^2}^2)^{1/2}\leq C \ve 
(\suml_l  \norm{P_l w}{L^\infty L^2}^2)^{1/2}\leq C\ve 
\norm{G}{L^1L^2}.
\end{eqnarray*}

By the Duhamel's 
formula  applied to \eqref{eq:3102} and \eqref{eq:3105}
\begin{eqnarray*}
& &(\suml_k \norm{P_{<k_0} P_k w}{L^q(0,T)L^r}^2)^{1/2}=
(\suml_k\left\|\int_0^t P_{<k_0} U(t,s) ( P_k
G(s,\cdot)) ds\right\|_{L^q(0,T)L^r}^2)^{1/2} +\\
& & (\suml_k\left\|\int_0^t P_{<k_0} U(t,s)( [P_k,\vec{A}] 
\nabla w)ds\right\|_{L^q(0,T)L^r}^2)^{1/2}\lesssim \\
& &\lesssim 
(\suml_k (\intl_0^T  \norm{P_k
G(s,\cdot)}{L^2}
ds)^2)^{1/2}+ (\suml_k (\intl_0^T  \norm{[P_k,\vec{A}] 
\nabla w(s,\cdot))}{L^2}
ds)^2)^{1/2}\lesssim \\
& &
\lesssim (\suml_k \norm{G_k}{L^1L^2}^2)^{1/2}+ 
(\suml_k \norm{[P_k, A] \nabla w}{L^1L^2}^2)^{1/2}
\lesssim \norm{G}{L^1L^2}.
\end{eqnarray*}
In particular,
\begin{equation}
\label{eq:3106}
\norm{P_{<k_0} w}{L^q(0,T)L^r}\leq C \norm{G}{L^1L^2}.
\end{equation}
Combining \eqref{eq:3105}, 
\eqref{eq:3106} yields that the solution to 
$\cl \psi =G, \q \psi(0,x)=f$, satisfies 
$$
\norm{P_{<k_0}\psi}{L^q(0,T)L^r}\leq C_n(\norm{f}{L^2}+\norm{G}{L^1L^2}),
$$
which is a contradiction with the maximality of $T^*$ 
(see \eqref{eq:iteration}), provided $\ve < 1/C_n$. 
Thus $T^*=\infty$ and one has the inequality 
$$
\norm{P_{<k_0}\psi}{L^qL^r}\leq C_n(\norm{\psi(0,\cdot)}{L^2}+
\norm{\cl \psi}{L^1L^2}).
$$ 
for all $k_0$. Taking a limit $k_0\to \infty$ establishes 
the Strichartz estimate \eqref{eq:1009}. Moreover,  
we have established a more general Besov  type estimate 
$$
(\suml_k \norm{\psi_k}{L^q L^r}^2)^{1/2}
\leq C_n(\norm{\psi(0,\cdot)}{L^2}+
\norm{\cl \psi}{L^1L^2}).
$$
\section{The construction of the parametrix}
\label{sec:par}
\noindent In this section, we show the existence of approximate solution 
(in  sense of Lemma \ref{le:parametrix}) to the 
equation $\cl g=0, \q g(x,0)=f_k$. \\
We first make some reductions. 
Our first  reduction is that 
it will suffice to show that there exists 
$v_k$ with the properties as listed in Lemma \ref{le:parametrix}, 
where $v_k$ satisfies $\norm{\cl^k v_k}{L^1 L^2}\leq C\ve\norm{f_k}{L^2}$ 
instead of $\norm{\cl v_k}{L^1 L^2}\leq C\ve \norm{f_k}{L^2}$. Indeed, 
suppose there is $v_k$ with $\norm{\cl^k v_k}{L^1 L^2}\leq 
C\ve\norm{f_k}{L^2}$. Then 
\begin{eqnarray*}
& &
\norm{\cl v_k}{L^1L^2}\leq 
\norm{\cl^k v_k}{L^1 L^2}+\norm{\vec{A}_{>k-4}}{L^1L^\infty}
\norm{\nabla v_k}{L^\infty L^2}\leq\\
& &\leq \norm{\cl^k v_k}{L^1L^2}+C \norm{\p \vec{A}}{L^1L^\infty} 
\norm{v_k}{L^\infty L^2}\leq C\ve \norm{f_k}{L^2}
\end{eqnarray*}
Suppose that one has already  a 
function  $v_k$ as in Lemma 
\ref{le:parametrix} without  
$supp\ \widehat{v_k}\subset \{|\xi|\sim 2^k\}$. 
Our claim is that $\tilde{P}_k v_k$ will satisfy 
all conditions in Lemma 
\ref{le:parametrix}, 
with $\cl$ replaced by $\cl^k$ according to our first 
reduction. \\
By construction $supp\ \widehat{\tilde{P}_k v_k}
\subset \{|\xi|\sim 2^k\}$.
Since  $\tilde{P}_k f_k=f_k$, 
$$
\norm{\tilde{P}_k v_k(0,x)-f_k}{L^2}=\norm{\tilde{P}_k 
(v_k(0,x)-f_k)}{L^2}\leq 
\norm{v_k(0,x)-f_k}{L^2}\leq \ve\norm{f}{L^2}.
$$
Next, 
$$
\|\tilde{P}_k v_k\|_{L^qL^r}\leq \norm{v_k}{L^qL^r}\leq C\norm{f_k}{L^2}.
$$
Finally, since $\cl^k \tilde{P}_k v_k=\cl^k v_k-
[ \tilde{P}_k, \vec{A}_{\leq k-4}] \nabla v_k$
\begin{eqnarray*}
& &
\norm{\cl \tilde{P}_k v_k}{L^1L^2}\leq \norm{\cl v_k}{L^1L^2}+
\norm{[ \tilde{P}_k, \vec{A}_{\leq k-4}] \nabla v_k}{L^1L^2}\leq \\
& &
\leq C\ve \norm{f_k}{L^2}+
\norm{\p \vec{A}}{L^1L^\infty}\norm{v_k}{L^\infty L^2}\leq 
C\ve \norm{f_k}{L^2}.
\end{eqnarray*}

Next, since all our estimates will be  
scale invariant, we may rescale and assume 
without loss of generality  $k=0$. Thus, 
matters  are reduced to the following 
\begin{lemma}
\label{le:parametrixmod}
Let $\ve>0$ and $\norm{\vec{A}}{Y_1\cap Y_2\cap Y_3}<\ve$ be a potential with \\
$supp\ \widehat{A}(t,\xi)\subset \{|\xi|\lesssim 1\}$. 
Then for every $T>0$ and 
for every function $f\in L^2(\rn)$ with $supp\ 
\widehat{f}\subset \{|\xi|\sim 1 \}$, 
one can find a 
function 
$v:[0,T)\times \rn\to \cc$, so that
\begin{eqnarray}
\label{eq:902}
& &\norm{v(0,x)-f}{L^2}\leq C\ve \norm{f}{L^2} \\
\label{eq:903}
& &\norm{v}{L^q(0,T)L^r}\leq C\norm{f}{L^2}\\
\label{eq:904}
& &\norm{\p_t v- i \De v+ \vec{A}\cdot \nabla v}{L^1L^2}\leq 
C\ve \norm{f}{L^2},
\end{eqnarray}
for some $C$ independent on $f,  T, \ve$.
\end{lemma}

\begin{proof}(Lemma \ref{le:parametrixmod})
We construct $v$ in the form 
$$
v(t,x)=\La f(t,x)= \int e^{i \si(t,x,\xi)} e^{-4\pi^2 i t|\xi|^2} e^{2\pi i \dpr{\xi}{x}}
\Omega(\xi) \widehat{f}(\xi) d\xi,
$$
where the $\Omega$ is a smooth cut-off of the annulus $|\xi|\sim 1$ and the 
phase correction $\si$ is to be selected momentarily. We have 
\begin{eqnarray*}
& &\cl v= \p_t v - i \De v+\vec{A}\cdot \nabla v =\\
& &=\int (i \p_t \si+\De \si+ 2\pi i(\dpr{\nabla\si}{\xi}+
 \vec{A}\cdot \xi)+ i 
[(\nabla\si)^2+\vec{A}\cdot \nabla \si])  \times \\
& &\times 
e^{i \si(t,x,\xi)} e^{-4\pi^2 i t|\xi|^2} e^{2\pi i \dpr{\xi}{x}}
\Omega(\xi) \widehat{f}(\xi) d\xi.
\end{eqnarray*}
We first comment on possible choices for $\si$. Clearly, 
since we are trying to almost solve $\cl v=0$, 
we should choose $\si$ in a way, so that the main terms are 
resolved  in  the formula for $\cl v$.  
We see that since the potential $\vec{A}$ is 
supported in the low frequencies and is small, the main terms are 
those, that are either  linear in $\vec{A}$ or linear in $\nabla \si$. 
It seems then  reasonable to choose $\si$, so that 
$$
\dpr{\nabla\si}{\xi}+ \vec{A}\cdot\xi=0.
$$
However, it turns out that when $|\dpr{\xi/|\xi|}{\eta}|\lesssim |\eta|^2$ 
(here $\eta$ is the Fourier variable for $\si$), 
one has that $\De\si$ is actually ``bigger'' 
compared to $\dpr{\nabla\si }{\xi}$. We therefore modify our 
choice as follows.\\
Set $\si=\si^0+\si^1$, where
\begin{eqnarray*}
& &
\si^0(t,x,\xi)=\suml_{k\leq -2}   \intl_0^\infty \vec{A}_k(t, x+z\xi/|\xi|) 
\cdot \f{\xi}{|\xi|} \chi(2^{2k}z) dz,\\
& &
\si^1(t,x,\xi)= 2\pi i \suml_{k\leq -2} 2^{2k} 
\int_0^\infty \De^{-1}\vec{A}_k(t, x+z\xi/|\xi|)\cdot \xi\chi'(2^{2k}z) dz.
\end{eqnarray*}
It is  elementary to see that 
$$\dpr{\nabla\si^0}{\xi}+ 
\vec{A}\cdot\xi=- \suml_k 2^{2k}\intl_0^\infty 
\dpr{\vec{A}_k(t, x+z\xi/|\xi|)}{\xi}
\chi'(2^{2k}z) dz,
$$
 whence 
$$
\De \si^1+2\pi i(\dpr{\nabla\si^0}{\xi}+ \vec{A}\cdot\xi)=0.
$$
Denote $\vec{\tilde{A}}_k:=2^{2k}\De^{-1} \vec{A}_k$. Then, 
we rewrite 
$$
\si^1(t,x,\xi)= 2\pi i \suml_{k\leq -2} 
\int_0^\infty \vec{\tilde{A}}_k(t, x+z\xi/|\xi|)\cdot \xi\chi'(2^{2k}z) dz.
$$
According to the choice of $\si^0, \si^1$, we get 
\begin{equation}
\label{eq:501}
\begin{array}{l}
\cl (\La f) (x,t) = \int (i \p_t (\si^0+\si^1) 
+\De \si^0 + 2\pi i \dpr{\nabla\si^1}{\xi}+ i 
[(\nabla\si)^2+\dpr{\vec{A}}{\nabla \si}])  \times \\
\times 
e^{i \si(t,x,\xi)} e^{-4\pi^2 i t|\xi|^2} e^{2\pi i \dpr{\xi}{x}}
\Omega(\xi) \widehat{f}(\xi) d\xi.
\end{array}
\end{equation}
Now every term in the formula for $\cl (\La f)$ 
except $\dpr{\nabla\si^1}{\xi}$ either has a 
two spatial derivatives or one time  derivative acting on it (recall 
that in our scaling time derivatives are worth two spatial derivatives) or is 
quadratic in $\nabla\si$ ( since $\nabla \si \sim A$ ). \\
However by our choice of $\si^1$, we have 
\begin{eqnarray}
& & \nonumber
\dpr{\nabla\si^1}{\xi} =  2\pi i
\suml_{k\leq -2} 
\intl (\suml_{j, m=1}^n \p_m \vec{\tilde{A}}^j_k(x+z \xi/|\xi|) \xi_j \xi_m )
\chi'(2^{2k} z) dz=\\
\nonumber 
& &=2\pi i
\suml_{k\leq -2}  \intl\left( \f{d}{d z} 
\vec{\tilde{A}}_k(x+z\xi/|\xi|) \cdot \xi\right)|\xi|  \chi'(2^{2k} z) dz=\\
\label{e:1}
& &= - 2\pi i
\suml_{k\leq -2}  2^{2k}\intl 
\vec{\tilde{A}}_k(x+z\xi/|\xi|) \cdot \xi |\xi|  \chi''(2^{2k}z) dz.
\end{eqnarray}
In this last expression, one has 
 multiplication by $2^{2k}$, which behaves 
like two spatial derivatives on 
$\vec{\tilde{A}}_k$. \\
Now that we have made our selection of 
 $\si$, we go back to the proof of Lemma \ref{le:parametrixmod}. 
First, expanding the exponential $e^{i \si}$ in 
Taylor series yields the representation \\
$\La = \suml_{\al\geq 0} i ^{\al}(\al! )^{-1} 
\La^{\al}$, where 
$$
\La^{\al}f(t,x)=\int (\si(t,x,\xi))^\al  
e^{-4\pi^2 i t|\xi|^2} e^{2\pi i \dpr{\xi}{x}}
\Om(\xi) \widehat{f}(\xi) d\xi.
$$
and similar for the expression for $\cl(\La f)$. It is clear now that in 
this formulation, it is convenient to think that  $\si$ is in the form 
$$
\si(t,x,\xi) \sim 
\suml_{k\leq -2}   \intl_0^\infty A_k(t, x+z\xi/|\xi|) 
 \chi(2^{2k}z) dz
$$
for some $C^\infty_0$ 
function 
$\chi$ supported in $(-2,2)$. 
 This is done by subsuming  the harmless terms $\xi/|\xi|$, $\xi$ and 
$\xi|\xi|$ in the multiplier $\Omega(\xi)$ and by considering the resulting 
expressions componentwise. This will be a good 
strategy for 
all estimates involving $L^2_x$ 
norms.\footnote{We will omit the interval of integration $(0,\infty)$ 
in the formula for $\si$. In other words, we will tacitly replace 
$\chi(2^{2k}z)$ by $\chi_+(2^{2k}z):=\chi(2^{2k}z)\chi_{(0,\infty)}(z)$.
This is not going to make any difference in the $L^2$ estimates, since 
no smoothness of the amplitude $\chi_+(2^{2k}z)$ is  needed.} 

We will however need to show also 
decay estimates for $\La^\al f$,  in which case,  it is better to 
think of $\si$ in the form  
$$
\si(t,x,\xi)\sim  \suml_{k\leq -2} \suml_{l\leq -2k} 
\intl_0^\infty A_k(t, x+z\xi/|\xi|) 
 \vp(2^{-l}z) dz.
$$
This is so since $\sum_{l\leq -2k} \vp(2^{-l} z)=\chi(2^{2k}z)$ and 
$\chi'(2^{2k}z)$ which enters in $\si^1$ has 
essentially the same form as $\vp(2^{2k}z)$. Note that since 
$\vp(2^{-l}z)$ has support away from the endpoints of the 
interval of integration, we can also  write 
$$
\si(t,x,\xi)\sim  \suml_{k\leq -2} \suml_{l\leq -2k} 
\intl A_k(t, x+z\xi/|\xi|) 
 \vp(2^{-l}z) dz,
$$
with some smooth compactly supported function $\vp$ with $supp\ \vp
\subset (1/2,4)$.
\section{Pointwise estimates on $\si$}
\label{sec:98}
 To start us off, we will need an additional 
angular decomposition for $\si$, which we describe next. Note that in this section, we completely ignore the $t$ dependence, since it is 
irrelevant in that setting. This is so, because in the spaces of 
interest to us (i.e. the Strichartz space and its dual)
 the $L^r_x$ norm always comes first.

To put ourselves in the desired framework, let us fix a family of unit 
vectors \\
$\{\theta_j^m\}_{j\in [1, c2^{m(n-1)}], 
m\in [-2, \infty)}$  with the property:
 For fixed $m$, 
the family of balls 
$\{B(\theta_j^m, 2^{-m})\}=:\{B_j^m\}$ forms a 
covering of $\sn$ 
with every ball in the family intersecting at most a 
fixed number
(depending only on the dimension n) 
 of other balls in the family.

In other words for fixed $m$, $\theta^m_j$ forms a $2^{-m}$ net over $\sn$ and 
the distance between any two $\theta^m_j, \theta^m_s$ is bounded below by 
$c2^{-m}$. 

For the locally finite covering $\{B(\theta_j^m, 2^{-m})\}$ of the annulus 
$\{\xi: 1-2^{-m-5}\leq |\xi|\leq 1+2^{-m-5}\}$, we find a smooth 
partition of unity subordinated to the covering. 
That is, there exists a familly 
$\{\psi_{j,m}\}\subset C^\infty_0$ with
$$
\suml_j \psi_{j,m}(2^m(\xi/|\xi|-\theta_j^m))=1.
$$
Moreover,  the functions $\psi_{j,m}$ 
satisfy uniform bound on their derivatives, 
i.e. 
$\supl_{j,m}|\p^\ga \psi_{j,m}(x)|\leq C_n^\ga$.

Fix $k,l: l>-k$ For any  function $H_k$, 
expand $H_k(x+z\xi/|\xi|)$ in Taylor series around $x+z\theta^{l+k}_j$, when 
$|\xi/|\xi|-\theta^{l+k}_j|\lesssim 2^{-l-k}$.
We get
\begin{eqnarray*}
& &
 H_k( x+z \xi/|\xi|)\psi_{j,l+k}(2^{l+k}(\xi/|\xi|-\theta_j^{l+k}))= \\
& &=
\suml_{\ga\geq 0} \f{1}{\ga!}
\p^\ga H_k( x+z \theta_j^{l+k})z^{|\ga|} (\xi/|\xi|-\theta^{l+k}_j)^\ga
\psi_{j,l+k}(2^{l+k}(\xi/|\xi|-\theta_j^{l+k}))
=\\
& &=\suml_{\ga\geq 0} \f{z^{|\ga|}}{2^{l|\ga|}\ga!} 
H^\ga_k( x+z \theta_j^{l+k})\psi^\ga_{j,l+k}(2^{l+k}(\xi/|\xi|-\theta_j^{l+k}))
\end{eqnarray*}
where $\psi_{j,m}^\ga (\mu)=\psi_{j,m}(\mu)\mu^{|\ga|}$. 
We will also adopt the convention of naming a  function
$H_k^\ga$, whenever  $H^\ga_k$ has the same Fourier support properties as 
$H_k$ and $\|H_k^\ga\|_{Z}\leq C^{|\ga|}
\|H_k\|_{Z}$ for all Banach spaces $Z$ used 
throughout the paper.
Using the formula  for \\ $H_k(x+z \xi/|\xi|)\psi_{j,l+k}(2^{l+k}
(\xi/|\xi|-\theta_j^{l+k}))$, 
we arrive at 
\begin{equation}
\label{eq:decompose}
\begin{array}{l}
\int H_k(x+z\xi/|\xi|)\psi_{j,l+k}(2^{l+k}(\xi/|\xi|-\theta_j^{l+k}))
\vp(2^{-l}z) dz= \\ 
=\suml_{\ga\geq 0} (\ga!)^{-1}
(\int H^\ga_k( x+z \theta^{l+k}_j)\vp^\ga
(2^{-l}z)dz) \psi_{j,l+k}^\ga(2^{l+k}(\xi/|\xi|-\theta_j^{l+k}))
\end{array}
\end{equation}
If $l\leq -k$,  
the above decomposition trivializes in the sense that we can write 
\begin{eqnarray*}
& & \int H_k(x+z\xi/|\xi|)
\vp(2^{-l}z) dz= \\
& & =\suml_{\ga\geq 0} (\ga!)^{-1}
(\int H^\ga_k( x+z e_1)\vp^\ga
(2^{-l}z)dz)(2^{l+k}(\xi/|\xi|-e_1))^\ga,
\end{eqnarray*}
that is, just the vector $e_1$ would suffice in that situation.

For the $L^2$ estimates, we use the decomposition 
\begin{equation}
\label{eq:54}
\begin{array}{l}
\int H_k(x+z\xi/|\xi|)\psi_{j,-k}(2^{-k}(\xi/|\xi|-\theta^{-k}_j))
\chi(2^{2k} z)
  dz= \\
 = \suml_{\ga\geq 0} (\ga!)^{-1}
(\int H^\ga_k( x+z \theta^{-k}_j)\chi^\ga(2^{2k}z) dz) 
\psi_{j,-k}^\ga(2^{-k}(\xi/|\xi|-\theta^{-k}_j)),
\end{array}
\end{equation}
which is derived similar to \eqref{eq:decompose}.

\noindent 
We have the following 
pointwise bound on  functions in the form  $\int H_k(x+\theta_j^{l+k}) 
\vp(2^{-l}z)dz$.
\begin{lemma}
\label{le:pointwise}
Let $k$ be fixed integer and $\vp$ be  a 
fixed Schwartz function with $supp \ \vp\subset (1/2, 4)$. 
Then 
\begin{equation}
\label{eq:254}
\supl_x \suml_{l>-k} \suml_j \int |H_k(x+z\theta_j^{l+k})| \vp(2^{-l}z) dz
\leq C_n
2^{k(n-1)}
\norm{H_k}{L^1_x}.
\end{equation}
In the case $l\leq -k$, we trivially have 
$$
\supl_x \suml_{l\leq -k} 
\suml_j \int |H_k(x+z e_1)| \vp(2^{-l}z) dz\leq C 2^{-k}\norm{H_k}{L^\infty}\leq C_n 2^{k(n-1)}\norm{H_k}{L^1_x}.
$$
\end{lemma}
\begin{proof} We concentrate on the case $l>-k$, since the other inequality 
is obvious. 

Represent $H_k(x) =2^{kn} \int
\varsigma(2^{2k}|x-y|^2) H_k(y) dy$ with some suitable Schwartz  function 
$\varsigma:\rone\to \rone$. Clearly 
\begin{eqnarray*}
& &2^{2k}|x+z\theta_j^{l+k}-y|^2=2^{2k}(|x-y|-z)^2+\\
& &+
2(2^{2k})(|x-y|z)(1-\dpr{\theta_j^{l+k}}{(y-x)/|y-x|}.
\end{eqnarray*}
We thus have 
$$
|\varsigma(2^{2k}|x+z\theta_j-y|^2)|\leq \f{C_N}{(1+2^{k}||x-y|-z|)^N
(1+2^{2k}(|x-y|z)\gamma^2)^N},
$$
where $\gamma$ is the angle between the unit vectors $\theta_j^{l+k}$ and 
$(y-x)/|y-x|$ and $N$ is arbitrary integer. It follows
that 
\begin{eqnarray*}
& & \suml_{l>-k} \suml_j \int |H_k(x+z\theta_j)| 
\vp(2^{-l}z) dz\leq \\
& &\leq \suml_{l>-k} \suml_j  \int  
\f{C_N 2^{k n}|H_k(y)|\vp(2^{-l}z) dz}{(1+2^{k}||x-y|-z|)^N
(1+2^{2k}(|x-y|z)\gamma^2)^N}  dydz=\\
& &=\suml_{l>-k} \suml_j  \int  
\f{C_N 2^{k n}
|H_k(x+r\theta)|r^{n-1}\vp(2^{-l}z) dz}{(1+2^{k}|r-z|)^N
(1+2^{2k}(r z)|\theta_j^{l+k}-\theta|^2)^N}  drd\theta dz
\end{eqnarray*}
The main term in the expression above is when the integration is 
over \\
$r\sim 2^{l}$, $|z-r|\leq 2^{-k}$ and $|\theta_j^{l+k}-\theta|\leq 2^{-k-l}$, 
with the corresponding decay away from this set. We  estimate by
\begin{eqnarray*}
& & C \suml_{l>-k} \suml_j 2^{k(n-1)} \int_{r\sim 2^l} 
\int_{|\theta_j^{l+k}-\theta|\leq 2^{-l-k}} |H_k(x+r\theta)||r|^{n-1} dr
\end{eqnarray*}
Clearly the summation in $j$ is 
sum of integrals over  (almost)  disjoint subsets of $\sn$ and as a result 
it  gives the integration over the whole $\sn$
 (since $\{\theta_j^{l+k}\}$ were chosen to be a $2^{-l-k}$ net of $\sn$). \\
We get 
\begin{eqnarray*}
& & \suml_{l>-k} \suml_j \int |H_k(x+z\theta_j^{l+k})| 
\vp(2^{-l}z) dz\lesssim 2^{k(n-1)} \int |H_k(x+y)| 
dy\leq 2^{k(n-1)} \norm{H_k}{L^1}.
\end{eqnarray*}
\end{proof}

\section{$L^qL^r$ estimates for the parametrix}
\label{sec:6}
We show that the parametrix is close to the initial 
data at $t=0$ and stays
 in the Strichartz space $L^qL^r$.
Taking into account that $\La^{0}f = e^{i t \De} f$, it is clear 
that \eqref{eq:902}, \eqref{eq:903} will follow from 
\begin{equation}
\label{eq:906}
\norm{\La^{\al}f}{L^qL^r}\leq C_{n}^\al
\norm{A}{\tilde{Y}_1}^{\al} \norm{f}{L^2}.
\end{equation}
The case $\al=0$ 
corresponds to the case of free solutions, which are  in $L^qL^r$ by the 
standard Strichartz estimates. 
We prove  \eqref{eq:906} by showing that 
 $\La^{\al}$ satisfies appropriate 
 decay and energy estimates.

\noindent We will show that for a fixed $s,t$ 
\begin{eqnarray}
\label{eq:90}
& &\norm{\La^{\al}f(t,\cdot)}{L^2}\leq 
C_n^\al\left(\suml_k \supl_{x\in\rn,\theta\in\sn}
\int |A_k(t,x+z\theta)|dz\right)^\al
\norm{f}{L^2}\\
\label{eq:91}
& & |\La^{\al}(t)\La^{\al}(s)^*f(x)|\leq C_n^{2\al}
(\suml_k 2^{k(n-1)}\norm{A_k}{L^\infty L^{1}})^{2 \al}
 |t-s|^{-n/2} \norm{f}{L^1},
\end{eqnarray}
whence by the abstract Strichartz estimates of 
Keel and Tao,  \cite{KeelTao}(see also 
Proposition \ref{prop:12}), one gets \eqref{eq:906}.

Note that by the Bernstein inequality  
\begin{eqnarray*}
& & \suml_k \supl_{x\in\rn,\theta\in\sn}\int |A_k(t,x+z\theta)|dz|
= \suml_k \supl_{U\in SU(n), x}\|A_k(t,x+Uz)\|_{L^\infty_t 
L^\infty_{z_2, \ldots, z_n}L^1_{z_1}} \lesssim \\
& &\lesssim  \suml_k 2^{k(n-1)}\|A_k\|_{L^\infty L^1_x}=\norm{A}{Y_1}
\end{eqnarray*} 
and therefore \eqref{eq:90} and \eqref{eq:91} hold for $A\in Y_1$. 

If however $n\geq 3$, we have by  complex multilinear 
 interpolation between 
\eqref{eq:90} and \eqref{eq:91} that  the ``modified decay estimate'' 
(see \eqref{e:115})
\begin{eqnarray*}
& &
\|\La^{\al}(t)\La^{\al}(s)^*f(x)\|_{L^{p'}}\leq C_n^{2\al}
(\suml_k 2^{k(n-1)/p}\supl_{U\in SU(n),x}
\norm{A_k(t,x+Uz)}{L^\infty_t L^{p}_{z_2, \ldots, z_n} L^1_{z_1}})^{2 \al}\times \\
& & \times 
 |t-s|^{-n/(2p)} \norm{f}{L^p}= C_n^\al \norm{A}{\tilde{Y}_1}^{2\al} \norm{f}{L^p}, 
\end{eqnarray*}
which suffices for \eqref{eq:906}. Thus, 
 \eqref{eq:90} and \eqref{eq:91} imply \eqref{eq:906}. 
\subsection{Energy estimates: Proof of \eqref{eq:90}} It is 
 more convenient to show 
$(\La^{\al})^*:L^2\to L^2$, which is equivalent to  \eqref{eq:90}. 
Clearly
$$
(\La^{\al})^*f(x,t)=\intl    
 e^{4\pi^2 i t|\xi|^2} e^{2\pi i \dpr{\xi}{x-y}}
\Om(\xi) f(y) \overline{\si^{\al}(t,y,\xi)} dy d\xi.
$$
Having in mind the specific form of $\si$, matters reduce to 
\begin{lemma}
\label{le:34}
Let $\al$ be an positive integer and $k_1,  k_2 \ldots,  k_\al$ 
are integers. Let 
$\{F^\mu_{k_\mu}\}$ be a collection of functions with 
$supp \widehat{F^\mu_{k_\mu}}\subset \{\xi:|\xi|\sim 2^{k_\mu}\}$. Then 
for the multilinear operator 
\begin{eqnarray*}
& &\Ve^{k_1, \ldots, k_\al}_{F^1, \ldots, F^\al} f(t,x)= 
\int 
e^{4\pi^2 i t|\xi|^2} e^{2\pi i \dpr{\xi}{x-y}}
\Omega(\xi)  d\xi f(y) \prod_{\mu=1}^\al
(\int F^\mu_{k_\mu}(t,y+z\xi/|\xi|) \chi(2^{2k_\mu}z) dz )dy,
\end{eqnarray*}
there is the estimate  
\begin{equation}
\label{a:11}
\norm{\Ve^{k_1, \ldots, k_\al}_{F^1, \ldots, F^\al} 
 f(t, \cdot)}{L^2_x}\leq C_n^\al \norm{f}{L^2}\prod_{\mu=1}^\al 
\supl_{x,\theta} \int |F^\mu_{k_\mu}(t, x+z\theta)| dz.
\end{equation}
\end{lemma}
The Lemma is applied to $(\La^\al)^*$ in an obvious way. That is, 
 write \\ 
$\si=\suml_{k\leq -2} \int A_k(t,x+z\xi/|\xi|) \chi(2^{2k}z) dz$ and 
\begin{eqnarray*}
& &
(\La^\al)^* f(t,x)= \\
& & =\suml_{k_1, \ldots, k_\al} 
\int e^{4\pi^2 i t|\xi|^2} e^{2\pi i \dpr{\xi}{x-y}}
\Omega(\xi)  f(y) \prod_{\mu=1}^\al
(\int A_{k_\mu}(t,y+z\xi/|\xi|) \chi(2^{2k_\mu}z) dz )dy d\xi =\\
& &= \suml_{k_1, \ldots, k_\al} 
 \Ve_{A, \ldots, A}^{k_1, \ldots k_\al} f(t,x)
\end{eqnarray*}
It follows that 
\begin{eqnarray*}
& & \norm{\La^\al f}{L^2}\leq \suml_{k_1, \ldots, k_\al} 
\norm{\Ve_{A, \ldots, A}^{k_1, \ldots k_\al} f}{L^2}\leq C_n^\al \norm{f}{L^2}
(\suml_{k} \supl_{x, \theta} \int |A_k(t,x+z\theta)|dz)^\al,
\end{eqnarray*}
as claimed. 
\begin{proof}(Lemma \ref{le:34})
Let us first present to the proof in the case $\al=1$, 
since the proof in the general case follows similar scheme, 
with somewhat cumbersome notations.
 
The basic idea is to ``pretend'' 
that $\int F_{k_\mu}(t,y+z\xi/|\xi|) \chi(2^{2k_\mu}z) dz$ is independent 
of $\xi$. 
We  show that this is almost true, 
modulo the angular  decomposition, that we have described earlier. 

For the given $k$,  introduce the partition of unity 
$$
\suml_j \psi_{j, -k}(2^{-k}(\xi/|\xi|-\theta_j^{-k}))=1
$$
Next, expand $F_k(y+z\xi/|\xi|)$ around $y+z\theta^{-k}_j$. We get 
\begin{eqnarray*}
& &
F_k(t,y+z\xi/|\xi|) \chi(2^{2k}z) 
\psi_{j, -k}(2^{-k}(\xi/|\xi|-\theta_j^{-k}))= \\
& &=
\suml_{\ga\geq 0} (\ga!)^{-1} 
F^\ga_k(t, y+z\theta_j^{-k}) \chi^\ga(2^{2k}z) 
\psi^\ga_{j, -k}(2^{-k}(\xi/|\xi|-\theta_j^{-k}))
\end{eqnarray*}
We drop the $\ga$'s, since in the end, we always 
add up with the help of $(\ga!)^{-1}$. We have 
\begin{equation}
\label{eq1}
\Ve f(t,x)=
e^{- i t \De} \suml_{j} 
P^{j} 
[f(\cdot)\int F_k(t,\cdot+z\theta_j^{-k}) \chi(2^{2k}z)dz ],
\end{equation}
where $\widehat{P^{j}g}(\xi)= \hat{g}(\xi) 
\psi_{j, -k}(2^{-k}(\xi/|\xi|-
\theta_{j}^{-k})) $. Note that the function \\
$y\to \int F_k(y+z\theta_j^{-k}) \chi(2^{2k}z)dz$ has a  
Fourier transform supported in $\{\xi: |\xi|\sim 2^k\}$. It follows
\begin{equation}
\label{eq2}
P^{j} 
[f(\cdot)\int F_k(t,\cdot+z\theta_j^{-k}) \chi(2^{2k}z)dz ]=P^{j} 
[(\int F_k(t,\cdot+z\theta_j^{-k}) \chi(2^{2k}z)dz )\tilde{P}^j(f)],
\end{equation}
for some $\tilde{P}^j: \tilde{P}^jP^j=P^j$ and $\tilde{P}^j$ 
has  multiplier with support inside $|\xi/|\xi|-
\theta_{j}^{-k}|\leq 2^{2k+2}$. 

This implies 
\begin{eqnarray*}
& & \norm{\Ve f(t,x)}{L^2_x}\lesssim \left(\suml_j \norm{P^{j} 
[f(\cdot)\int F_k(t,\cdot+z\theta_j^{-k}) 
\chi(2^{2k}z)dz ]}{L^2}^2\right)^{1/2}= \\
& &= \left(\suml_j \norm{P^j[(\int F_k(t,\cdot+z\theta_j^{-k}) 
\chi(2^{2k}z)dz) \tilde{P}^j f]}{L^2}^2\right)^{1/2} \lesssim \\
& &\lesssim  \left(\suml_j \norm{(\int F_k(t,\cdot+z\theta_j^{-k}) 
\chi(2^{2k}z)dz) \tilde{P}^j f}{L^2}^2\right)^{1/2} \lesssim\\
& &\lesssim \left(\suml_j \norm{\int F_k(t,\cdot+z\theta_j^{-k}) 
\chi(2^{2k}z)dz}{L^\infty}^2\norm{\tilde{P}^j f}{L^2}^2\right)^{1/2}
\lesssim \\
& &\leq C_n  \supl_{y, \theta} |\int F_k(t,y+z\theta) 
\chi(2^{2k}z)dz| \norm{f}{L^2}.
\end{eqnarray*}
Note that the previous calculation requires slightly augmenting $P^j$ to 
$\tilde{P}^j$ for an additional constant $C_n$ comming 
on the account of the extra overlap of the supports of different 
$\tilde{P}^j$. 

Unfortunately, we have to be extra careful 
for the case of $\al>1$, because the constants (following this argument) 
are estimated by $\al^\al$, whereas we need (and can manage) constants of 
magnitude  $C_n^\al$. This problems occurs, when too many of the 
frequencies $k_1, \ldots , k_\al$ are equal. 

Thus, we start our considerations for $\al>1$ 
by ordering these frequencies. Without loss of generality 
let us assume 
$k_1=\ldots=k_{s_1}<k_{s_1+1}=\ldots=k_{s_2}<\ldots 
k_{s_{m-1}+1}=\ldots=k_{s_m}=k_\al$, where we have 
set $s_0=0$ for convenience. 

 Since for every $1\leq r\leq m$,  and 
$\mu\in[s_{r-1}+1,s_r]$
$$
\suml_{j_r} \psi_{j_r, -k_{\mu}}(2^{-k_{\mu}}(\xi/|\xi|-
\theta_{j_r}^{-k_{\mu}}))=1,
$$
we expand  $F^\mu_{k_{\mu}}(t, y+z\xi/|\xi|)$ around  
$y+z \theta_{j_r}^{-k_{\mu}}$. Just as in the case $\al=1$, this 
 allows us to write 
\begin{equation}
\label{eqn12}
\Ve^{k_1, \ldots, k_\al}_{F_1, \ldots, F_\al} f(t,x)= e^{- i t \De}
\suml_{j_1, \ldots, j_m}
\suml_{\ga_1, \ldots 
\ga_\al} 
\f{ P^{j_1}\ldots P^{j_m} [
f(\cdot)\prod_{\mu=1}^\al g_\mu^\ga(t,\cdot)]}{(\ga_1! \ldots 
\ga_\al! )},
\end{equation}
where for $\mu\in[s_{r-1}+1,s_r]$
\begin{eqnarray*}
& & \widehat{P^{j_r}g}(\xi)= \hat{g}(\xi) 
\psi_{j_r, -k_{\mu}}^{\ga}(2^{-k_\mu}(\xi/|\xi|-
\theta_{j_r}^{-k_\mu})) ,\\
& & g_\mu^\ga(t,y)=\int F^{\mu, \ga}_{k_\mu}(t, y+z\theta_{j_r}^{-k_\mu})
\chi(2^{2k_\mu}z)dz.
\end{eqnarray*}
As always we drop the $\ga$'s and concentrate on the case 
$\ga_1=\ldots=\ga_\al=0$.  

Next,  observe that since $e^{-i t \De}$ is an isometry on $L^2$, 
we can dispose of it immediately. We have 
$$
\norm{\Ve^{k_1, \ldots, k_\al}_{F_1, \ldots, F_\al} f}{L^2_x}^2 
\leq \suml_{j_1, \ldots, j_m} \norm{P^{j_1}\ldots 
P^{j_m} [f(\cdot)\prod_{\mu=1}^\al g_\mu^\ga(t,\cdot)]}{L^2}^2. 
$$
For technical reasons, it is more convenient to 
replace $P^{j_r}$ by ``rough'' versions of the same. Namely, 
 introduce the Fourier restriction operators $Q^{j_r}$, 
which act  via \\
$\widehat{Q^{j_r} g}(\xi)= \hat{g}(\xi) 1_{\{|\xi/|\xi|-
\theta_{j_r}^{k_{s_r}}|\leq 2^{k_{s_r}}\}}$. More generally, for any $a>0$, 
$$
\widehat{Q^{j_r}_a g}(\xi)= \hat{g}(\xi) 1_{\{\xi: |\xi/|\xi|-\theta_{j_r}^{k_{s_r}}|\leq a\}}
$$
Since the Fourier supports of the multipliers of $P^{j_r}$ are in 
$\{\xi: |\xi/|\xi|-\theta_{j_r}^{k_{s_r}}|\leq 2^{k_{s_r}}\}$,  we have 
\begin{equation}
\label{eqn:7}
\norm{\Ve^{k_1, \ldots, k_\al}_{F_1, \ldots, F_\al} f}{L^2_x}^2 
\leq  \suml_{j_1, \ldots, j_m} \norm{Q^{j_1}\ldots 
Q^{j_m} [f(\cdot)\prod_{\mu=1}^\al g_\mu(t,\cdot)]}{L^2}^2. 
\end{equation}
Fix $j_1$. Clearly the summation in any $j_r$ runs only on 
$j_r: |\theta_{j_r}^{-k_{s_r}}-\theta_{j_1}^{-k_{s_1}}|\leq 
2^{k_{s_1}}+2^{k_{s_r}}$, since  otherwise $Q^{j_1} Q^{j_r}=0$. Thus, 
\begin{equation}
\label{eqn:1}
\norm{\Ve^{k_1, \ldots, k_\al}_{F_1, \ldots, F_\al} f}{L^2_x}^2 
\leq \suml_{j_1} 
\suml_{j_2,  \ldots, j_m:|\theta_{j_r}^{-k_{s_r}}-\theta_{j_1}^{-k_{s_1}}|\leq 
2^{k_{s_1}}+2^{k_{s_r}}}
\norm{Q^{j_1} [f(\cdot)\prod_{\mu=1}^\al g_\mu(t,\cdot)]}{L^2}^2. 
\end{equation}
Note that since $supp\ \hat{g_\mu}\subset \{\xi: |\xi|\leq 2^{k_\mu+1}\}$, 
we have that $supp\ \widehat{g_1\ldots g_{s_1}}
\subset \{\xi: |\xi|\leq 2 s_1 2^{k_{s_1}}\}$. It follows that 
\begin{equation}
\label{eqn14}
\begin{array}{l}
Q^{j_1} [f(\cdot)\prod_{\mu=1}^\al g_\mu(t,\cdot)]= 
Q^{j_1} [(g_1\ldots g_{s_1}) 
(f(\cdot)\prod_{\mu=s_1+1}^\al g_\mu(t,\cdot))]=\\
= Q^{j_1} [(g_1\ldots g_{s_1}) 
Q^{j_1}_{2^{k_{s_1}}+2s_1 
2^{k_{s_1}}}(f(\cdot)\prod_{\mu=s_1+1}^\al g_\mu(t,\cdot))]
\end{array}
\end{equation}
Plugging this back in \eqref{eqn:1} and taking into account 
$\norm{Q^{j_1}}{L^2\to L^2}=1$ yields 
\begin{eqnarray*}
& &
\norm{\Ve^{k_1, \ldots, k_\al}_{F_1, \ldots, F_\al} f}{L^2_x}^2 
\leq  \\
& &\leq \suml_{j_1} 
\suml_{j_2,  \ldots, j_m:|\theta_{j_r}^{-k_{s_r}}-\theta_{j_1}^{-k_{s_1}}|\leq 
2^{k_{s_1}}+2^{k_{s_r}}}
 \norm{g_1 \ldots g_{s_1} Q^{j_1}_{2^{k_{s_1}}+2s_1 2^{k_{s_1}}}
(f(\cdot)\prod_{\mu=s_1+1}^\al g_\mu(t,\cdot))}{L^2}^2\leq \\
& &
\leq \suml_{j_1} 
\suml_{j_2,  \ldots, j_m:|\theta_{j_r}^{-k_{s_r}}-\theta_{j_1}^{-k_{s_1}}|\leq 
2^{k_{s_1}}+2^{k_{s_r}}}
\norm{g_1 \ldots g_{s_1}}{L^\infty}^2
 \norm{Q^{j_1}_{2^{k_{s_1}}+2s_1 2^{k_{s_1}}}
(f(\cdot)\prod_{\mu=s_1+1}^\al g_\mu(t,\cdot))}{L^2}^2\leq \\
& &\leq \prod_{\mu=1}^{s_1} 
\sup_{y, \theta} |\int F^\mu_{k_\mu}(t,y+z\theta) \chi(2^{2k_\mu}z) dz|^2
\times \\
& &\times \suml_{j_1} 
\suml_{j_2,  \ldots, j_m:|\theta_{j_r}^{-k_{s_r}}-\theta_{j_1}^{-k_{s_1}}|\leq 
2^{k_{s_1}}+2^{k_{s_r}}}
 \norm{Q^{j_1}_{2^{k_{s_1}}+2s_1 2^{k_{s_1}}}
(f(\cdot)\prod_{\mu=s_1+1}^\al g_\mu(t,\cdot))}{L^2}^2. 
\end{eqnarray*}
It is now time to reintroduce the $Q^{j_r}$ multipliers. 
Since $|\theta_{j_r}^{-k_{s_r}}-\theta_{j_1}^{-k_{s_1}}|\leq 
2^{k_{s_1}}+2^{k_{s_r}}$, we conclude that for every $r: 1\leq r\leq m$,   
$Q^{j_1}_{2^{k_{s_1}}+2s_1 2^{k_{s_1}}}=Q^{j_1}_{2^{k_{s_1}}+2s_1 2^{k_{s_1}}} 
Q^{j_r}_{2^{k_{s_r}}+(2s_1+2)2^{k_{s_1}}}$ and we have 
\begin{eqnarray*}
& & \suml_{j_1} 
\suml_{j_2,  \ldots, j_m:|\theta_{j_r}^{-k_{s_r}}-\theta_{j_1}^{-k_{s_1}}|\leq 
2^{k_{s_1}}+2^{k_{s_r}}}
 \norm{Q^{j_1}_{2^{k_{s_1}}+2s_1 2^{k_{s_1}}}
(f(\cdot)\prod_{\mu=s_1+1}^\al g_\mu(t,\cdot))}{L^2}^2\leq \\
& &\leq \suml_{j_2,  \ldots, j_m}\suml_{j_1}
\norm{ Q^{j_1}_{2^{k_{s_1}}+2s_1 2^{k_{s_1}}} 
Q^{j_2}_{2^{k_{s_2}}+(2s_1+2)2^{k_{s_1}}}\ldots 
Q^{j_m}_{2^{k_{s_m}}+(2s_1+2)2^{k_{s_1}}} (f(\cdot)
\prod_{\mu=s_1+1}^\al g_\mu(t,\cdot))}{L^2}^2
\end{eqnarray*}
Note that
\begin{equation}
\label{eqn17}
\suml_{j_1} \norm{Q^{j_1}_{2^{k_{s_1}}+2s_1 2^{k_{s_1}}} g}{L^2}^2\leq C_n 
\left(\f{2^{k_{s_1}}+2s_1 2^{k_{s_1}}}{2^{k_{s_1}}}\right)^n 
\norm{g}{L^2}^2=C_n(1+2s_1)^n \norm{g}{L^2}^2 ,
\end{equation}
because of the extra overlap created by passing from $Q^{j_1}$ to 
$Q^{j_1}_{2^{k_{s_1}}+2s_1 2^{k_{s_1}}}$. It follows that 
\begin{eqnarray*}
& &\norm{\Ve^{k_1, \ldots, k_\al}_{F_1, \ldots, F_\al} f}{L^2_x}^2  \leq 
C_n \prod_{\mu=1}^{s_1} 
\sup_{y, \theta} |\int F^\mu_{k_\mu}(t,y+z\theta) 
\chi(2^{2k_\mu}z) dz|^2 (1+2s_1)^n \times \\
& & \times \suml_{j_2,  \ldots, j_m} 
\norm{Q^{j_2}_{2^{k_{s_2}}+(2s_1+2)2^{k_{s_1}}}\ldots 
Q^{j_m}_{2^{k_{s_m}}+(2s_1+2)2^{k_{s_1}}} (f(\cdot)
\prod_{\mu=s_1+1}^\al g_\mu(t,\cdot))}{L^2}^2.
\end{eqnarray*} 
This is very similar to \eqref{eqn:7}, except that the sum in $j_1$ is taken care of and 
$Q^{j_2}\ldots Q^{j_m}$ gets 
replaced by the slightly larger  
$Q^{j_2}_{2^{k_{s_2}}+(2s_1+2)2^{k_{s_1}}}\ldots 
Q^{j_m}_{2^{k_{s_m}}+(2s_1+2)2^{k_{s_1}}}$. Continuing in this fashion yields 
the estimate 
\begin{eqnarray*}
& &\norm{\Ve^{k_1, \ldots, k_\al}_{F_1, \ldots, F_\al} f}{L^2_x}^2  \leq 
C_n^m \prod_{\mu=1}^{\al} 
\sup_{y, \theta} |\int F^\mu_{k_\mu}(t,y+z\theta) 
\chi(2^{2k_\mu}z) dz|^2 \norm{f}{L^2}^2 (1+2s_1)^n  \times \\
& & \times  \prod_{r=2}^m 
\left(\f{2^{k_{s_r}}+(2s_1+2) 2^{k_{s_1}}+ (2(s_2-s_1)+2) 2^{k_{s_2}}+\ldots +
(2(s_r-s_{r-1})+2)2^{k_{s_r}}}{2^{k_{s_r}}}\right)^n.
\end{eqnarray*} 
Note that for every $1\leq i \leq r$, there is 
 $k_{s_r}-k_{s_i}\geq r-i$, since 
there are $r-i$ strict inequalities in the chain  $k_{s_i}<\ldots<k_{s_r}$. 
We therefore need an estimate for 
$$
G=\prod_{r=1}^m (1+\suml_{i=1}^r (2(s_i-s_{i-1})+2)2^{i-r})^n
$$
For $\ln(G)$, note that 
 since $\ln(1+x)\leq x$ and $m\leq \al$,  we can  estimate 
$$
\ln(G)\leq n \suml_{r=1}^m \suml_{i=1}^r (2(s_i-s_{i-1})+2)2^{i-r}
\leq 4n (m+\suml_{i=1}^{m} (s_i - s_{i-1}))=4n(m+\al)\leq  8n\al. 
$$
It follows that $G\leq e^{8 n\al}$ and 
$$
\norm{\Ve^{k_1, \ldots, k_\al}_{F_1, \ldots, F_\al} f}{L^2_x}^2 \leq  C_n^\al 
\prod_{\mu=1}^{\al} 
\sup_{y, \theta} |\int F^\mu_{k_\mu}(t,y+z\theta) 
\chi(2^{2k_\mu}z) dz|^2 \norm{f}{L^2}^2
$$

\end{proof}

\subsection{Dispersive estimates: Proof of \eqref{eq:91}}
For the dispersive estimate,  write 
\begin{eqnarray*}
& &\La^{\al}(t)\La^{\al}(s)^*f = \intl   \si^{\al}(t,x,\xi)
 \overline{\si^{\al}(s,y,\xi)}
 e^{-4\pi^2 i (t-s)|\xi|^2} e^{2\pi i \dpr{\xi}{x-y}}
\Om^2(\xi) f(y) dy d\xi
\end{eqnarray*}
Clearly for $\al=0$, we have 
\eqref{eq:91} by the decay  estimates for the free solution. 
Consider the case $\al=1$ for simplicity, the general case to 
be addressed momentarily.\\
Expand $\si(t,x,\xi)$ to get 
\begin{eqnarray*}
& &\La^{1}(t)\La^{1}(s)^*f= \suml_{k_1 \leq -2}
\suml_{l_1\leq -2k_1}
\suml_{\ga_1\geq 0} (\ga_1!)^{-1}\suml_{j_1}  
(\int A^{\ga_1}_k(t, x+z \theta^{l_1+k_1}_{j_1})\vp^{\ga_1}
(2^{-l_1}z)dz) \times \\
& &\times \int 
\psi_{j_1,l_1+k_1}^\ga(2^{l_1+k_1}(\xi/|\xi|-\theta_{j_1}^{l_1+k_1}))  
 e^{-4\pi^2 i (t-s)|\xi|^2} e^{2\pi i \dpr{\xi}{x-y}}
\Om^2(\xi) f(y) \overline{(\si)(s,y,\xi)} dy d\xi=\\
& &=\suml_{k_1\leq -2}
\suml_{l_1\leq -2 k_1}
\suml_{\ga_1\geq 0} (\ga_1!)^{-1}\suml_{j_1}  
(\int A^{\ga_1}_{k_1}(t, x+z \theta^{l_1+k_1}_{j_1})\vp^{\ga_1}
(2^{-l_1}z)dz)\Ga^1_{j_1, l_1+k_1} f,
\end{eqnarray*}
where\footnote{In the case $l_1\leq -k_1$ the summation in $j_1$ collapses to 
a single term and $\theta_{j_1}^{l_1+k_1}=e_1$ as pointed out in the previous section.}
$$
\Ga^\al_{j, l+k} f(x)= \intl   
 \overline{\si^{\al}(s,y,\xi)}
 e^{-4\pi^2 i (t-s)|\xi|^2} e^{2\pi i \dpr{\xi}{x-y}}
\psi_{j,l+k}^\ga(2^{l+k}(\xi/|\xi|-\theta_{j}^{l+k})) 
\Om^2(\xi) f(y) dy d\xi
$$
It is easy to see that
\begin{eqnarray*}
& &\norm{\La^{1}(t)\La^{1}(s)^*f}{L^\infty_x}\lesssim \\
& &\lesssim  \supl_x 
(\suml_{k_1\leq -2}
\suml_{l_1\leq -2 k_1} \suml_{\ga_1\geq 0} (\ga_1!)^{-1}\suml_{j_1}  
\int |A^{\ga_1}_{k_1}(t, x+z \theta^{l_1+k_1}_{j_1})|\vp^{\ga_1}
(2^{-l_1}z)dz )\times \\
& &\times \supl_{j_1, l_1, k_1} \norm{\Ga^1_{j_1, l_1+k_1} f}{L^\infty}.
\end{eqnarray*}
Bythe pointwise estimates of 
 Lemma \ref{le:pointwise} and $A\in Y_3$, the last expression is 
bounded by $C\norm{A}{Y_1}^2 |t-s|^{-n/2}$, provided one can show 
$\norm{\Ga^1_{j_1, l_1+k_1} f}{L^\infty}\leq C\norm{A}{Y_1}|t-s|^{-n/2}\norm{f}{L^1}$. 

\noindent 
More generally, it is easy to see that by 
iterating the argument above, we have 
\begin{eqnarray*}
& &\La^{\al}(t)\La^{\al}(s)^*f(x)= \\
& &=
\prod_{\mu=1}^\al (\suml_{k_\mu\leq -2}
\suml_{l_\mu\leq -2k_\mu} 
\suml_{\ga_\mu} (\ga_\mu!)^{-1} 
\suml_{j_\mu}  
\int A^{\ga_\mu}_{k_\mu}(t, x+z \theta^{l_\mu+k_\mu}_{j_\mu})\vp^{\ga_\mu}
(2^{-l_\mu}z)dz )
\Ga^\al_{j_1,l_1+k_1;\ldots; j_\al, l_\al+k_\al} f,
\end{eqnarray*}
where 
\begin{eqnarray*}
& & \Ga^\al_{j_1,l_1+k_1;\ldots; j_\al, l_\al+k_\al} f(x)  = 
\intl   
 \overline{\si^{\al}(s,y,\xi)}
 e^{-4\pi^2 i (t-s)|\xi|^2} e^{2\pi i \dpr{\xi}{x-y}}
 \times \\ 
& & \times  \prod_{\mu=1}^\al \psi_{j_\mu,l_\mu+k_\mu}^{\ga_\mu}
(2^{l_\mu+k_\mu}(\xi/|\xi|-\theta_{j_\mu}^{l_\mu+k_\mu}))
\Om^2(\xi) f(y) dy d\xi.
\end{eqnarray*}
Thus, 
\begin{eqnarray*}
& &\norm{\La^{\al}(t)\La^{\al}(s)^*f}{L^\infty}\leq \\
& & \leq  \prod_{\mu=1}^\al \supl_x[
\suml_{k_\mu\leq -2}
\suml_{l_\mu\leq -2k_\mu} 
\suml_{\ga_\mu\geq 0} (\ga_\mu)^{-1}  \suml_{j_\mu} 
\int |A^{\ga_\mu}_{k_\mu}(t, x+z \theta^{l_\mu+k_\mu}_{j_\mu})|\vp^{\ga_\mu}
(2^{-l_\mu}z) dz )]\times \\
& & \times \supl_{j_1,l_1,k_1;\ldots; j_\al, l_\al,k_\al}
\norm{\Ga^\al_{j_1,l_1+k_1;\ldots; j_\al, l_\al+k_\al} f}{L^\infty}.
\end{eqnarray*}
By Lemma \ref{le:pointwise}, we have 
\begin{eqnarray*}
& &\norm{\La^{\al}(t)\La^{\al}(s)^*f}{L^\infty} 
\leq C_n^{2 \al} \norm{A}{Y_1}^{2\al} |t-s|^{-n/2}\norm{f}{L^1}, 
\end{eqnarray*}
provided one can show 
\begin{equation}
\label{eq:340}
 \supl_{j_1,l_1,k_1;\ldots; j_\al, l_\al,k_\al}
\norm{\Ga^\al_{j_1,l_1+k_1;\ldots; j_\al, l_\al+k_\al} f}{L^\infty}
\leq C_n^\al \norm{A}{Y_1}^\al |t-s|^{-n/2}\norm{f}{L^1}
\end{equation}
Note that $\Ga^\al_{j_1, l_1+k_1; \ldots ; j_\al, l_\al+k_\al} f$ 
looks somewhat like $\La(s)^* f$, with the important difference 
that it has the  ``multiplier'' 
$\prod_{\mu=1}^\al \psi_{j_\mu,l_\mu+k_\mu}(2^{l_\mu+k_\mu}(\xi/|\xi|-
\theta_{j_\mu}^{l_\mu+k_\mu}))$ in its definition. \\
Dualizing \eqref{eq:340} leads us to showing that 
\begin{eqnarray*}
& & [(\Ga^\al_{j_1, l_1+k_1; \ldots; 
j_\al, l_\al+k_\al})^* f] (x) = \intl   
 \si^{\al}(s,x,\xi)
 e^{ 4\pi^2 i (t-s)|\xi|^2} e^{2\pi i \dpr{\xi}{x-y}}
 \times \\ 
& & \times \prod_{\mu=1}^\al \psi_{j_\mu,l_\mu+k_\mu}^{\ga_\mu}
(2^{l_\mu+k_\mu}(\xi/|\xi|-\theta_{j_\mu}^{l_\mu+k_\mu}))
\Om^2(\xi) f(y) dy d\xi.
\end{eqnarray*}
maps $L^1\to L^\infty$ with norm no bigger than 
$C_n^\al \norm{A}{Y_1}^\al |t-s|^{-n/2}$. Expand again to get 
\begin{eqnarray*}
& & [(\Ga^\al_{j_1, l_1+k_1; \ldots; 
j_\al, l_\al+k_\al})^* f](x)= \\
& &=
\prod_{\nu=1}^\al (\suml_{\tilde{k}_\nu \leq -2}
\suml_{\tilde{l}_\nu\leq -2\tilde{k}_\nu} 
\suml_{\ga_\nu} (\ga_\nu!)^{-1} \suml_{j_\nu} 
\int A^{\ga_\nu}_{\tilde{k}_\nu}(s, 
x+z \theta^{\tilde{l}_\nu+\tilde{k}_\nu}_{j_\nu})\vp^{\ga_\nu}
(2^{-\tilde{l}_\nu}z)dz )\times \\
& &\times \intl   
 e^{ 4\pi^2 i (t-s)|\xi|^2} e^{2\pi i \dpr{\xi}{x-y}}
  \prod_{\nu=1}^\al \psi_{j_\nu,\tilde{l}_\nu+\tilde{k}_\nu}^{\ga_\nu}
(2^{\tilde{l}_\nu+\tilde{k}_\nu}(\xi/|\xi|-
\theta_{j_\nu}^{\tilde{l}_\nu+\tilde{k}_\nu})) \times \\
& &\times 
\prod_{\mu=1}^\al \psi_{j_\mu,l_\mu+k_\mu}^{\ga_\mu}
(2^{l_\mu+k_\mu}(\xi/|\xi|-\theta_{j_\mu}^{l_\mu+k_\mu}))
\Om^2(\xi) f(y) dy d\xi.
\end{eqnarray*}
By the Krein-Milman 
theorem, $span[\de_b]$ is  $w^*$ dense in 
$M(\rn)\supset L^1(\rn)$. This, together with 
Lemma \ref{le:pointwise} allows us to estimate 
\begin{eqnarray*}
& & \norm{(\Ga^\al_{j_1, l_1+k_1; \ldots; 
j_\al, l_\al+k_\al})^* f}{L^\infty_x}\leq  C_n^\al (\norm{A}{Y_1})^\al 
\times \\
& & \supl_{x,b\in \rn; m_1, \ldots, m_{2\al}\geq 0;\tilde{\theta}_1, 
\ldots, \tilde{\theta}_{2\al} \in \sn}|\intl   
 e^{ 4\pi^2 i (t-s)|\xi|^2} e^{2\pi i \dpr{\xi}{x-b}}
\prod_{\mu=1}^{2\al} \psi_{j_\mu,m_\mu}
(2^{m_\mu}(\xi/|\xi|-\tilde{\theta}_{\mu}))
\Om^2(\xi)d\xi|.
\end{eqnarray*}
The oscillatory integral above is bounded by $C_n^\al |t-s|^{-n/2}$ by 
  Lemma \ref{le:osc} in the Appendix, and the required 
dispersive estimates hold. 
\end{proof}

\section{The parametrix almost satisfies the equation}
\label{sec:7}
In this section, we  show
 that the parametrix satisfies \eqref{eq:904}.
We have several types of terms that arise according to 
 \eqref{eq:501}. 

First, we take on the terms  
$$
\int (i \p_t (\si^0+\si^1)+\De \si^0+2\pi i \dpr{\nabla\si^1}{\xi})
e^{i \si(t,x,\xi)} e^{-4\pi^2 i t|\xi|^2} e^{2\pi i \dpr{\xi}{x}}
\Omega(\xi) \widehat{f}(\xi) d\xi.
$$
These are {all terms linear in either $\si^0$ or $\si^1$} 
with either a time derivative or two spatial derivative acting on them
(recall from the expression 
\eqref{e:1} that $\dpr{\nabla\si^1}{\xi})$ is also of this form).
\subsection{Main result}
\begin{lemma}
\label{le:10}
 Let $\Ve^{k_1, \ldots,k_\al}_{F^1, \ldots, F^\al}$
 be as in Lemma \ref{le:34}, i.e. 
\begin{eqnarray*}
& &\Ve^{k_1, \ldots,k_\al}_{F^1, \ldots, F^\al}f(t,x)= 
\int \prod_{\mu=1}^\al 
(\int F^\mu_{k_\mu}(t,x+z\xi/|\xi|) \chi(2^{2k_\mu}z) dz )
e^{-4\pi^2 i t|\xi|^2} e^{2\pi i \dpr{\xi}{x}}
\Omega(\xi) \widehat{f}(\xi) d\xi.
\end{eqnarray*}
Then  there exists a constant $C_n$, so that  
\begin{eqnarray*}
& &
\norm{\Ve^{k_1, \ldots,k_\al}_{F^1, \ldots, F^\al}f}{L^1_tL^2_x}\leq 
C_n^\al \norm{f}{L^2} \times \\
& &\times 
\prod_{\mu=1}^{\al-1}
(2^{k_\mu(n-1)/2}\supl_{U\in SU(n)} \supl_t \supl_x 
\| F^\mu_{k_\mu}(t,x+ Uz)\|_{L^2_{z_2, \ldots, z_n}L^1_{z_1}}+
\\
& &+ 2^{k_\mu(n+3)/2}\supl_{U\in SU(n)} \| \supl_x 
\| F^\mu_{k_\mu}(t,x+Uz)\|_{L^2_{z_2, \ldots, z_n}L^1_{z_1}}\|_{L^1_t})
 \times \\
& &(2^{k_\mu(n-1)/2}\supl_{U\in SU(n)} \| \supl_x 
\| F^\al_{k_\al}(t,x+Uz)\|_{L^2_{z_2, \ldots, z_n}L^1_{z_1}}\|_{L^1_t}+ \\
& &+
2^{k_\mu(n-5)/2}\supl_{U\in SU(n)} \supl_t \supl_x 
\| F^\al_{k_\al}(t,x+ Uz)\|_{L^2_{z_2, \ldots, z_n}L^1_{z_1}}).
\end{eqnarray*}
In a different form, 
\begin{eqnarray*}
& &
\norm{\Ve^{k_1, \ldots,k_\al}_{F^1, \ldots, F^\al}f}{L^1_tL^2_x}\leq 
C_n^\al \norm{f}{L^2} \times \\
& &\times 
\prod_{\mu=1}^{\al-1}
(2^{k_\mu(n-1)/2}\supl_{U\in SU(n)} \supl_x 
\| F^\mu_{k_\mu}(t,x+ Uz)\|_{L^\infty_t L^2_{z_2, \ldots, z_n}L^1_{z_1}}+ \\
& &
+ 2^{k_\mu(n+3)/2}\supl_{U\in SU(n), t\to x(t)}
\| F^\mu_{k_\mu}(t,x(t)+Uz)\|_{L^1_t L^2_{z_2, \ldots, z_n}L^1_{z_1}})
 \times \\
& &(2^{k_\mu(n-1)/2}\supl_{U\in SU(n), t\to x(t)} 
\| F^\al_{k_\al}(t,x(t)+Uz)\|_{L^1_t 
L^2_{z_2, \ldots, z_n}L^1_{z_1}}+ \\
& &+
2^{k_\mu(n-5)/2}\supl_{U\in SU(n)}\supl_x 
\| F^\al_{k_\al}(t,x+ Uz)\|_{L^\infty_t L^2_{z_2, \ldots, z_n}L^1_{z_1}}),
\end{eqnarray*}
where $\supl_{t\to x(t)}$ is taken over all measurable
 functions $x(\cdot): \rone\to \rn$. 
\end{lemma}
\subsection{Application of Lemma \ref{le:10} to various terms of $\norm{\cl v}{L^1L^2}$.}
\noindent Assuming the validity of Lemma \ref{le:10}, one can easily handle the first type of terms. 

Expand  $e^{i \si(t,x,\xi)}$ in powers of $\si$  
as in the previous section yields
\begin{eqnarray*}
& &
\int (i \p_t (\si^0+\si^1)+\De \si^0+2\pi 
i\dpr{\nabla\si^1}{\xi})
e^{i \si(t,x,\xi)} e^{-4\pi^2 i t|\xi|^2} e^{2\pi i \dpr{\xi}{x}}
\Omega(\xi) \widehat{f}(\xi) d\xi =\\
& &=\suml_{\al\geq 1} \f{i^{\al-1}}{(\al-1)!} 
\int \si(t,x,\xi)^{\al-1} (i \p_t (\si^0+\si^1)+\De \si^0+2\pi 
i\dpr{\nabla\si^1}{\xi})\times \\
& &\times 
 e^{-4\pi^2 i t|\xi|^2} e^{2\pi i \dpr{\xi}{x}}
\Omega(\xi) \widehat{f}(\xi) d\xi =\suml_{\al\geq 1} \f{i^{\al-1}}{(\al-1)!}
 \suml_{k_1, \ldots, k_\al}
 \\
& &
=\int \prod_{\mu=1}^{\al-1} 
(\int A_{k_\mu}(t, x+z\xi/|\xi|) \chi(2^{2k_\mu}z) dz)
(\int \p^2 A_{k_\al}(t, x+z\xi/|\xi|) \chi(2^{2k_\al}z) dz)\times \\
& &\times 
e^{-4\pi^2 i t|\xi|^2} e^{2\pi i \dpr{\xi}{x}}
\Omega(\xi) \widehat{f}(\xi) d\xi =\\
& &= \suml_{\al\geq 1} \f{i^{\al-1}}{(\al-1)!}
 \suml_{k_1, \ldots, k_\al} \Ve^{k_1, \ldots, k_\al}_{A, \ldots,A, \p^2 A},
\end{eqnarray*}
where we have denoted by 
$\int \p^2 A_{k_\al} (t, x+z\xi/|\xi|) \chi(2^{2k}z) dz$ 
all the terms 
$P_{k_\al}(i \p_t (\si^0+\si^1)+\De \si^0+2\pi 
i\dpr{\nabla\si^1}{\xi})$. \\
Thus, an application of Lemma \ref{le:10} yields  
\begin{eqnarray*}
& &\norm{\int (i \p_t (\si^0+\si^1)+\De \si^0+2\pi 
i\dpr{\nabla\si^1}{\xi})
e^{i \si(t,x,\xi)} e^{-4\pi^2 i t|\xi|^2} e^{2\pi i \dpr{\xi}{x}}
\Omega(\xi) \widehat{f}(\xi) d\xi}{L^1 L^2} \\
& & \leq C_n \norm{f}{L^2} 
exp(C_n(\suml_k 2^{k(n-1)/2}\supl_{U\in SU(n)} \supl_x 
\| A_{k}(t,x+ Uz)\|_{L^\infty_t L^2_{z_2, \ldots, z_n}L^1_{z_1}}+ \\
& &
+ 2^{k(n+3)/2}\supl_{U\in SU(n), t\to x(t)}
\| A_{k}(t,x(t)+Uz)\|_{L^1_t L^2_{z_2, \ldots, z_n}L^1_{z_1}}))\times \\
& &\times \suml_k (2^{k(n-1)/2}\supl_{U\in SU(n), x(t)} 
\| (|\p^2 A_k|+|\p_t A_k|) (t,x(t)+Uz)\|_{L^1_t L^2_{z_2, 
\ldots, z_n}L^1_{z_1}}+ \\
& &+
2^{k(n-5)/2}\supl_{U\in SU(n), x(t)}
\| (|\p^2 A_k|+|\p_t A_k|) (t,x(t) + Uz)\|_{L^\infty_t 
L^2_{z_2, \ldots, z_n}L^1_{z_1}}). 
\end{eqnarray*}

Next, we show how to use Lemma \ref{le:10} to control terms in the form 
$$
\int  (\nabla\si)^2
e^{i \si(t,x,\xi)} e^{-4\pi^2 i t|\xi|^2} e^{2\pi i \dpr{\xi}{x}}
\Omega(\xi) \widehat{f}(\xi) d\xi.
$$
in $L^1L^2$ norm. We expand the same way as before 
$e^{i \si}=\sum_{\al\geq 2}i^{\al-2}/(\al-2)! \si^{\al-2}$ to get 
\begin{eqnarray*}
& & \int  (\nabla\si)^2
e^{i \si(t,x,\xi)} e^{-4\pi^2 i t|\xi|^2} e^{2\pi i \dpr{\xi}{x}}
\Omega(\xi) \widehat{f}(\xi) d\xi = \\
& &= \sum_{\al\geq 2} \f{i^{\al-2}}{(\al-2)!}
 \suml_{k_1, \ldots, k_\al} \int 
\prod_{\mu=1}^{\al-2} 
(\int A_{k_\mu}(t, x+z\xi/|\xi|) \chi(2^{2k_\mu}z) dz)\times \\
& &\times 
(\int \p A_{k_{\al-1}}(t, x+z\xi/|\xi|) \chi(2^{2k_{\al-1}}z) dz)
(\int \p A_{k_{\al}}(t, x+z\xi/|\xi|) \chi(2^{2k_\al}z) dz)\times \\
& &\times 
e^{-4\pi^2 i t|\xi|^2} e^{2\pi i \dpr{\xi}{x}}
\Omega(\xi) \widehat{f}(\xi) d\xi= \sum_{\al\geq 2} \f{i^{\al-2}}{(\al-2)!}
 \suml_{k_1, \ldots, k_\al} \Ve^{k_1, 
\ldots, k_\al}_{A, \ldots, A, \p A, \p A}.
\end{eqnarray*}
By the symmetry in the last entries $\p A_{k_{\al-1}}$ and $\p A_{k_\al}$, 
we will without loss of generality assume $k_{\al-1}\leq k_\al$. 
Applying Lemma \ref{le:10} yields 
\begin{eqnarray*}
& & \norm{\int  (\nabla\si)^2
e^{i \si(t,x,\xi)} e^{-4\pi^2 i t|\xi|^2} e^{2\pi i \dpr{\xi}{x}}
\Omega(\xi) \widehat{f}(\xi) d\xi}{L^1L^2} \leq  \\
& &\leq C_n \norm{f}{L^2}
 exp(C_n(\suml_k 2^{k(n-1)/2}\supl_{U\in SU(n)} \supl_x 
\| A_{k}(t,x+ Uz)\|_{L^\infty_t L^2_{z_2, \ldots, z_n}L^1_{z_1}}+ \\
& &
+ 2^{k(n+3)/2}\supl_{U\in SU(n), t\to x(t)}
\| A_{k}(t,x(t)+Uz)\|_{L^1_t L^2_{z_2, \ldots, z_n}L^1_{z_1}}))\times\\
& & \suml_{k_\al} 
\suml_{k_{\al-1}\leq k_\al} 2^{k_{\al-1}(n-1)/2}\supl_{U\in SU(n)} \supl_x 
\| \p A_{k_{\al-1}}(t,x+ 
Uz)\|_{L^\infty_t L^2_{z_2, \ldots, z_n}L^1_{z_1}}+ \\
& &
+ 2^{k_{\al-1}(n+3)/2}\supl_{U\in SU(n), t\to x(t)}
\| \p A_{k_{\al-1}}(t,x(t)+Uz)\|_{L^1_t L^2_{z_2, \ldots, z_n}L^1_{z_1}}))\times \\
& &\times (2^{k_\al(n-1)/2}\supl_{U\in SU(n)} \supl_x 
\| \p A_{k_\al}(t,x+ Uz)\|_{L^1_t L^2_{z_2, \ldots, z_n}L^1_{z_1}}+ \\
& &
+ 2^{k_\al(n-5)/2}\supl_{U\in SU(n), t\to x(t)}
\| \p A_{k_\al}(t,x(t)+Uz)\|_{L^\infty_t L^2_{z_2, \ldots, z_n}L^1_{z_1}}).
\end{eqnarray*}
We simply now replace (at the expense of a constant) 
 $\p A_{k_{\al-1}}$ by $2^{k_{\al-1}} A_{k_{\al-1}}$, 
which is in turn smaller than $2^{k_{\al}} A_{k_{\al-1}}$ 
in all the norms above involving $A_{k_{\al-1}}$.  We get an estimate 
\begin{eqnarray*}
& & \norm{\int  (\nabla\si)^2
e^{i \si(t,x,\xi)} e^{-4\pi^2 i t|\xi|^2} e^{2\pi i \dpr{\xi}{x}}
\Omega(\xi) \widehat{f}(\xi) d\xi}{L^1L^2} \leq  \\
& &\leq C_n \norm{f}{L^2} exp(C_n(\suml_k 2^{k(n-1)/2}\supl_{U, x} 
\| A_{k}(t,x+ Uz)\|_{L^\infty_t L^2_{z_2, \ldots, z_n}L^1_{z_1}}+ \\
& &
+ 2^{k(n+3)/2}\supl_{U,x(t)}
\| A_{k}(t,x(t)+Uz)\|_{L^1_t L^2_{z_2, \ldots, z_n}L^1_{z_1}}))\times \\
& & 
(\suml_{k} 2^{k(n-1)/2}\supl_{U, x} 
\| A_{k}(t,x+ 
Uz)\|_{L^\infty_t L^2_{z_2, \ldots, z_n}L^1_{z_1}}+ \\
& &\suml_k 2^{k (n+3)/2}\supl_{U,x(t)}
\| A_{k}(t,x(t)+Uz)\|_{L^1_t L^2_{z_2, \ldots, z_n}L^1_{z_1}})\times \\
& &\times (\suml_k 2^{k(n+1)/2}\supl_{U,x}  
\| \p A_{k}(t,x+ Uz)\|_{L^1_t L^2_{z_2, \ldots, z_n}L^1_{z_1}}+ \\
& &
+ 2^{k(n-3)/2}\supl_{U,x(t)}
\| \p A_{k}(t,x(t)+Uz)\|_{L^\infty_t L^2_{z_2, \ldots, z_n}L^1_{z_1}}).
\end{eqnarray*}
\noindent Next,  using the fact that 
$\norm{A}{L^2L^\infty}\leq \ve$, one sees that to control 
the term
$$
\int (i \dpr{\vec{A}}{\nabla\si})
e^{i \si(t,x,\xi)} e^{-4\pi^2 i t|\xi|^2} e^{2\pi i \dpr{\xi}{x}}
\Omega(\xi) \widehat{f}(\xi) d\xi.
$$
in $L^1L^2$, we need to control $\|\int \nabla\si(t,x,\xi) 
e^{i \si(t,x,\xi)} e^{-4\pi^2 i t|\xi|^2} e^{2\pi i \dpr{\xi}{x}}
\Omega(\xi) \widehat{f}(\xi) d\xi\|_{L^2_t L^2_x}$. 

For the proof of that, we can proceed by interpolation 
between the $L^1L^2$ estimates of Lemma \ref{le:10} and the $L^\infty L^2$ estimates
 of Lemma \ref{le:34}.  An even easier way is the following. By inspection of the 
proof of Lemma \ref{le:10}, one sees that $L^2_t L^2_x$ estimate  is 
much easier to prove (since $L^2_t$ norm 
 commutes with $L^2_x$ norm)  and one gets 
\begin{eqnarray*}
& & \|\int \nabla\si(t,x,\xi) 
e^{i \si(t,x,\xi)} e^{-4\pi^2 i t|\xi|^2} e^{2\pi i \dpr{\xi}{x}}
\Omega(\xi) \widehat{f}(\xi) d\xi\|_{L^2_t L^2_x}\leq \\
& & \leq C_n \norm{f}{L^2} exp(C_n\suml_k 2^{k(n-1)/2}\supl_{U, x} 
\| A_{k}(t,x+ Uz)\|_{L^\infty_t L^2_{z_2, \ldots, z_n}L^1_{z_1}} )\times\\
& & \times \suml_k 2^{k(n-1)/2}\supl_{U,x}
\| \p  A_{k}(t,x+ Uz)\|_{L^2_t L^2_{z_2, \ldots, z_n}L^1_{z_1}}.
\end{eqnarray*}
Note that by the convexity of the norms 
\begin{eqnarray*}
& & 2^{k(n-1)/2} \supl_{U,x}
\| \p  A_{k}(t,x+ Uz)\|_{L^2_t L^2_{z_2, \ldots, z_n}L^1_{z_1}} \leq \\
& & 
\leq (2^{k(n-1)/2} \supl_{U,x}
\|A_{k}(t,x+ Uz)\|_{L^\infty_t L^2_{z_2, \ldots, z_n}L^1_{z_1}})^{1/2}\times \\
& &\times
(2^{k(n-1)/2} \supl_{U,x}
\|\p^2 A_{k}(t,x+ Uz)\|_{L^1_t L^2_{z_2, \ldots, z_n}L^1_{z_1}})^{1/2},
\end{eqnarray*}
whence 
\begin{eqnarray*}
& & \suml_k 2^{k(n-1)/2} \supl_{U,x}
\| \p  A_{k}(t,x+ Uz)\|_{L^2_t L^2_{z_2, \ldots, z_n}L^1_{z_1}} \leq \\
& &\leq (\suml_k 2^{k(n-1)/2} \supl_{U,x}
\|A_{k}(t,x+ Uz)\|_{L^\infty_t L^2_{z_2, \ldots, z_n}L^1_{z_1}})^{1/2}\times \\
& &\times (\suml_k 2^{k(n-1)/2} \supl_{U,x}
\|\p^2 A_{k}(t,x+ Uz)\|_{L^1_t L^2_{z_2, \ldots, z_n}L^1_{z_1}})^{1/2}
\leq \norm{A}{Y_2\cap Y_3}.
\end{eqnarray*}

\noindent Thus, it remains to prove Lemma \ref{le:10}.

\subsection{Proof of Lemma \ref{le:10}.}
To outline the main ideas, 
we start the proof with the simpler case $\al=1$. 
\subsubsection{The case $\al=1$}
We follow the method of Lemma \ref{le:34}. 
According to \eqref{eq1}, \eqref{eq2}, 
$$
\Ve^{k_1}_F f(t,x)= e^{-i t \De} \suml_j P^j[
(\int F_k(t,\cdot+z\theta_j^{-k})\chi(2^{2k}z) dz)\tilde{P}^j (f)],
$$ 
where the summation in $j$ is over the family $\{\theta_j^{-k}\}$. 
Taking $L^2_x$ norms and taking into account the (almost) 
orthogonality of the $j$ sum, we have 
\begin{eqnarray*}
& &\norm{\Ve^{k_1}_F f(t,\cdot)}{L^2}\leq C_n  \left(\suml_j 
\norm{(\int F_k(t,\cdot+z\theta_j^{-k})
\chi(2^{2k}z) dz)\tilde{P}^j (f)}{L^2}^2\right)^{1/2}\leq \\
& &\leq C_n \left(\suml_j \supl_x|h_j(t,x)|^2 
\norm{\tilde{P}^j f}{L^2}^2\right)^{1/2},
\end{eqnarray*}
where $h_j(t,x)= \int F_k(t,x+z\theta_j^{-k})
\chi(2^{2k}z) dz$. 
Take  $L^1_t$ norm. We will show  
\begin{equation}
\label{eqn20}
\begin{array}{l}
\norm{\left(\suml_j \supl_x|h_j(t,x)|^2 
\norm{\tilde{P}^j f}{L^2}^2\right)^{1/2}}{L^1_t}\leq \\
\leq C_n
2^{k(n-1)/2} \supl_{U\in SU(n)} 
\|\supl_x\|F_k(t, x+Uz)\|_{L^{2}_{z_2,\ldots, z_n} 
L^1_{z_1}}\|_{L^1_t}\norm{f}{L^2}.
\end{array}
\end{equation}
This follows by complex interpolation between 
\begin{equation}
\label{eqn10}
\begin{array}{l}
\norm{\suml_j \supl_x|h_j(t,x)|
\norm{\tilde{P}^j f}{L^2}}{L^1_t}\leq \\
\leq 
\supl_{U\in SU(n)} 
\norm{\supl_x
\norm{F_k(t,x+Uz)}{L^{\infty}_{z_2,\ldots, z_n}L^1_{z_1}}}{L^1_t} 
\suml_j \norm{\tilde{P}^j f}{L^2}
\end{array}
\end{equation}
and
\begin{equation}
\label{eqn11}
\norm{\supl_j[ \supl_x|h_j(t,x)| 
\|\tilde{P}^j f\|_{L^2}]}{L^1_t}\leq  C_n 
2^{k(n-1)} \norm{F_k}{L^1_t L^1_z}\sup_j \norm{\tilde{P}^j f}{L^2},
\end{equation}
because $\norm{F_k}{L^1_t L^1_x}=\sup_{U\in SU(n)} 
\norm{\sup_x
\norm{F_k(t,x+Uz)|}{L^{1}_z}}{L^1_t}$. \\
The $l^1_j$ estimate,  \eqref{eqn10} is straightforward, since 
$\tilde{P}^j f$ is indendent of $t$. We have 
\begin{eqnarray*}
& & \norm{\suml_j \supl_x|h_j(t,x)|
\norm{\tilde{P}^j f}{L^2}}{L^1_t}=  \suml_j
\norm{\supl_x |h_j(t,x)|}{L^1_t} 
 \norm{\tilde{P}^j f}{L^2}\leq \\
& &\leq \supl_j \norm{\supl_x |h_j(t,x)|}{L^1_t} 
\suml_j \norm{\tilde{P}^j f}{L^2}
\end{eqnarray*}
It remains to observe that 
\begin{eqnarray*}
& & \supl_j \norm{\supl_x |h_j(t,x)|}{L^1_t}\leq 
\supl_\theta \norm{\supl_x \int |F_k(t,x+z\theta)|dz}{L^1_t}= \\
& &=
\supl_{U\in SU(n)} 
\norm{\supl_x\norm{F_k(t,x+Uz)|}{L^\infty_{z_2, \ldots, z_n}L^1_{z_1}}}{L^1_t}.
\end{eqnarray*}
For the $l^\infty_j$ estimate,  \eqref{eqn11},  write 
\begin{eqnarray*}
& & \norm{\supl_j [\supl_x|h_j(t,x)| 
\norm{\tilde{P}^j f}{L^2}]}{L^1_t}=  \norm{\supl_x\supl_j |h_j(t,x)| }{L^1_t} 
\supl_j \norm{\tilde{P}^j f}{L^2}.
\end{eqnarray*}
Note that 
\begin{eqnarray*}
& & 
\supl_x \supl_j |h_j(t,x)|\leq \supl_\theta \supl_x 
\suml_{l\leq -2k} \int |F_k(t,x+z\theta)|\vp(2^{-l}z) dz\leq \\
& &\leq 
\suml_{l\leq -2k} \supl_\theta \supl_x 
 \int |F_k(t,x+z\theta)|\vp(2^{-l}z) dz.
\end{eqnarray*}
Since for every fixed $\theta$, 
$2^{(l+k)(n-1)}\int_{\theta_0\in\sn: |\theta_0-\theta|<2^{-l-k}}
d\theta_0 \gtrsim  1$, write 
\begin{eqnarray*}
& &  \int |F_k(t,x+z\theta)|\vp(2^{-l}z) dz\leq \\
& &\leq C_n  
2^{(l+k)(n-1)} \int_{\theta_0\in\sn: |\theta_0-\theta|<2^{-l-k}} 
\int |F_k(t,x+z\theta)|\vp(2^{-l}z) dz  
 d\theta_0
\end{eqnarray*}
Expand out as before $F_k(t,x+z\theta)$ around $x+z\theta_0$. We get 
\begin{eqnarray*}
& &  \supl_\theta\supl_x  \int |F_k(t,x+z\theta)|\vp(2^{-l}z) dz\leq \\
& &\leq C_n 
2^{(l+k)(n-1)} \suml_{\ga\geq 0}(\ga!)^{-1} 
\int_{\theta_0\in\sn: |\theta_0-\theta|<2^{-l-k}} 
\int |F_k^\ga(t,x+z\theta_0)|\vp^\ga(2^{-l}z) dz d\theta_0\leq \\
& & \leq C_n 
2^{(l+k)(n-1)} \suml_{\ga\geq 0}(\ga!)^{-1} 
\int_{\theta_0\in\sn} 
\int |F_k^\ga(t,x+z\theta_0)|\vp^\ga(2^{-l}z) dz d\theta_0\leq \\
& &\leq 
C_n 
2^{k(n-1)} \suml_{\ga\geq 0}(\ga!)^{-1} 
\int_{\theta_0\in\sn} 
\int |F_k^\ga(t,x+z\theta_0)||z|^{n-1}\vp^\ga(2^{-l}z) dz d\theta_0.
\end{eqnarray*}
Summing the last inequalities in $l$ implies 
\begin{equation}
\label{eqn18}
\begin{array}{l}
|\supl_x \supl_j h_j(t,x)|\leq C_n 
2^{k(n-1)} \suml_{\ga\geq 0}(\ga!)^{-1}  \\\suml_l  \int_{\theta_0\in\sn} 
\int |F_k^\ga(t,x+z\theta_0)||z|^{n-1}\vp^\ga(2^{-l}z) dz d\theta_0 
\leq 
C_n 
2^{k(n-1)} \norm{F_k(t, \cdot)}{L^1_x}
\end{array}
\end{equation}
Taking $L^1_t$ norm  yields 
$$
 \norm{\supl_x \supl_j |h_j(t,x)|}{L^1_t}\leq C_n 
2^{k(n-1)} \norm{F_k}{L^1_{t x}},
$$
as required. 

\subsubsection{The case $\al>1$}
The strategy here is to start ``peeling off'' the functions 
$\int F^\mu_{k_\mu}(t,x+z\xi/|\xi|) \chi(2^{2k_\mu}z) dz$ in a way 
 similar to Lemma \ref{le:34}. Recall that the method presented 
in Lemma \ref{le:34} starts the ``peeling'' argument with the terms 
with the lowest frequency \\ 
$k_{\min}=\min(k_1, \ldots, k_\al)$. 

We have here 
the extra complication of having to take $L^1_t$. Moreover, we must  
measure  the  
last  term $\int F^\al_{k_\al}(t,x+z\xi/|\xi|) \chi(2^{2k_\al}z) dz$ 
in the $L^1_t$ norm, while all the other terms should be 
 measured in $L^\infty_t$. 

As we have alluded to above, the order of the frequencies is not 
insignificant.  
The case $\al=1$,  considered in the 
previous section roughly corresponds to the case when the last  frequency 
$k_\al$ is maximal, i.e. $k_\al=\max(k_1, \ldots, k_\al)$. 
We consider this case, and then 
we indicate the necessary changes 
when $k_\al<\max(k_1, \ldots, k_\al)$. \\
{\bf SubCase 1: $k_\al=\max(k_1, \ldots, k_\al)$}.\\
By the symmetry of  the terms $1, \ldots, (\al-1)$, 
assume without loss of generality $k_1\leq \ldots\leq  k_{\al-1}$. 
In fact, following Lemma \ref{le:34}, let 
$k_1=\ldots=k_{s_1}<k_{s_1+1}=\ldots=k_{s_2}<\ldots 
k_{s_{m-1}+1}=\ldots=k_{s_m}=k_\al$. 

According to \eqref{eqn12}, we can write 
$\Ve^{k_1, \ldots, k_\al}_{F_1, \ldots, F_\al}$ 
as a free solution (with time  depending data). Recall also 
\eqref{eqn:7} and \eqref{eqn:1}, which give an 
estimate for its $L^2$ norm (for a fixed time $t$). 
Take into account \eqref{eqn14}, to get 
 \begin{eqnarray*}
& &
\norm{\Ve^{k_1, \ldots, k_\al}_{F_1, \ldots, F_\al} f}{L^2_x} 
\leq  \\
& &\leq C_n  (\suml_{j_1} 
\suml_{j_2,  \ldots, j_m:|\theta_{j_r}^{-k_{s_r}}-\theta_{j_1}^{-k_{s_1}}|\leq 
2^{k_{s_1}}+2^{k_{s_r}}}
 \norm{g_1^{j_1} \ldots g_{s_1}^{j_1} Q^{j_1}_{2^{k_{s_1}}+2s_1 2^{k_{s_1}}}
(f(\cdot)\prod_{\mu=s_1+1}^\al g_\mu^{j_{r(\mu)}}(t,\cdot))}{L^2_x}^2)^{1/2},
\end{eqnarray*}
where for $\mu\in [s_{r-1}+1, s_r]$, we have introduced 
$$
g_\mu^{j_r}(t,x)= \int F_{k_\mu}^\mu(t,x+z
\theta_{j_r}^{-k_\mu})\chi(2^{2k_\mu}z) dz.
$$
Furthermore, by support considerations 
(as discussed in the proof of Lemma \ref{le:34}),  
$$
Q^{j_1}_{2^{k_{s_1}}+2s_1 2^{k_{s_1}}}=Q^{j_1}_{2^{k_{s_1}}+2s_1 2^{k_{s_1}}} 
Q^{j_2}_{2^{k_{s_2}}+(2s_1+2)2^{k_{s_1}}}\ldots 
Q^{j_r}_{2^{k_{s_r}}+(2s_1+2)2^{k_{s_1}}}.
$$
We conclude 
\begin{equation}
\label{eqn16}
\begin{array}{l}
\norm{\Ve^{k_1, \ldots, k_\al}_{F_1, \ldots, F_\al} f}{L^2_x} 
\leq C_n  (\suml_{j_1, j_2, \ldots, j_m} 
 \|g_1^{j_1}(t,x) \ldots g_{s_1}^{j_1}(t,x) \times \\
\times   Q^{j_1}_{2^{k_{s_1}}+2s_1 2^{k_{s_1}}}\ldots 
Q^{j_r}_{2^{k_{s_r}}+(2s_1+2)2^{k_{s_1}}}
(f(\cdot)\prod_{\mu=s_1+1}^\al g_\mu^{j_r(\mu)}(t,\cdot))\|_{L^2_x}^2)^{1/2}.
\end{array}
\end{equation}
Take $L^1_t$ norm. Our  claim is that 
$$
\begin{array}{l}
\norm{\Ve^{k_1, \ldots, k_\al}_{F_1, \ldots, F_\al} f}{L^1_t L^2_x}
\leq C_n (1+2s_1)^n 
\prod_{\mu=1}^{s_1} 2^{k_\mu(n-1)/2}\supl_t \supl_{U\in SU(n)} 
\supl_x \norm{F_{k_\mu}^\mu(t, x+Uz)}{L^2_{z_2, \ldots, z_n}L^1_{z_1}} \\
\times 
 \|(\suml_{j_2,  \ldots, j_m} 
\norm{Q^{j_2}_{2^{k_{s_2}}+(2s_1+2)2^{k_{s_1}}}\ldots 
Q^{j_m}_{2^{k_{s_m}}+(2s_1+2)2^{k_{s_1}}} (f(\cdot)
\prod_{\mu=s_1+1}^\al g_\mu^{j_{r(\mu)}}(t,\cdot))}{L^2}^2)^{1/2}\|_{L^1_t}.
\end{array}
$$
On account of \eqref{eqn16}, this follows  
by a complex multilinear  interpolation between the $l^1_j$ estimate
\begin{eqnarray*}
& & \|\suml_{j_1, j_2,  \ldots, j_m} 
 \|g_1^{j_1} \ldots g_{s_1}^{j_1} Q^{j_1}_{2^{k_{s_1}}+2s_1 2^{k_{s_1}}}\ldots 
Q^{j_m}_{2^{k_{s_m}}+(2s_1+2)2^{k_{s_1}}}
(f(\cdot)\prod_{\mu=s_1+1}^\al g_\mu^{j_{r(\mu)}}(t,\cdot))\|_{L^2_x}\|_{L^1_t} \leq \\
& &\leq \prod_{\mu=1}^{s_1} (\supl_t \supl_{\theta, x} 
\int |F_{k_\mu}^\mu(t,x+z\theta)|dz) \times \\
& &\times \|\suml_{ 
j_1, j_2,  \ldots, j_m} \| Q^{j_1}_{2^{k_{s_1}}+2s_1 2^{k_{s_1}}}
Q^{j_2}_{2^{k_{s_2}}+(2s_1+2)2^{k_{s_1}}}\ldots 
Q^{j_m}_{2^{k_{s_m}}+(2s_1+2)2^{k_{s_1}}} (f(\cdot)
\prod_{\mu=s_1+1}^\al g_\mu^{j_{r(\mu)}}(t,\cdot))\|_{L^2}\|_{L^1_t}.
\end{eqnarray*}
and the $l^\infty_j$ estimate
\begin{eqnarray*}
& &\| \supl_{j_1,j_2, \ldots, j_m} 
 \|g_1^{j_1} \ldots g_{s_1}^{j_1} Q^{j_1}_{2^{k_{s_1}}+2s_1 2^{k_{s_1}}}
\ldots 
Q^{j_m}_{2^{k_{s_m}}+(2s_1+2)2^{k_{s_1}}}
(f(\cdot)\prod_{\mu=s_1+1}^\al g_\mu^{j_{r(\mu)}}(t,\cdot))\|_{L^2_x}\|_{L^1_t}
 \leq \\
& &\leq C_n \prod_{\mu=1}^{s_1} 2^{k_\mu(n-1)}
\norm{F_{k_\mu}^\mu}{L^\infty_t L^1_x} \times \\
& &\times 
\|\supl_{j_1, j_2, \ldots, j_m}\|Q^{j_1}_{2^{k_{s_1}}+2s_1 2^{k_{s_1}}} Q^{j_2}_{2^{k_{s_2}}+
(2s_1+2)2^{k_{s_1}}}\ldots 
Q^{j_m}_{2^{k_{s_m}}+(2s_1+2)2^{k_{s_1}}} (f(\cdot)
\prod_{\mu=s_1+1}^\al g_\mu^{j_{r(\mu)}}(t,\cdot))\|_{L^2}\|_{L^1_t}.
\end{eqnarray*}
Note that after the complex interpolation, the sum in $j_1$ disappears (at the expense of a constant $C_n(1+2s_1)^n$), 
since we have estimated by \eqref{eqn17}. 

The proof of the $l^1_j$ estimate is immediate, by just pulling  out \\
$$
\prod_{\mu=1}^{s_1}
\sup_{j_1,t,x} |g_1^{j_1}(t,x)\ldots g_{s_1}^{j_1}(t,x)|\leq 
\prod_{\mu=1}^{s_1} \supl_t \supl_{\theta, x} 
\int |F_{k_\mu}^\mu(t,x+z\theta)|dz.
$$
For the $l^\infty_j$ estimate, estimate by 
\begin{eqnarray*}
& &\| \supl_{j_1,j_2, \ldots, j_m} 
 \|g_1^{j_1} \ldots g_{s_1}^{j_1} Q^{j_1}_{2^{k_{s_1}}+2s_1 2^{k_{s_1}}}
\ldots 
Q^{j_m}_{2^{k_{s_m}}+(2s_1+2)2^{k_{s_1}}}
(f(\cdot)\prod_{\mu=s_1+1}^\al g_\mu^{j_{r(\mu)}}(t,\cdot))\|_{L^2_x}\|_{L^1_t}
 \leq \\
& &\leq \prod_{\mu=1}^{s_1} \supl_t \supl_{j_1,x} |g_\mu^{j_1}(t,x)| \times \\
& &\times \|\supl_{j_1,j_2, \ldots, j_m} 
\|Q^{j_1}_{2^{k_{s_1}}+2s_1 2^{k_{s_1}}}
\ldots 
Q^{j_m}_{2^{k_{s_m}}+(2s_1+2)2^{k_{s_1}}}
(f(\cdot)\prod_{\mu=s_1+1}^\al 
g_\mu^{j_{r(\mu)}}(t,\cdot))\|_{L^2_x}\|_{L^1_t}.
\end{eqnarray*}
By \eqref{eqn18} however, 
$$
\supl_t \supl_{j_1,x} |g^{j_1}_\mu(t,x)|\leq C_n 
2^{k_\mu(n-1)}\supl_t \norm{F^\mu_{k_\mu}(t,\cdot)}{L^1_x},
$$
as required. 

This shows the main step in the argument. We continue in this 
fashion (and as in Lemma \ref{le:34}, we keep incurring 
constants comming from the increased 
overlap of the supports of  the operators $Q^j$) 
until we reach the maximal frequency $k_\al$. In this final step, 
we finally have to  take the $L^1_t$ norm on 
$\int F_{k_\al}(t, x+z\theta)\chi(2^{2k_\al}z) dz$. We have 
\begin{eqnarray*}
& & \norm{\Ve^{k_1, \ldots, k_\al}_{F^1, \ldots, F^\al}}{L^1_t L^2_x}\leq \\
& & \leq C_n^\al \prod_{\mu=1}^{k_{s_{m-1}}} 
2^{k_\mu(n-1)/2}\supl_t \supl_{U\in SU(n)} 
\supl_x \norm{F_{k_\mu}^\mu(t, x+Uz)}{L^2_{z_2, \ldots, z_n}L^1_{z_1}}\times \\
& & \times  \|(\suml_{j_m} 
\|\prod_{\mu=s_{m-1}+1}^\al g_\mu^{j_m}(t,\cdot) 
Q^{j_m}_{2^{k_{s_m}}+(2s_1+2)2^{k_{s_1}}+
\ldots (2(s_{m}-s_{m-1})+2)2^{k_{s_m}}} [f(\cdot)]\|_{L^2}^2)^{1/2}\|_{L^1_t}.
\end{eqnarray*}
Observe that since $p(\al,m)=2^{k_{s_m}}+(2s_1+2)2^{k_{s_1}}+
\ldots (2(s_{m}-s_{m-1})+2)2^{k_{s_m}}\leq 4\al 2^{k_{s_m}}$, we have
$$
(\suml_{j_m} \norm{Q^{j_m}_{p(\al,m)}g}{L^2}^2)^{1/2}\leq 
C_n (p(\al,m)/2^{k_{s_m}})^n \norm{g}{L^2}\leq C_n \al^{n}\norm{g}{L^2} . 
$$  
It remains  to show 
\begin{eqnarray*}
& & \|(\suml_{j_m} 
\|\prod_{\mu=s_{m-1}+1}^\al g_\mu^{j_m}(t,\cdot) 
Q^{j_m}_{p(\al, m)} [f(\cdot)]\|_{L^2}^2)^{1/2}\|_{L^1_t}
\leq C_n^\al (\prod_{\mu=s_{m-1}+1}^{\al-1} \supl_{t, \theta, x}\int 
|F_{k_\mu}^\mu(t,x+z \theta)|dz)\times \\
& &\times   2^{k_\al(n-1)/2}
\supl_{U\in SU(n)} 
\|\supl_x \|F_{k_\al}(t, x+Uz)
\|_{L^2_{z_2, \ldots, z_n}L^1_{z_1}}\|_{L^1_t}\norm{f}{L^2}. 
\end{eqnarray*}
This follows as in \eqref{eqn20}, by interpolating between 
\begin{eqnarray*}
& & \|\suml_{j_m} 
\|\prod_{\mu=s_{m-1}+1}^\al g_\mu^{j_m}(t,\cdot) 
Q^{j_m}_{p(\al,m)} [f(\cdot)]\|_{L^2}\|_{L^1_t}
\leq \\
& &\leq C_n (\prod_{\mu=s_{m-1}+1}^{\al-1} \supl_{t, \theta, x}\int 
|F_{k_\mu}^\mu(t,x+z \theta)|dz)
 \suml_{j_m} 
\|Q^{j_m}_{p(\al,m)} [f(\cdot)]\|_{L^2}\times\\
& &\times \supl_{\theta}  \|\supl_x\|\int 
|F_{k_\al}^\al(t,x+z \theta)|dz\|_{L^1_t} 
\end{eqnarray*}
and 
\begin{eqnarray*}
& & \|\supl_{j_m} 
\|\prod_{\mu=s_{m-1}+1}^\al g_\mu^{j_m}(t,\cdot) 
Q^{j_m}_{p(\al,m)} [f(\cdot)]\|_{L^2}\|_{L^1_t}
\leq \\
& &\leq C_n (\prod_{\mu=s_{m-1}+1}^{\al-1} 2^{k_\mu(n-1)} 
\norm{F_{k_\mu}^\mu}{L^\infty_t L^1_{x}}) 
 \supl_{j_m} 
\|Q^{j_m}_{p(\al,m)} [f(\cdot)]\|_{L^2}\times\\
& &\times 2^{k_\al(n-1)} 
\norm{F_{k_\al}^\al}{L^1_t L^1_{x}}) .
\end{eqnarray*}
We omit the details, as they are exactly as in the proof of \eqref{eqn20}.
This completes the proof of Subcase I.  \\ 
{\bf Subcase II: $k_\al<\max(k_1, \ldots, k_\al)$}\\
Let for some $1\leq r_0 <\al$, we have $k_{r_0}= 
\max(k_1, \ldots, k_\al)$. 
The estimates proved in Subcase I yield  
\begin{eqnarray*}
& & \norm{\Ve^{k_1, \ldots, k_\al}_{F^1, \ldots, F^\al}}{L^1_t L^2_x}\leq 
 C_n^\al 
(\prod_{\mu=1,\mu\neq r_0}^{\al-1} 
2^{k_\mu(n-1)/2}\supl_t \supl_{U\in SU(n)} 
\supl_x \norm{F_{k_\mu}^\mu(t, x+Uz)}{L^2_{z_2, \ldots, z_n}L^1_{z_1}})
 \times \\
& & \times 2^{k_\al(n-1)/2}\supl_t \supl_{U\in SU(n)} 
\supl_x \norm{F_{k_\al}^\al(t, x+Uz)}{L^2_{z_2, \ldots, z_n}L^1_{z_1}}\times\\
& &\times 2^{k_{r_0}(n-1)/2}
\supl_{U\in SU(n)} 
\|\supl_x \|F_{k_{r_0}}^{r_0}(t, x+Uz)\|_{L^2_{z_2, \ldots, 
z_n}L^1_{z_1}}\|_{L^1_t} 
\norm{f}{L^2}
\end{eqnarray*}
Since $k_\al<\max(k_1, \ldots, k_\al)=k_{r_0}$, we have $2^{k_\al(n-1)/2}2^{k_{r_0}(n-1)/2}\leq 2^{k_\al(n-5)/2}2^{k_{r_0}(n+3)/2}$. We get an estimate 
\begin{eqnarray*}
& & \norm{\Ve^{k_1, \ldots, k_\al}_{F^1, \ldots, F^\al}}{L^1_t L^2_x}\leq  
 C_n^\al 
(\prod_{\mu=1,\mu\neq r_0}^{\al-1} 
2^{k_\mu(n-1)/2}\supl_t \supl_{U\in SU(n)} 
\supl_x \norm{F_{k_\mu}^\mu(t, x+Uz)}{L^2_{z_2, \ldots, z_n}L^1_{z_1}})
 \times \\
& & \times 2^{k_\al(n-5)/2}\supl_t \supl_{U\in SU(n)} 
\supl_x \norm{F_{k_\al}^\al(t, x+Uz)}{L^2_{z_2, \ldots, z_n}L^1_{z_1}}\times\\
& &\times 2^{k_{r_0}(n+3)/2}
\supl_{U\in SU(n)} 
\|\supl_x \|F_{k_{r_0}}^{r_0}(t, x+Uz)\|_{L^2_{z_2, \ldots, 
z_n}L^1_{z_1}}\|_{L^1_t} 
\norm{f}{L^2}
\end{eqnarray*}
This finishes the Proof of Lemma \ref{le:10}.

\section{Strichartz estimates in  
$L^2_t L^{\f{2(n-1)}{(n-3)}}_{x_2, \ldots, x_n}L^2_{x_1}$} 
\label{sec:8}
In this section, we establish Proposition \ref{prop:3}. We concentrate 
on the case $n\geq 4$, since this is what we use anyway. 
On the  other hand, 
minor changes are needed in the proofs  for 
$n=2$ and $n=3$. 

Let us show first 
Proposition \ref{prop:3} for the case $A=0$, that is for free solutions. 
\begin{lemma}
\label{le:890}
Let $n\geq 4$. Then 
\begin{eqnarray}
\label{eq:600}
& & \supl_{U\in SU(\rn), x(t)} \norm{e^{i t \De} f(x(t)+U z)}{L^2_t L^{2(n-1)/(n-3)}_{z_2, \ldots, z_n}L^2_{z_1}}\leq C_n \norm{f}{L^2} \\
\label{eq:601}
& & 
\supl_{U\in SU(\rn), x(t)} \norm{ \intl_0^t e^{i (t-s) \De} 
F(s,(x(s)+U z)ds)}{L^2_t L^{2(n-1)/(n-3)}_{z_2, \ldots, z_n}L^2_{z_1}}\leq C_n
\norm{F}{L^1L^2}.
\end{eqnarray}
\end{lemma}
\begin{proof}
For a fixed $U\in SU(\rn)$ and fixed measurable function $x(t)$, 
denote \\ $W(t) f(y):=(e^{i t \De} f)(x(t)+Uy)$. According 
to Proposition \ref{prop:14}, we need to verify that $W(t):L^2\to L^2$ and 
$W(t)W(s)^*:L^1_{y_2, \ldots, y_n} L^2_{y_1}\to 
L^\infty_{y_2, \ldots, y_n} L^2_{y_1}$ with norm no larger than \\
$C_n |t-s|^{(n-1)/2}$. \\
The $L^2$ boundedness is obvious since, 
$$
\norm{W(t)f}{L^2}=\norm{(e^{it \De}f)(x(t)+U y)}{L^2_y}=
\norm{(e^{it \De}f)(x(t)+y)}{L^2_y}= 
\norm{e^{it \De}f}{L^2_y}=\norm{f}{L^2}.
$$ 

For the dispersive estimates, note that 
$$
W(t)W(s)^* f(y)= (e^{i (t-s)\De} f)(y+U^*x(s)-U^* x(t)).
$$
It suffices to verify 
the decay estimate for a familly of extreme points, whose convex  
span is $w^*$ dense in the unit ball of 
$L^1_{y_2, \ldots, y_n} L^2_{y_1}$. Since the $\de$ functions 
provide such a set in $L^1_X$ for any measure space $X$, it will 
suffice to take $f(y)=\de(\bar{y}-b) g(y_1)$, where 
$\bar{y}=(y_2, \ldots, y_n)$, $b\in {\mathbf R}^{n-1}$ and $g\in L^2(\rone)$. 
Fix $s$ and $t$ and denote $z=U^*x(s)-U^* x(t)$, which is a fixed vector in $\rn$. 
Compute $W(t)W(s)^* f$
$$
W(t)W(s)^* f(y)=(e^{i (t-s) \p_1^2} g)(y_1+z_1) 
[e^{i(t-s)\De_{n-1}}\de(\cdot-b)](\bar{y}+\bar{z}).
$$
Immediately, 
$\norm{W(t)W(s)^* f(\bar{y})}{L^2_{y_1}}=
|e^{i(t-s)\De_{n-1}}\de(\cdot-b)(\bar{y}+\bar{z})|\|g\|_{L^2_{y_1}}$ and thus  
$$
\norm{W(t)W(s)^* f}{L^\infty_{y_2, \ldots, y_n} 
L^2_{y_1}}\leq C_n|t-s|^{-(n-1)/2}\|g\|_{L^2_{y_1}}= 
C_n|t-s|^{-(n-1)/2}\|f\|_{L^1_{z_2, \ldots, z_n}L^2_{y_1}}.
$$
Note that the constant $C_n$ is independent of $g$, and $z\in \rn$ 
and depends only on the dimension $n$. 
\end{proof}
 
For the case $A\neq 0$, we proceed as in the proof of Theorem 
\ref{theo:pisa}. We first establish a ``naive'' Strichartz estimate similar to Proposition \ref{prop:7}. Then, we show that to extend the 
naive Strichartz estimate indefinitely in time (under the assumption 
that $A$ is small), we need an estimate of the parametrix, similar to Lemma \ref{le:parametrix} in the appropriate norms. 

We start with the 
naive Strichartz estimate. This is done 
exactly the same way as Proposition \ref{prop:7}, 
given that we already have Lemma \ref{le:890}. 
\begin{proposition}
\label{prop:8} For a fixed integer  $k_0$, there exists a time 
$T_0=T_0(k_0)\leq \infty$, so that 
whenever $0<T< T_0$, $\psi\in \cs$  
\begin{equation}
\label{eq:1820}
\supl_{U\in SU(\rn), x(t)} \norm{(P_{<k_0}\psi)(t, x(t)+U z)}{L^2(0,T)
L^{2(n-1)/(n-3)}_{z_2, \ldots, z_n}L^2_{z_1}}\leq 
C(T,k_0)(\norm{\psi(0,\cdot)}{L^2}+
\norm{\cl\psi}{L^1L^2}).
\end{equation}
Moreover, $C(T,k_0)$ depends on  $T$ in a continuous way.
\end{proposition}
Fix $k_0$ and a small $\ve$. 
Set as before $0<T^*\leq \infty$ to be the maximum time, 
so that for every $0<T<T^*$, one has 
$$
\supl_{U\in SU(\rn), x(t)} \norm{(P_{<k_0}\psi)(t, x(t)+U z)}{L^2(0,T)
L^{2(n-1)/(n-3)}_{z_2, \ldots, z_n}L^2_{z_1}}\leq 
\ve^{-1}(\norm{\psi(0,\cdot)}{L^2}+
\norm{\cl\psi}{L^1L^2}),
$$
for all Schwartz functions $\psi$. 

The goal would be again to show that for small 
enough $\ve$, we have that $T^*=\infty$. This is reduced in a standard way (recall that $P_{<k_0}$ has an integrable kernel and is therefore bounded on 
$L^2(0,T)
L^{2(n-1)/(n-3)}_{z_2, \ldots, z_n}L^2_{z_1}$) to the following 
estimate for the parametrix, constructed in Section \ref{sec:par}. This 
 needs to be compared to Lemma \ref{le:parametrix}. 
\begin{lemma}
\label{le:9}
Given 
$\norm{\vec{A}}{\tilde{Y}_1\cap Y_2\cap Y_3}\leq \ve$, integer $k$ and $T>0$, 
and for every function $f_k\in L^2(\rn)$ with $supp\ \widehat{f_k}
\subset \{|\xi|\sim 2^k\}$, 
one can find a 
function 
$v_k:[0,T)\times \rn\to \cc$, so that $supp\ 
\widehat{v_k}\subset \{|\xi|\sim 2^k\}$ and 
\begin{eqnarray*}
& &\norm{v_k(0,x)-f_k}{L^2}\leq C\ve \norm{f_k}{L^2} \\
& & \supl_{U\in SU(n), x(t)}  
\norm{v_k(t, x(t)+Uz)}{L^2_tL^{2(n-1)/(n-3)}_{z_2, \ldots, z_n}L^2_{z_1}}
\leq C\norm{f_k}{L^2}\\
& &\norm{\cl v_k}{L^1L^2}\leq C\ve \norm{f_k}{L^2},
\end{eqnarray*}
for some $C$ independent of $f, k, T$.
\end{lemma}
\noindent It remains to prove Lemma \ref{le:9}. 
\begin{proof}(\ref{le:9})
We reduce as in Section \ref{sec:par} to the case, when $k=0$ and 
$supp \hat{A}\subset 
\{\xi: |\xi|<<1\}$ and {\it without the condition}
 $supp\ \hat{v}\subset \{\xi: |\xi|\sim 1\}$.
We again use the function $v$ constructed in Section \ref{sec:par}. \\
Since we have already verified 
$\norm{v(0,x)-f}{L^2}\leq C\ve \norm{f}{L^2}$ and 
$\norm{\cl v}{L^1L^2}\leq C\ve \norm{f}{L^2}$ under the consitions imposed 
in Lemma \ref{le:9}, it remains to check 
$$
\supl_{U\in SU(n), x(t)}  
\norm{v(t, x(t)+Uz)}{L^2(0,T)
L^{2(n-1)/(n-3)}_{z_2, \ldots, z_n}L^2_{z_1}}\leq C\norm{f}{L^2}.
$$
for 
$$
v(t,x)=\La f(t,x) = \int e^{i \si(t,x,\xi)} e^{-4\pi^2 i t|\xi|^2} e^{2\pi i \dpr{\xi}{x}}\Om(\xi) \hat{f}(\xi) d\xi.
$$
Expanding $e^{i \si}=\suml_\al (i^\al \si^\al)/(\al!)$ reduces 
matters to showing 
\begin{equation}
\label{eq:op1}
\supl_{U\in SU(n), x(t)}   \norm{[\La^\al f](t,x(t)+Ux)}{L^2(0,T)
L^{2(n-1)/(n-3)}_{x_2, \ldots, x_n}L^2_{x_1}}\leq C_n^\al \ve^\al 
\norm{f}{L^2}.
\end{equation}
Fix $U\in SU(n)$ and a measurable function $x(t)$. Set 
$$
W^\al(t) f(x)= [\La^\al f](t,x(t)+Ux)=
\int \si^\al(t,x(t)+Ux,\xi) e^{-4\pi^2 i t|\xi|^2} 
e^{2\pi i \dpr{\xi}{x(t)+Ux}}\Om(\xi) \hat{f}(\xi) d\xi.
$$
By Proposition \ref{prop:14}, \eqref{eq:op1} follows from  
 the energy estimate 
$$
\norm{W^\al(t) f}{L^2}\leq C_n^\al(\suml_k 
\supl_{U,x} \|A_k(t, x+Uz)\|_{L^\infty_t 
L^\infty_{z_2, \ldots, z_n}L^1_{z_1}})\al 
\norm{f}{L^2}
$$ 
and 
the decay estimate 
$$
\norm{W^\al(t)W^\al(s)^* f}{L^{\infty}_{x_2, 
\ldots, x_n} L^2_{x_1}}\leq C_n^{\al} (\suml_k 2^{k(n-1)}
\|A_k\|_{L^\infty_t L^1_{x}}) ^{2\al} |t-s|^{-(n-1)/2} 
\norm{f}{L^1_{x_2, 
\ldots, x_n} L^2_{x_1}}.
$$
By interpolation between the last two estimates, we obtain the ``modified decay estimate''
\begin{eqnarray*}
& &
\norm{W^\al(t)W^\al(s)^* f}{L^{p'}_{x_2, 
\ldots, x_n} L^2_{x_1}}\leq \\
& &\leq C_n^{\al} (\suml_k 2^{k(n-1)/p_0}
\supl_{U,x} \|A_k(t, x+Uz)\|_{L^\infty_t 
L^{p_0}_{z_2, \ldots, z_n}L^1_{z_1}}) ^{2\al} |t-s|^{- (n-1)/(2p_0)} 
\norm{f}{L^p_{x_2, 
\ldots, x_n} L^2_{x_1}}= \\
& &= C_n^\al |t-s|^{-(n-1)/(2p_0)} \norm{A}{\tilde{Y}_1}^{2\al}
\norm{f}{L^p_{x_2, 
\ldots, x_n} L^2_{x_1}},
\end{eqnarray*}
which implies \eqref{eq:op1}, as long as $p_0<(n-1)/2$, 
according to the Remark before  Proposition \ref{prop:14}. 
\subsection{Energy estimates for $W^\al$}
Observe that the case $\al=0$ is 
simply the energy estimate in Lemma \ref{le:890}, while for 
the general case  
observe\footnote{Recall  that $\Om$ is a radial  function.}
\begin{eqnarray*}
& &
W^\al(t) f(x)= [\La^\al f](t,x(t)+Ux)= \\
& &= \int \si^\al(t,x(t)+Ux,U^* \xi) e^{-4\pi^2 i t|\xi|^2} 
e^{2\pi i \dpr{\xi}{U^* x(t)+x-y}}\Om(\xi) d\xi f(Uy) dy.
\end{eqnarray*}
By a simple translational invariance, 
$\norm{f(U\cdot)}{L^2}=\norm{f}{L^2}$ and 
$\norm{A(t, x(t)+ U\cdot)}{Y_j}=\norm{A}{Y_j}$, 
we have that by 
\eqref{eq:906}, the energy estimate $\norm{W^\al(t) f}{L^2}\leq C_n^\al\ve^\al \norm{f}{L^2}$  is satisfied for every $\al>0$. 
\subsection{Decay estimates for $W^\al$}
The case $\al=0$ is  the decay estimate in Lemma \ref{le:890}. 

For $\al\geq 1$, we are following the approach of the dispersive estimates in 
Section \ref{sec:6}. Write 
\begin{eqnarray*}
& &W^{\al}(t)W^{\al}(s)^*f(x)= \\
& &=
\prod_{\mu=1}^\al (\suml_{k_\mu\leq -2}
\suml_{l_\mu\leq -2k_\mu} 
\suml_{\ga_\mu} (\ga_\mu!)^{-1} 
\suml_{j_\mu}  
\int A^{\ga_\mu}_{k_\mu}(t, x(t)+Ux +z \theta^{l_\mu+k_\mu}_{j_\mu})
\vp^{\ga_\mu}
(2^{-l_\mu}z)dz )\times \\
& &\times 
H^\al_{j_1,l_1+k_1;\ldots; j_\al, l_\al+k_\al} f(x),
\end{eqnarray*}
where 
\begin{eqnarray*}
& & H^\al_{j_1,l_1+k_1;\ldots; j_\al, l_\al+k_\al} f(x)  = 
\intl   
 \overline{\si^{\al}(s,x(s)+U y,\xi)}
 e^{-4\pi^2 i (t-s)|\xi|^2} e^{2\pi i \dpr{\xi}{x(t)-x(s)+Ux-Uy}}
 \times \\ 
& & \times  \prod_{\mu=1}^\al \psi_{j_\mu,l_\mu+k_\mu}^{\ga_\mu}
(2^{l_\mu+k_\mu}(\xi/|\xi|-\theta_{j_\mu}^{l_\mu+k_\mu}))
\Om^2(\xi) f(y) dy d\xi.
\end{eqnarray*}
Hence by Lemma \ref{le:pointwise}, we conclude
\begin{eqnarray*}
& & 
\norm{W^{\al}(t)W^{\al}(s)^*f}{L^\infty_{x_2, 
\ldots, x_n} L^2_{x_1}}\leq  \\
& & \leq \supl_x \prod_{\mu=1}^\al (\suml_{k_\mu\leq -2}
\suml_{l_\mu\leq -2k_\mu} 
\suml_{\ga_\mu} (\ga_\mu!)^{-1} 
\suml_{j_\mu}  
\int A^{\ga_\mu}_{k_\mu}(t, x(t)+Ux +z \theta^{l_\mu+k_\mu}_{j_\mu})
\vp^{\ga_\mu}
(2^{-l_\mu}z)dz ) \times \\
& & \times \supl_{j_1,l_1,k_1;\ldots; j_\al, l_\al,k_\al}
\norm{H^\al_{j_1,l_1+k_1;\ldots; j_\al, l_\al+k_\al} f}{L^\infty_{x_2, 
\ldots, x_n} L^2_{x_1}}.\\
& & 
\leq C_n^\al  (\norm{A}{Y_1})^\al 
\supl_{j_1,l_1,k_1;\ldots; j_\al, l_\al,k_\al}
\norm{H^\al_{j_1,l_1+k_1;\ldots; j_\al, l_\al+k_\al} f}{L^\infty_{x_2, 
\ldots, x_n} L^2_{x_1}}.
\end{eqnarray*}
Dualizing the needed estimate for 
$\supl_{j_1,l_1,k_1;\ldots; j_\al, l_\al,k_\al}
\norm{H^\al_{j_1,l_1+k_1;\ldots; j_\al, l_\al+k_\al} f}{L^\infty_{x_2, 
\ldots, x_n} L^2_{x_1}}$, reduces matters to showing that 
\begin{eqnarray*}
& & (H^\al_{j_1,l_1+k_1;\ldots; j_\al, l_\al+k_\al})^* f(x)= 
\intl   \si^{\al}(s,x(s)+U x,\xi)
 e^{4\pi^2 i (t-s)|\xi|^2} e^{2\pi i \dpr{\xi}{x(s)-x(t)+Ux-Uy}}
 \times \\ 
& & \times  \prod_{\mu=1}^\al \psi_{j_\mu,l_\mu+k_\mu}^{\ga_\mu}
(2^{l_\mu+k_\mu}(\xi/|\xi|-\theta_{j_\mu}^{l_\mu+k_\mu}))
\Om^2(\xi) f(y) dy d\xi
\end{eqnarray*}
is a mapping $L^1_{x_2, 
\ldots, x_n} L^2_{x_1}\to L^\infty_{x_2, 
\ldots, x_n} L^2_{x_1}$ with norm no bigger than $C_n^\al \norm{A}{Y_1}^{\al}
|t-s|^{-(n-1)/2}$. 

Again, as in Section \ref{sec:6}, we expand 
$$
\si^\al(s, x(s)+Ux,\xi)=\prod_{\nu=1}^\al (\suml_{k_\nu\leq -2, 
l_\nu\leq -2k_\nu} 
\sum_{\ga_\nu, j_\nu} (\ga_\nu!)^{-1}  
\int A^{\ga_\nu}_{k_\nu}(t, x(s)+Ux +z \theta^{l_\nu+k_\nu}_{j_\nu})
\vp^{\ga_\nu}
(2^{-l_\nu}z)dz )
$$ and 
estimate by Lemma \ref{le:pointwise}. We get 
\begin{eqnarray*}
& &
 \norm{(H^\al_{j_1,l_1+k_1;\ldots; j_\al, l_\al+k_\al})^* f}{L^\infty_{x_2, 
\ldots, x_n} L^2_{x_1}}\leq C_n^\al\norm{A}{Y_1}^\al \times \\
& &\times \intl   
 e^{4\pi^2 i (t-s)|\xi|^2} e^{2\pi i \dpr{\xi}{x(s)-x(t)+Ux-Uy}}
  \prod_{\mu=1}^\al \psi_{j_\mu,l_\mu+k_\mu}^{\ga_\mu}
(2^{l_\mu+k_\mu}(\xi/|\xi|-\theta_{j_\mu}^{l_\mu+k_\mu})) \times \\
& &\times  \prod_{\nu=1}^\al \psi_{j_\nu,l_\nu+k_\nu}^{\ga_\nu}
(2^{l_\nu+k_\nu}(\xi/|\xi|-\theta_{j_\nu}^{l_\nu+k_\nu}))
\Om^2(\xi) f(y) dy d\xi
\end{eqnarray*}
By the Krein-Milman theorem, the linear span of 
$\{\de_b(x_2, \ldots, x_n) g(x_1): g\in L^2(\rone)\}$ is $w^*$ dense in 
$L^1_{x_2, \ldots, x_n} L^2_{x_1}$. Thus, 
it will suffice to verify the estimate 
\begin{eqnarray*}
& &\supl_{b,\bar{x} \in R^{n-1}, z\in \rn, \tilde{\theta}_1, \ldots, 
\tilde{\theta}_{2\al}\in\sn}
\|\intl   
 e^{4\pi^2 i (t-s)|\xi|^2} e^{2\pi i \dpr{\xi}{z+U(x_1, \bar{x})-U(y_1, b)}}\times \\
& & 
  \prod_{\mu=1}^{2\al} 
\psi_{\mu}(2^{m_\mu}(\xi/|\xi|-\tilde{\theta}_\mu))
\Om^2(\xi) g(y_1) dy_1 d\xi\|_{L^2_{x_1}} \leq C_n^\al |t-s|^{-(n-1)/2} 
\norm{g}{L^2(\rone)}. 
\end{eqnarray*}
Clearly, by rotational invariance, we can assume $U=Id$. 
But then the expression above is equal to 
$$
\intl   
 e^{4\pi^2 i (t-s)|\xi_1|^2} e^{2\pi i (\xi_1)(z_1+x_1-y_1)} K(\xi_1) 
g(y_1) dy_1= e^{-i(t-s)\p_1^2}K(\p_1/2i \pi)[g(\cdot)](z_1+x_1), 
$$
where 
$$
 K(\xi_1) = \int \intl   
 e^{4\pi^2 i (t-s)|\bar{\xi}|^2} 
e^{2\pi i \dpr{\bar{\xi}}{\bar{z}+\bar{x}-b}}
  \prod_{\mu=1}^{2\al} 
\psi_{\mu}(2^{m_\mu}(\xi/|\xi|-\tilde{\theta}_\mu))
\Om^2(\xi) d\xi_2\ldots d \xi_n.
$$
By Lemma \ref{le:osc} (see its second statement), we have 
$$
\supl_{\xi_1} |K(\xi_1)|\leq C_n^{\al}|t-s|^{-(n-1)/2}. 
$$
We get 
$$
\supl_{b,\bar{x} \in R^{n-1}, z\in \rn, \tilde{\theta}_1, \ldots, 
\tilde{\theta}_{2\al}} \norm{e^{-i(t-s)\p_1^2}K(\p_1/2i \pi)[g(\cdot)]}{L^2}
\lesssim \norm{g}{l^2}\supl_{\xi_1} |K(\xi_1)|\leq C_n^{\al}|t-s|^{-(n-1)/2} 
\norm{g}{l^2},
$$
as required. 

\end{proof}
\section{Global regularity for  Schr\"odinger maps}
\label{sec:11}
In this section, we sketch the proof of Theorem \ref{theo:MSM}. 
As it was discussed earlier, it will suffice to show that $A$ 
stays small in the space of vector potentials 
 $Y$, given the {\it a priori} information 
that $u$ is small in (a portion of ) 
the solution space to be described below.  
Let $\dot{X}^\al$ be the completion of all Schwartz functions in the norm 
$$
\|u\|_{\dot{X}^\al} = (\suml_k 2^{2 \al k}\supl_{q,r- Str.} 
\|u_k\|_{L^q L^r}^2 +2^{2 \al k} \supl_{U\in SU(n), x(t)} 
\|u_k(t, x(t)+Ux)\|_{L^2_t L^{2(n-1)/(n-3)}_{x_2, 
\ldots, x_n}L^2_{x_1}}^2 )^{1/2}
$$
Let $X^s:=\dot{X}^{s}\cap \dot{X}^{0}$. 
We will generally 
measure the solution in $X^s$, but 
moreover, we will show it is small in $\dot{X}^{s_0}$. Note that 
since $s>s_0$, 
$\norm{u}{\dot{X}^{s_0}}\lesssim \norm{u}{X^s}$. 

The space $Y$ of acceptable vector potentials 
is on the level of smoothness of $\dot{X}^{s_0}$.\\
 Fix $\de>0$, so that 
$s>(n+1)/2+\de$. We will only assume that $\norm{f}{\dot{H}^{(n+1)/2+\de}}
\norm{f}{\dot{H}^{(n-5)/2-\de}}<<1$. Clearly\footnote{Note that in the fromulation of the theorem, we have asked for a lot more, namely  
$\norm{f}{\dot{H}^s}=\ve\norm{g}{\dot{H}^s}<<1$ for all $s\in[0,(n+1)/2+\de)$.} 
$\norm{f}{\dot{H}^{s_0}}\leq (\norm{f}{\dot{H}^{(n+1)/2+\de}}\norm{f}{\dot{H}^{(n-5)/2-\de}})^{1/2}<<1$.  \\
 Under this assumptions, it will suffice to check 
\begin{itemize}
\item ($A=A(u)$ is controlled by  
$\norm{u}{\dot{X}^{s_0}}$ and $\norm{u}{X^s}$) 
\begin{eqnarray}
\label{eqnm:3}
& & \norm{A(u)}{Y}\leq C_n 
\norm{u}{\dot{X}^{(n+1)/2+\de}} \norm{u}{\dot{X}^{(n-5)/2-\de}}, \\
\label{eqnm:4}
& & (\suml_k 2^{2ks} \norm{A_k}{L^2_t L^{2n/(n-2)}}^2)^{1/2}\leq 
C_n \norm{u}{\dot{X}^s}\norm{u}{\dot{X}^{s_0}}\q\q \textup{for every} \q s\geq 0. 
\end{eqnarray}
\item ($N(u)$ is controlled by $\norm{u}{\dot{X}^{s_0}}$ and $\norm{u}{X^s}$)
\begin{eqnarray}
\label{eqnm:6}
& & \norm{\p^{s} N(u)}{L^1 L^2} \leq C_n 
\norm{u}{\dot{X}^{s}}(\norm{u}{\dot{X}^{s_0}}^2+\norm{u}{\dot{X}^{s_0}}^4)
\end{eqnarray}
\end{itemize}
Le ts us first show how  Theorem \ref{theo:MSM} follows from \eqref{eqnm:3}, 
\eqref{eqnm:4}, \eqref{eqnm:6}. 

To that end, we know that the Strichartz estimates  hold for at 
least for some time $T$, so that $\norm{A}{Y_T}\leq \ve$. Fix one such $T$. 
We have for every $s\geq 0$, by  \eqref{eqnm:4} and \eqref{eqnm:6}
\begin{eqnarray*}
& & \norm{u}{\dot{X}_T^{s}}\leq 
C_n (\norm{f}{\dot{H}^{s}}+ 
\norm{\p^s N(u)}{L^1_TL^2})+  C_n \norm{\nabla u}{L^2_T L^n_x}
(\suml_k 2^{2ks}\norm{A_k}{L^2_T L^{2n/(n-2)}}^2)^{1/2}\leq \\
& &\leq C_n \norm{f}{\dot{H}^{s}}+C_n \norm{u}{\dot{X}_T^s} 
(\norm{u}{\dot{X}_T^{s_0}}^2+\norm{u}{\dot{X}_T^{s_0}}^4)+
C_n \norm{u}{\dot{X}_T^{s_0}}
(\suml_k 2^{2ks}\norm{A_k}{L^2_T L^{2n/(n-2)}}^2)^{1/2}\leq \\
& & \leq C_n (\norm{f}{\dot{H}^{s}}+ 
\norm{u}{\dot{X}_T^s} 
(\norm{u}{\dot{X}_T^{s_0}}^2+
\norm{u}{\dot{X}_T^{s_0}}^4)).
\end{eqnarray*}
In particular, for $s=s_0$ and the smallness of 
$\norm{f}{\dot{H}^{s_0}}$, it follows that 
$$
\norm{u}{\dot{X}_T^{s_0}}\leq C_n \norm{f}{\dot{H}^{s_0}}.
$$
This means that $\norm{u}{\dot{X}_T^{s_0}}$ is 
small (independently of $T$), which  in 
turn implies that 
$$
\norm{u}{\dot{X}_T^{s}}\leq C_n \norm{f}{\dot{H}^{s}},
$$
for every $s\geq 0$.  
But how far can we really push that? Recall \eqref{eqnm:3}, which 
gives us a control of $A$ back in terms of $u$. Namely, since 
$$
 \norm{A(u)}{Y_T}\leq C_n 
\norm{u}{\dot{X}_T^{(n+1)/2+\de}} \norm{u}{\dot{X}_T^{(n-5)/2-\de}}\leq C_n
\norm{f}{\dot{H}^{(n+1)/2+\de}} \norm{f}{\dot{H}^{(n-5)/2-\de}}<<1.
$$
This implies  that $\norm{A}{Y_T}$ is small and one could 
apply back the Strichartz estimates, which means 
that $T$ could be taken to be $\infty$. 
Theorem \ref{theo:MSM} follows.

\subsection{Proof of \eqref{eqnm:3}, \eqref{eqnm:4}}
We will not give the full details of 
 \eqref{eqnm:3}, \eqref{eqnm:4} , since 
these are standard Besov type estimates for products. 

Let us for example consider the estimate for  $\norm{A}{Y_3}$. 
First, it is not hard to see that the terms containing $\p_t A$ one  uses the equation 
\eqref{eqnm:1}, to write it 
like 
$$
\p_t A = \p_t Q_1(u,\bar{u})\sim  \tilde{Q}_1 (u,u_t)\sim 
 \tilde{Q}_1 (u,\De u)+\tilde{Q}_1 (u,N(u)).
$$
Thus, everything is reduced to the terms containing 
$\p^2 A$ and $N(u)$, the latter being easy to treat. 

So, we concentrate on the terms involving $\p^2 A$. 
For those, 
take into account $\p^2A_k \sim 2^{2k}A_k$ and $A\sim \p^{-1}Q(u, \bar{u})$, to conclude 
\begin{eqnarray*}
& &\norm{A}{Y_3}\sim \suml_k  2^{k(n+3)/2} \supl_{U\in SU(n), x(t)}
\norm{A_k(t, x(t)+Ux)}{L^1_t L^2_{x_2, \ldots, x_n}L^1_{x_1}}\sim \\
& &\sim 
\suml_k  2^{k(n+1)/2} 
\supl_{U\in SU(n), x(t)}
\norm{(u v)_k (t, x(t)+Ux)}{L^1_t L^2_{x_2, \ldots, x_n}L^1_{x_1}}.
\end{eqnarray*}
Following Lemma 3.1 and Lemma 3.2 in \cite{rodnianski}, 
we have to split 
$(u v)_k$
in\footnote{Here $v$ might be either $u$ or $\bar{u}$.} 
two types of terms  - high low interactions $u_{\sim k} v_{<k+5}$ 
and 
high-high interactions 
$\suml_{l>k-5}(u_{l} v_{l-2<\cdot< l+2})_k$. 

The high-low interactions are more difficult to 
handle in this context\footnote{For the high-high interactions, one can actually split the $(n+1)/2$ derivatives between the two terms and get a better more 
balanced estimate.}, so let us concentrate on these. 
We have by Cauchy-Schwartz and 
Bernstein inequalities 
\begin{eqnarray*}
& & 2^{k(n+1)/2}\supl_{U\in SU(n), x(t)}
\norm{u_{k} v_{ k-m} 
(t, x(t)+Ux)}{L^1_t L^2_{x_2, \ldots, x_n}L^1_{x_1}}\leq \\
& & \leq 2^{k(n+1)/2} 
\supl_{U, x(t)}
\norm{u_{ k} (t, x(t)+Ux)}{L^2_t L^{2(n-1)/(n-3)}_{\bar{x}}L^2_{x_1}}
\norm{v_{k-m } (t, x(t)+Ux)}{L^2_t L^{n-1}_{\bar{x}}L^2_{x_1}}\\
& &\leq C_n 2^{-\de m}
\supl_{U, x(t)} 2^{k((n+1)/2+\de)}
\norm{u_{k} (t, x(t)+Ux)}{L^2_t L^{2(n-1)/(n-3)}_{\bar{x}}L^2_{x_1}}\times \\
& & 
\supl_{U, x(t)} 2^{(k-m)((n-5)/2-\de)} 
\norm{v_{k-m } (t, x(t)+Ux)}{L^2_t L^{2(n-1)/(n-3)}_{\bar{x}}L^2_{x_1}}
\end{eqnarray*}
Note that at this stage to apply the Bernstein inequality in the variables 
$x_2, \ldots, x_n$ \\ 
$\|v_{k-m}\|_{L^{n-1}_{x_2, \ldots, x_n}L^2_{x_1}}\leq C_n 
2^{(k-m)(n-5)/2}\|v_{k-m}\|_{L^{2(n-1)/(n-3)}_{x_2, \ldots, x_n} L^2_{x_1}}$, 
we needed $n-1>2(n-1)/(n-3)$, which is the dimensional restriction $n>5$. \\
Summing in $m>-5$ yields  
\begin{eqnarray*}
& & \suml_k 2^{k(n+1)/2}\supl_{U\in SU(n), x(t)}
\norm{u_{ k} v_{k-m}(t, x(t)+Ux)}{L^1_t 
L^2_{x_2, \ldots, x_n} L^1_{x_1}}\leq \\
& & \leq 
(\suml_k 2^{k((n+1)+2\de)} \supl_{U, x(t)} 
\norm{u_k(t, x(t)+Ux)}{L^2_t L^{2(n-1)/(n-3)}_{\bar{x}}
L^2_{x_1}}^2)^{1/2}\times \\
& & \times 
(\suml_k  2^{k((n-5)-2\de)}  \supl_{U, x(t)}
\norm{v_k (t, x(t)+Ux)}{L^2_t 
L^{2(n-1)/(n-3)}_{\bar{x}}L^2_{x_1}}^2)^{1/2}\leq \\
& & \leq C_n 
\norm{u}{\dot{X}^{(n+1)/2+\de}}\norm{v}{\dot{X}^{(n-5)/2-\de}}.
\end{eqnarray*}

The proof of \eqref{eqnm:4} is in fact very similar 
 and boils down to the same Besov space estimates for products. 
\subsection{Proof of \eqref{eqnm:6}} 
The estimates for the nonlinearities are the easiest ones. 
It basically suffice to apply the Kato-Ponce type-estimates
$$
\norm{\p^s(u v)}{L^r}\leq C\norm{\p^s u}{L^p}\norm{v}{L^q}+
C\norm{\p^s v}{L^p}\norm{u}{L^q}
$$
whenever $1/r=1/p+1/q$. We omit the details. 

\section{Appendix}
\subsection{Decay estimates for the free Schr\"odinger equation 
with initial data Fourier supported in a small cap} 
\begin{lemma}
\label{le:osc}
Let $k_1, \ldots, k_\mu$ be positive integers. 
Let also $\{\theta_j\}\in \sn$ and 
$\psi_j$ be smooth cutoff  functions adapted to 
$\{|\xi|\leq 1\}$, whose smoothness bounds are uniform in $j$. 
 Then there exists a constant $C$ 
depending only on the dimension, so that 
\begin{equation}
\label{eq:514}
\supl_{\theta_1, k_1, \ldots, \theta_\mu\in\sn, k_\mu}\supl_x | \intl   
 e^{-4\pi^2 i t |\xi|^2} e^{2\pi i \dpr{\xi}{x}}
\prod_{j=1}^\mu \psi_{j}(2^{k_j}(\xi/|\xi|-\theta_{j})) 
\vp(\xi) d\xi|\leq  C_n^\mu |t|^{-n/2}.
\end{equation}
Also, if one fixes $\xi_1$ and integrates with respect 
to $\xi_2, \ldots, \xi_n$,  
$$
\supl_{\theta_1, k_1, \ldots, \theta_\mu \in \sn, k_\mu}
\supl_{x,\xi_1} | \intl   
 e^{-4\pi^2 i t |\bar{\xi}|^2} e^{2\pi i \dpr{\bar{\xi}}{\bar{x}}}
\prod_{j=1}^\mu \psi_{j}(2^{k_j}(\xi/|\xi|-\theta_{j})) 
\vp(\xi) d\xi_2\ldots d\xi_n|\leq  C_n^\mu |t|^{-(n-1)/2}.
$$
\end{lemma}
\begin{proof}
We prove only \eqref{eq:514}, since the second statement in Lemma 
\ref{le:osc} requires  only a slight adjustment of the argument. 

This a standard stationary phase argument, 
except that we have to keep track of the 
derivatives that may pile up from the cutoffs 
$\psi_{j}(2^{k_j}(\xi/|\xi|-\theta_{j}))$. Fix $k_1, \ldots, k_\mu$ and let 
$k=\max(k_1, \ldots, k_\mu)=k_{j_0}$. 
If $2^k\geq \sqrt{t}$, we pass to polar coordinates 
and estimate by 
\begin{eqnarray*}
& &C \intl_{|\theta-\theta_{j_0}|\lesssim  2^{-k}} |\int 
e^{-4\pi^2 i t \rho^2} e^{2\pi \rho i \dpr{\theta}{x}}
\vp(\rho)d\rho|d\theta \leq 
C 2^{-k(n-1)} |t|^{-1/2}\leq C |t|^{-n/2},
\end{eqnarray*}
where we have used the decay estimate by $C t^{-1/2}$ for the 1 D Schr\"odinger 
equation. 

Thus assume $2^k\leq \sqrt{t}$.
We (smoothly) split the region of integration, according to the
 size of the derivative of the phase. If 
 $|-8\pi^2 t \xi+2\pi x|\leq \sqrt{t}$, 
we estimate by absolute values and obtain the desired estimate 
by the volume of the $\xi$ support, which is $C_n |t|^{-n/2}$.

 It remains to show that 
\begin{eqnarray*}
& &| \intl   
 e^{-4\pi^2 i t |\xi|^2} e^{2\pi i \dpr{\xi}{x}} 
\vp(2^{-m} t^{-1/2} (-8\pi^2 t \xi+2\pi x))
\prod_{j=1}^\mu \psi_{j}(2^{k_j}(\xi/|\xi|-\theta_{j})) 
\vp(\xi) d\xi|\\
& & \leq C 2^{-m} C_{n}^\mu |t|^{-n/2},
\end{eqnarray*}
Now that we have $2^{k_j}\leq 2^{k}\leq \sqrt{t}$. Then, the  
argument goes as in the classical estimate, that is after $n+1$ 
integration by parts with the phase function 
$\varrho(\xi)=i(-4\pi^2  t |\xi|^2+2\pi  \dpr{\xi}{x})$, 
(at  each step one gains at least a  
factor of $2^{-m}|t|^{-1/2}$ 
and  loses a factor of $C_{n} \sqrt{|t|}$ at the most), 
we put in absolute values. Taking 
into account the volume of the support 
$\leq C_n 2^{mn} t^{-n/2}$, we estimate by 
$$
 D_{n}2^\mu  2^{-m} |t|^{-n/2}\leq C_n^\mu 2^{-m} |t|^{-n/2},
$$
whence summation by $m\geq 0$ yields the result. 
\end{proof}
\subsection{Estimates on the error term}
In this section we give an estimate on the error terms $E^k$, defined in 
\eqref{eq:1201}. 
\begin{lemma}
Let  $\norm{A}{Y_0}\leq \ve$ and 
$s\geq 0$. Let also $1=1/p_1+1/p_2; 1/2=1/q_1+1/q_2$. Then  
\begin{eqnarray}
\label{eq:error}
\suml_k 2^{2ks} \norm{E^k}{L^1 L^2}^2\lesssim \ve^2
\suml_k 2^{2ks} \norm{u_k}{L^\infty L^2}^2+
\norm{u}{L^{p_1}L^{q_1}}^2 
\suml_k 2^{2k(s+1)} \norm{A_k}{L^{p_2}_t L^{q_2}_x}^2,\\
\label{e:11}
\suml_k 2^{2ks} \norm{E^k}{L^1 L^2}^2\lesssim \ve^2
\suml_k 2^{2ks} \norm{u_k}{L^\infty L^2}^2+
\norm{\nabla u}{L^{p_1}L^{q_1}}^2 
\suml_k 2^{2k s} \norm{A_k}{L^{p_2}_t L^{q_2}_x}^2.
\end{eqnarray}
In particular, from \eqref{eq:error}, with 
$s=0, q_1=2, p_1=\infty, p_2=1, q_2=\infty$, 
and since $\norm{A_k}{L^1L^\infty_x}
\leq 2^{kh} \norm{A_k}{L^1 L^{n/h}}$
\begin{equation}
\label{e:10}
\suml_k \norm{E^k}{L^1 L^2}^2\lesssim \ve^2
\suml_k  \norm{u_k}{L^\infty L^2}^2
\end{equation}
\end{lemma}
\begin{proof}
We estimate on a term-by-term basis in formula \eqref{eq:1201}. 

For the first term, we have 
\begin{eqnarray*}
\suml_k  2^{2ks} \norm{[P_k,\vec{A}_{\leq k-4}] 
\nabla u_{k-1\leq \cdot\leq k+1}}{L^1L^2}^2
\lesssim 
\norm{\nabla A}{L^1L^\infty}^2 \suml_k 2^{2ks} \norm{u_k}{L^\infty L^2}^2.
\end{eqnarray*}

The second and third terms in \eqref{eq:1201} are treated in a similar 
fashion, so we 
concentrate on the second one. For any positive $h\leq 1$, we have
\begin{eqnarray*}
& & \suml_k 2^{2ks} \left(\suml_{l\geq k-2} \norm{P_k( A_l \nabla 
u_{l-2\leq \cdot \leq l+2})}{L^1L^2}\right)^2\lesssim\\
& &\lesssim \suml_k 2^{2ks} 
\left(\suml_{l\geq k-2} 2^{h k} \norm{P_k( A_l \nabla 
u_{l-2\leq \cdot \leq l+2})}{L^1L^{2n/(n+h)}}\right)^2\lesssim \\
& &\lesssim\suml_k 2^{2 k(s+h)} \left(\suml_{l\geq k-2} 2^l 
\norm{A_l}{L^1 L^{2 n/h}}
\norm{u_{l-2\leq \cdot \leq l+2}}{L^\infty L^{2}}\right)^2\lesssim \\
& & \lesssim\suml_k 2^{2 k(s+h)} \left(\suml_{l\geq k-2} 2^{-h l} 
(2^{l+lh}\norm{A_l}{L^1 L^{n/h}})
\norm{u_l}{L^2L^{2n/(n-2)}}\right)^2
\end{eqnarray*}
One obtains by
Lemma \ref{le:sequences} an estimate by 
\begin{eqnarray*}
\suml_l 2^{2l(1+h)}\norm{A_l}{L^1 L^{n/h}}^2
\suml_l 2^{2ls} \norm{u_l}{L^\infty L^2}^2.
\end{eqnarray*}
The last fourth term in \eqref{eq:1201} can be bounded in two ways. 
\begin{eqnarray*}
& & \suml_k 2^{2ks} \norm{A_{k-1\leq \cdot\leq k+1}\nabla 
u_{\leq k-4})}{L^1 L^2}^2\lesssim \norm{u}{L^{p_1} L^{q_1}}^2
\suml_k 2^{2k (s+1)} \norm{A_k}{L^{p_2} L^{q_2}}^2,\\
& & \suml_k 2^{2ks} \norm{A_{k-1\leq \cdot\leq k+1}\nabla 
u_{\leq k-4})}{L^1 L^2}^2\lesssim \norm{\nabla u}{L^{p_1} L^{q_1}}^2
\suml_k 2^{2k s} \norm{A_k}{L^{p_2} L^{q_2}}^2
\end{eqnarray*}

\end{proof}

\end{document}